\documentclass[a4paper,review]{cas-sc}

\usepackage[section]{placeins} 
\usepackage{amssymb}
\usepackage{amsthm}
\usepackage{amsmath}
\usepackage{mathtools}
\usepackage{upgreek}
\usepackage{pdflscape}
\usepackage{listings}
\usepackage{multirow}
\usepackage[percent]{overpic}
\usepackage{multicol}
\usepackage{color}
\usepackage{stmaryrd}
\usepackage{xfrac}
\usepackage{enumerate}
\usepackage[capitalise]{cleveref}
\usepackage{booktabs}
\usepackage[section]{placeins} 
\usepackage[section]{algorithm}
\usepackage{calc}
\usepackage[titletoc]{appendix}
\usepackage{siunitx}
\usepackage[sort&compress,square,numbers]{natbib}

\usepackage{graphicx}
\usepackage{graphics}
\usepackage{wrapfig}
\usepackage{float}
\usepackage{subfig}
\usepackage[percent]{overpic}
\usepackage{varwidth}
\usepackage{tikz}
\usepackage{pgfplots}
\usepackage{adjustbox}
\usetikzlibrary{arrows,matrix,positioning,fit}
\usetikzlibrary{shapes,positioning}
\usetikzlibrary{backgrounds}
\usepackage{tikz-layers}
\usepgfplotslibrary{fillbetween}
\usetikzlibrary{intersections}
\pgfplotsset{compat=1.14}
\usepgfplotslibrary{colorbrewer}
\usepgfplotslibrary{patchplots}
\usepgfplotslibrary[colorbrewer]
\usetikzlibrary{pgfplots.colorbrewer}
\usetikzlibrary[pgfplots.colorbrewer]
\usepgfplotslibrary{units}
\usetikzlibrary{spy}
\usepackage{pgfplotstable}
\usepackage{arrayjobx}
\graphicspath{ {./figs/} }
\usetikzlibrary{external}

\newlength\myheight
\newlength\mydepth
\settototalheight\myheight{Xygp}
\newcommand*\inlinegraphics[1]{%
  \settototalheight\myheight{Xygp}%
  \settodepth\mydepth{Xygp}%
  \raisebox{-\mydepth}{\includegraphics[height=\myheight]{#1}}%
}
\newcommand\orcid[1]{\href{https://orcid.org/#1}{\inlinegraphics{orcid_16x16.png}}}

\makeatletter
\def\BState{\State\hskip-\ALG@thistlm}
\makeatother

\newtheorem{theorem}{Theorem}[section]
\newtheorem{corollary}{Corollary}[theorem]

\newdefinition{definition}{Definition}[section]
\newtheorem*{remark}{Remark}

\newcommand\pxvar[2]{\partial_{#2} #1}

\begin{document}

\title[mode=title]{High-order limiting methods using maximum principle bounds derived from the Boltzmann equation I: Euler equations}
\shorttitle{High-order limiting methods using maximum principle bounds derived from the Boltzmann equation I: Euler equations}
\shortauthors{T. Dzanic and L. Martinelli}

\author[1]{Tarik Dzanic}[orcid=0000-0003-3791-1134]
\cormark[1]
\cortext[cor1]{Corresponding author}
\ead{dzanic1@llnl.gov}
\author[2]{Luigi Martinelli}[orcid=0000-0001-6319-7646]

\address[1]{Center for Applied Scientific Computing, Lawrence Livermore National Lab, Livermore, CA 94550, USA}
\address[2]{Department of Mechanical and Aerospace Engineering, Princeton University, Princeton, NJ 08544, USA}

\begin{abstract}
The use of limiting methods for high-order numerical approximations of hyperbolic conservation laws generally requires defining an admissible region/bounds for the solution. In this work, we present a novel approach for computing solution bounds and limiting for the Euler equations through the kinetic representation provided by the Boltzmann equation, which allows for extending limiters designed for linear advection directly to the Euler equations. Given an arbitrary set of solution values to compute bounds over (e.g., numerical stencil) and a desired linear advection limiter, the proposed approach yields an analytic expression for the admissible region of particle distribution function values, which may be numerically integrated to yield a set of bounds for the density, momentum, and total energy. These solution bounds are shown to preserve positivity of density/pressure/internal energy and, when paired with a limiting technique, can robustly resolve strong discontinuities while recovering high-order accuracy in smooth regions without any ad hoc corrections (e.g., relaxing the bounds). This approach is demonstrated in the context of an explicit unstructured high-order discontinuous Galerkin/flux reconstruction scheme for a variety of difficult problems in gas dynamics, including cases with extreme shocks and shock-vortex interactions. Furthermore, this work presents a foundation for limiting techniques for more complex macroscopic governing equations that can be derived from an underlying kinetic representation for which admissible solution bounds are not well-understood. 
\end{abstract}




\maketitle

\section{Introduction}
\label{sec:intro}
The understanding of the dynamics of inviscid fluid flow governed by the compressible Euler equations plays a critical role in predicting the behavior of fluids in high-speed regimes. These gas dynamics equations are not only fundamental as a governing model for the study of compressible flow phenomena such as shock waves but also serve as a stepping-stone to more complex simulations of viscous compressible flows. Numerical methods for simulating transport-dominated physics have widely relied on the Euler equations as the prototypical model problem, largely due to the challenges posed by the strongly nonlinear interactions between shock waves, discontinuities, and vortical features in the flow. These challenges have driven the development of sophisticated numerical techniques to yield robust and high-fidelity resolution of these complex flow features. 

Various numerical approaches have been utilized for gas dynamics simulations, with the brunt of applications focusing on second-order finite volume and finite difference schemes. Although well-established in the literature and broadly robust for many practical uses, their relatively high numerical dissipation and low accuracy make them prohibitively expensive for high-fidelity (i.e., scale-resolving) simulations of complex unsteady flows. As a result, a particular class of numerical schemes which have more recently grown in popularity are high-order methods, which broadly encompass finite difference, finite volume, and finite element methods that can achieve arbitrarily high-order accuracy. Due to this increased accuracy paired with lower numerical dissipation, these methods can generally offer better resolution of small-scale flow features and transport-dominated flow physics at a lower relative computational cost than low-order schemes~\citep{Vermeire2017}. However, the high-order approximations which underpin these methods are primarily designed to yield high accuracy in smooth regions of the flow. Around discontinuities or severe gradients, the use of high-order schemes can cause spurious oscillations in the numerical solution, often referred to as Gibbs phenomena. These spurious oscillations result in a numerically ill-behaved solution and, in many cases, the failure of the numerical scheme altogether.

To extend the use of high-order method to problems that exhibit discontinuities and steep gradients, it is often necessary to augment the scheme with an additional numerical stabilization approach to increase its robustness in the vicinity of such features. One such stabilization approach is through \textit{a priori} or \textit{a posteriori} limiting on the discrete solution, which relies on using a secondary numerical approximation that is more robust (but generally less accurate). The solutions as predicted by the high-order scheme and the secondary low-order scheme are then combined in a manner which ideally ensures sufficient numerical stabilization around discontinuities without sacrificing the accuracy of the high-order scheme in smooth regions of the flow. Various forms of limiting have been applied to hyperbolic conservation laws, typically falling into variants of the flux-corrected transport method of \citet{Boris1997} (e.g.,  \citep{Zalesak1979,Pazner2021,Lin2023,RuedaRamrez2022,Peyvan2023,Kuzmin2020}) or through limiting/modification of the underlying solution/operator approximations (e.g., \citep{Harten1983,Barth1989,Jameson1995,Shu1998,Zhang2011b,Zhu2013,Zhang2015,Dumbser2016,Dzanic2022,Ching2023,Pandare2023,Zhang2023}). 
We remark here that these examples only cover a small subset of the numerical approaches in the literature and do not by any means encompass all the possible limiting techniques used over the decades.

One of the primary difficulties in designing performant limiting methods for high-order schemes stems not necessarily from the development of the limiting schemes themselves but instead from the difficulty of defining exactly what a ``well-behaved'' numerical solution is for nonlinear systems of equations. Limiting between the low-order and high-order solution is typically dictated by some ad hoc heuristic (e.g., indicator functions) or through a set of bounds (i.e., admissible regions/sets) for the solution. While heuristic approaches can be effective for some applications, they often rely on empirical problem-dependent parameters which can misidentify regions that require limiting and degrade the accuracy and robustness of the approach. On the other hand, limiting based on a set of admissible bounds for the solution is typically a highly-robust alternative, often providing guarantees that the limited solution will abide by certain physical constraints~\citep{Guermond2018,Kuzmin2020} (although typically at a larger computational expense). However, it is unclear what the optimal bounds are to enforce on the solution, tying in to the problem of defining what a well-behaved numerical solution is. Bounds on the numerical properties of the solution (e.g., total variation diminishing~\citep{Harten1983}, local extremum diminishing methods~\citep{Jameson1995}, etc.) have shown to be robust for many applications, but these bounds are often quite restrictive and unnecessarily degrade the accuracy of the high-order scheme in smooth regions without ad hoc corrections, whereas more advanced stencil modification techniques (e.g., ENO/WENO/TENO-type schemes and limiters~\citep{Shu1998,Fu2016}) can retain this accuracy but are typically prohibitively expensive and algorithmically complex for arbitrary mixed unstructured meshes. Bounds based on attributes of the underlying governing equations/systems have shown promise, relying on either physical constraints (e.g., positivity of density/pressure~\citep{Zhang2011b,Kuzmin2020}) or additional mathematical properties (e.g., minimum entropy principle~\citep{Dzanic2022,Ching2023}, invariant region~\citep{Guermond2016,Pazner2021}). Nevertheless, they still exhibit some difficulties, such as sometimes also reducing to first-order accuracy around smooth extrema (requiring ad hoc corrections such as relaxing the bounds), requiring nonlinear limiting procedures which may be computationally expensive, and, most importantly, relying on mathematical properties which may not easily extend to systems which are of more practical interest (e.g., Navier--Stokes equations, multi-physics applications, etc.).

The methods to be presented in this work are motivated by two primary observations. The first observation is that numerical approximations of the nonlinear Boltzmann equation, a general kinetic description of molecular gas dynamics which underpins the behavior of fluids governed by the Euler and Navier--Stokes equations, are significantly more robust than approximations of the associated continuum conservation laws. In fact, high-order schemes for the Boltzmann equation have been effectively used to directly resolve shock structures without any numerical shock capturing approaches (or with just simple limiting techniques) while accurately resolving smooth features in the flow~\citep{Wilde2021,Guo2015,Dzanic2023b, Dzanic2024}. The second observation is that linear advection, which governs particle transport in the Boltzmann equation, is a significantly easier governing equation to design high-order limiting schemes for. Stabilization methods for high-order approximations of the advection equation have been developed for decades, and due to the linear scalar nature of the equation, simple yet effective bounds can be formed for the solution (e.g., local maximum principle). As such, there exist a variety of high-order limiting schemes for the advection equation that one can choose from, with varying degrees of robustness and accuracy.

With these observations in mind, the goal of this work is to introduce a novel approach to computing and enforcing bounds for limiting schemes used in high-order approximations of the Euler equations. The basis of this work is the connection between the Euler equations and the underlying kinetic representation provided by the Boltzmann equation, which has used to derive a variety of numerical schemes for gas dynamics from kinetic theory, e.g.~\citep{Pullin1980,  Deshpande1986, Perthame1992, Prendergast1993, Xu1996}.
This connection opens up the possibility of extending limiting techniques designed for linear advection to the Euler equations. The proposed approach relies on the existence of an exact representation of the particle distribution function of the Boltzmann equation in the inviscid limit, for which the desired linear advection limiting technique can analytically be applied to yield a set of admissible bounds for the \emph{conserved variables} (i.e., density, momentum, energy) in the Euler equations. We demonstrate this approach in the context of an explicit unstructured high-order discontinuous Galerkin/flux reconstruction scheme~\citep{Hesthaven2008DG, Huynh2007} utilizing an associated local maximum principle preserving limiter~\citep{Zhang2010} for the underlying advection equation --- however, this approach is broadly applicable in the sense of computing bounds given a set of discrete solution values (e.g., numerical stencil, local neighborhood, etc.), which can be extended to other bounds-preserving limiting methods such as flux-corrected transport and/or slope-limiting techniques. This approach is used to simulate a variety of difficult problems in gas dynamics, including cases with extreme shocks, near-vacuum conditions, and small-scale smooth flow features, showcasing its robustness around strong discontinuities and accuracy in smooth regions of the flow. Furthermore, this work presents a foundation for limiting techniques for more complex macroscopic governing equations that can be derived from an underlying kinetic representation (e.g., Navier--Stokes equations, magnetohydrodynamics, radiation-hydrodynamics, etc.) for which admissible solution bounds are not necessarily well-understood.

The remainder of this manuscript is organized as follows. In \cref{sec:prelim}, some preliminaries on the Euler and Boltzmann equations as well as limiting methods for high-order discontinuous Galerkin schemes are presented. The proposed approach is then introduced in \cref{sec:methodology}, followed by numerical analysis in \cref{sec:analysis} and implementation details in \cref{sec:implementation}. The results of a variety of numerical experiments are then presented in \cref{sec:results} followed by conclusions in \cref{sec:conclusion}.

\section{Preliminaries}
\label{sec:prelim}

\subsection{Euler equations}
The system of interest in this work is represented by the compressible Euler equations for gas dynamics in $d$ dimensions without source terms, written in conservation form as
\begin{equation}\label{eq:gen_hype}
        \pxvar{\mathbf{w}(\mathbf{x}, t)}{t} + \boldsymbol{\nabla}\cdot\mathbf{F}(\mathbf{w}) = 0,
\end{equation}
where the solution $\mathbf{w}$ and flux $\mathbf{F}(\mathbf{w})$ are defined as 
\begin{equation}\label{eq:euler}
    \mathbf{w} = \begin{bmatrix}
            \rho \\ \mathbf{m} \\ E
        \end{bmatrix} \quad  \mathrm{and} \quad \mathbf{F}(\mathbf{w}) = \begin{bmatrix}
            \mathbf{m}^T\\
            \mathbf{m}\otimes\mathbf{U} + P\mathbf{I}\\
        (E+P)\mathbf{U}^T
    \end{bmatrix}.
\end{equation}
Here, $\rho$ is the density, $\mathbf{m}$ is the momentum, and $E$ is the total energy. The symbol $\mathbf{I}$ denotes the identity matrix in $\mathbb R^{d \times d}$, $\mathbf{U} = \mathbf{m}/\rho$ denotes the velocity, and $P$ denotes the pressure, computed through the ideal gas assumption as 
\begin{equation}
    P = (\gamma - 1)\rho e = (\gamma - 1)\left (E - \frac{1}{2} \mathbf{m}{\cdot}\mathbf{m}/ \rho \right),
\end{equation}
where $\gamma$ is the specific heat ratio and $e$ is the specific internal energy. It is sometimes convenient to represent the solution in terms of the primitive variables $\mathbf{q} = [\rho, \mathbf{U}, P]^T$. The solution of the Euler equations is endowed with a convex invariant set corresponding to the strict positivity of density ($\rho > 0$) and strict positivity of pressure/internal energy ($P, \rho e > 0$). 

\subsection{Boltzmann equation}
The Euler equations can be recovered as the limiting behavior of the Boltzmann equation, given as 
\begin{equation}\label{eq:boltzmann}
    \partial_t f (\mathbf{x}, \mathbf{u}, t) + \mathbf{u} {\cdot} \nabla f = \frac{1}{Kn}\mathcal C(f, f'),
\end{equation}
where $f (\mathbf{x}, \mathbf{u}, t)\in \mathbb R^d$ is a scalar particle distribution function, $Kn$ is the Knudsen number, $\mathcal C(f, f')$ is the collision operator which models particle interactions \citep{Cercignani1988}, and $f' (\mathbf{x}, \mathbf{u}, t)\in \mathbb R^d$ is a post-collision scalar particle distribution function. This distribution function is defined over a higher-dimensional space, which includes not only the spatial position $\mathbf{x}$ but also the microscopic velocity $\mathbf{u}$, and represents the probability density of a particle existing at some location $\mathbf{x}$ with some velocity $\mathbf{u}$. For brevity, we neglect the internal energy domain here, but the effects of internal energy are later presented. The conserved flow variables $\mathbf{w}(\mathbf{x}, t)$ can be recovered through the moments of the distribution function as 
\begin{equation}\label{eq:moments}
    \mathbf{w}(\mathbf{x}, t) = \left[\rho, \mathbf{m}, E \right]^T = 
    \int_{\mathbb R^d} f (\mathbf{x}, \mathbf{u}, t)\ \boldsymbol{\psi} (\mathbf{u}) \ \mathrm{d}\mathbf{u},
\end{equation}
where $\boldsymbol{\psi} (\mathbf{u}) \coloneqq [1, \mathbf{u}, (\mathbf{u}\cdot\mathbf{u})/2]^T$ is the vector of collision invariants. We utilize the notation that the lowercase symbol $\mathbf{u}$ refers to the microscopic velocity and the uppercase symbol $\mathbf{U}$ refers to the macroscopic velocity.

In the limit of zero Knudsen number (i.e., zero viscosity), the distribution function $f(\mathbf{x}, \mathbf{u}, t)$ tends towards a unique equilibrium distribution function $g(\mathbf{x}, \mathbf{u}, t)$ which minimizes (across the velocity domain) the entropy $H(z)$ for some argument $z$, i.e.,
\begin{equation}
    g = \underset{z}{\mathrm{arg\ min}}\ H(z), \quad \quad H(z) = \int_{\mathbb R^d} z \log(z) \ \mathrm{d} \mathbf{u}.
\end{equation}
and possesses identical primary moments, i.e.,
\begin{equation}\label{eq:equilibriummoments}
    \int_{\mathbb R^d} g (\mathbf{x}, \mathbf{u}, t)\ \boldsymbol{\psi} (\mathbf{u}) \ \mathrm{d}\mathbf{u} = \int_{\mathbb R^d} f (\mathbf{x}, \mathbf{u}, t)\ \boldsymbol{\psi} (\mathbf{u}) \ \mathrm{d}\mathbf{u}.
\end{equation}
From the H-theorem, the equilibrium distribution function (for a monatomic particle) is a Maxwell--Boltzmann distribution of the form
\begin{equation}\label{eq:maxwell}
    g(\mathbf{x}, \mathbf{u}, t) = \frac{\rho (\mathbf{x}, t)}{\left[2 \pi \theta(\mathbf{x}, t) \right]^{d/2}}\exp \left [-\frac{ \|\mathbf{u} - \mathbf{U}(\mathbf x, t) \|_2^2}{2 \theta (\mathbf x, t)} \right],
\end{equation}
where $\theta = P/\rho$ is a scaled temperature and $\|\cdot \|_2^2$ denotes the squared norm along the spatial dimension. We sometimes use the notation $g(\mathbf{w})$ to denote the equilibrium distribution corresponding to the macroscopic state $\mathbf{w}$. The effects of particle collision drive the distribution function $f$ to this equilibrium state $g$. In particular, the collision operator is known to conserve primary moments, i.e.,
\begin{equation} 
    \int_{\mathbb R^d}  C(f,f')\ \boldsymbol{\psi} (\mathbf{u})  \ \mathrm{d}\mathbf{u} = \boldsymbol{0},
\end{equation}
and is identically zero at equilibrium, i.e.,
\begin{equation}\label{eq:cgg}
    C(g,g') = 0.
\end{equation}

\subsection{Discontinuous Galerkin methods}

   \begin{figure}[tbhp]
        \centering
        \adjustbox{width=0.45\linewidth, valign=b}{     \begin{tikzpicture}[spy using outlines={rectangle, height=3cm,width=2.3cm, magnification=3, connect spies}]
		\begin{axis}[name=plot1,
		    axis line style={draw=none},
		    tick style={draw=none},
		    axis x line=left,
            axis y line=left,
            axis equal image,
            clip mode=individual,
    		xmin=-2.3,
    		xmax=2.3,
    		xticklabels={,,},
    		ymin=-1.1,
    		ymax=1.2,
    		yticklabels={,,},
    		style={font=\Large},
    		scale = 1]

    	\addplot[black,thick]
        coordinates {(-1.0, -1.0) (-0.895, -0.95) (-0.789, -0.906) (-0.684, -0.867) (-0.579, -0.834) (-0.474, -0.806) (-0.368, -0.784) (-0.263, -0.767) (-0.158, -0.756) (-0.053, -0.751) (0.053, -0.751) (0.158, -0.756) (0.263, -0.767) (0.368, -0.784) (0.474, -0.806) (0.579, -0.834) (0.684, -0.867) (0.789, -0.906) (0.895, -0.95) (1.0, -1.0) (1.0, -1.0) (0.947, -0.832) (0.895, -0.672) (0.842, -0.518) (0.789, -0.371) (0.737, -0.231) (0.684, -0.098) (0.632, 0.028) (0.579, 0.147) (0.526, 0.259) (0.474, 0.364) (0.421, 0.463) (0.368, 0.554) (0.316, 0.639) (0.263, 0.716) (0.211, 0.787) (0.158, 0.85) (0.105, 0.907) (0.053, 0.957) (0.0, 1.0) (0.0, 1.0) (-0.053, 0.847) (-0.105, 0.7) (-0.158, 0.558) (-0.211, 0.421) (-0.263, 0.289) (-0.316, 0.163) (-0.368, 0.042) (-0.421, -0.074) (-0.474, -0.185) (-0.526, -0.29) (-0.579, -0.39) (-0.632, -0.485) (-0.684, -0.574) (-0.737, -0.658) (-0.789, -0.737) (-0.842, -0.811) (-0.895, -0.879) (-0.947, -0.942) (-1.0, -1.0) };

    	\addplot[black,thick]
        coordinates {(1.2, -1.0) (1.147, -0.832) (1.095, -0.672) (1.042, -0.518) (0.989, -0.371) (0.937, -0.231) (0.884, -0.098) (0.832, 0.028) (0.779, 0.147) (0.726, 0.259) (0.674, 0.364) (0.621, 0.463) (0.568, 0.554) (0.516, 0.639) (0.463, 0.716) (0.411, 0.787) (0.358, 0.85) (0.305, 0.907) (0.253, 0.957) (0.2, 1.0) (0.2, 1.0) (0.305, 0.98) (0.411, 0.962) (0.516, 0.947) (0.621, 0.934) (0.726, 0.922) (0.832, 0.914) (0.937, 0.907) (1.042, 0.902) (1.147, 0.9) (1.253, 0.9) (1.358, 0.902) (1.463, 0.907) (1.568, 0.914) (1.674, 0.922) (1.779, 0.934) (1.884, 0.947) (1.989, 0.962) (2.095, 0.98) (2.2, 1.0) (2.2, 1.0) (2.147, 0.935) (2.095, 0.865) (2.042, 0.791) (1.989, 0.712) (1.937, 0.629) (1.884, 0.541) (1.832, 0.449) (1.779, 0.353) (1.726, 0.252) (1.674, 0.147) (1.621, 0.037) (1.568, -0.077) (1.516, -0.196) (1.463, -0.319) (1.411, -0.446) (1.358, -0.578) (1.305, -0.714) (1.253, -0.855) (1.2, -1.0) };

    	\addplot[black,thick]
        coordinates {(-0.2, 1.0) (-0.253, 0.847) (-0.305, 0.7) (-0.358, 0.558) (-0.411, 0.421) (-0.463, 0.289) (-0.516, 0.163) (-0.568, 0.042) (-0.621, -0.074) (-0.674, -0.185) (-0.726, -0.29) (-0.779, -0.39) (-0.832, -0.485) (-0.884, -0.574) (-0.937, -0.658) (-0.989, -0.737) (-1.042, -0.811) (-1.095, -0.879) (-1.147, -0.942) (-1.2, -1.0) (-1.2, -1.0) (-1.253, -0.915) (-1.305, -0.828) (-1.358, -0.738) (-1.411, -0.646) (-1.463, -0.552) (-1.516, -0.456) (-1.568, -0.357) (-1.621, -0.256) (-1.674, -0.153) (-1.726, -0.048) (-1.779, 0.059) (-1.832, 0.169) (-1.884, 0.281) (-1.937, 0.395) (-1.989, 0.512) (-2.042, 0.63) (-2.095, 0.751) (-2.147, 0.875) (-2.2, 1.0) (-2.2, 1.0) (-2.095, 1.02) (-1.989, 1.038) (-1.884, 1.053) (-1.779, 1.066) (-1.674, 1.078) (-1.568, 1.086) (-1.463, 1.093) (-1.358, 1.098) (-1.253, 1.1) (-1.147, 1.1) (-1.042, 1.098) (-0.937, 1.093) (-0.832, 1.086) (-0.726, 1.078) (-0.621, 1.066) (-0.516, 1.053) (-0.411, 1.038) (-0.305, 1.02) (-0.2, 1.0) };


            \addplot [color=black, style={thick}, only marks, mark=*, mark options={scale=1, fill=red!80}] 
            coordinates {(-0.579, -0.834) (-0.053, -0.751) (0.579, -0.834) (0.842, -0.518) (0.632, 0.028) (0.316, 0.639) (-0.158, 0.558) (-0.421, -0.074) (-0.737, -0.658) (1.042, -0.518) (0.832, 0.028) (0.516, 0.639) (0.621, 0.934) (1.147, 0.9) (1.779, 0.934) (1.937, 0.629) (1.621, 0.037) (1.358, -0.578) (-0.358, 0.558) (-0.621, -0.074) (-0.937, -0.658) (-1.411, -0.646) (-1.726, -0.048) (-2.042, 0.63) (-1.779, 1.066) (-1.253, 1.1) (-0.621, 1.066) 
            };

            \addplot [color=black, style={thick}, only marks, mark=*, mark options={scale=1, fill=blue!80}] 
            coordinates {(-0.56, -0.66) (0.04, -0.51) (0.64, -0.66) (0.4, 0.0) (0.04, 0.54) (-0.14, 0.0) (-0.56, -0.66) (1.213, -0.44) (0.973, 0.22) (0.613, 0.76) (1.213, 0.7) (1.813, 0.76) (1.513, 0.28) (1.213, -0.44) (-0.58, 0.773)  (-0.76, 0.193) (-1.18, -0.427) (-1.51, 0.173) (-1.78, 0.773) (-1.18, 0.833) (-0.58, 0.773) 
            };

            
		\end{axis}

	\end{tikzpicture}}
        \caption{\label{fig:dg_scheme} Schematic of an unstructured discontinuous Galerkin approximation with volume (blue) and surface (red) quadrature nodes.  
        }
    \end{figure}
While the proposed approach to be presented is somewhat agnostic to the numerical method of choice, this work particularly focuses on nodal discontinuous Galerkin (DG)-type~\citep{Hesthaven2008DG} approximations to \cref{eq:gen_hype}, including strong-form representations such as flux reconstruction~\citep{Huynh2007} schemes and collocated discontinuous spectral element methods. In these approaches, the domain $\Omega$ is partitioned into $N$ elements $\Omega_k$ such that $\Omega = \bigcup_N\Omega_k$ and $\Omega_i\cap\Omega_j=\emptyset$ for $i\neq j$. With a slight abuse of notation, we use $\mathbf{w}(\mathbf{x})$  to denote the discrete solution within each element $\Omega_k$, which is approximated by a set of $N_s$ nodal basis functions as 
\begin{equation}
    \mathbf{w}(\mathbf{x}) = \sum_{i=1}^{N_s} \mathbf{w}_{i}\phi_i(\mathbf{x})\subset V_h,
\end{equation}
where $\phi_i(\mathbf{x})$ are the nodal basis functions, $\mathbf{w}_i$ are their associated coefficients, and $V_h$ is the piece-wise polynomial space spanned by the basis functions. The order of the approximation is denoted as $\mathbb P_p$, where $p$ is the maximal order of $\mathbf{w}(\mathbf{x})$. We neglect the temporal domain here for brevity to focus on the spatial discretization and treat temporal derivatives by a method-of-lines approach.

By integrating \cref{eq:gen_hype} with respect to a test function $\mathbf{w}_h(\mathbf{x}) \subset V_h$, the DG formulation is obtained as
\begin{equation}\label{eq:semi-disc}
    \sum_{k=1}^{N} \left \{ \int_{\Omega_k} \partial_t \mathbf{w}{\cdot} \mathbf{w}_h\ \mathrm{d}V  + \int_{\partial \Omega_k} \overline{\mathbf{F}}(\mathbf{w}^-, \mathbf{w}^+, \mathbf{n}) {\cdot} \mathbf{w}_h\ \mathrm{d}\mathbf{s}  -  \int_{\Omega_k} \mathbf{F}(\mathbf{w}){\cdot} \nabla \mathbf{w}_h\ \mathrm{d}V\right \} = \mathbf{0}.
\end{equation}
These integrals are typically performed by quadrature over volume and surface quadrature nodes, shown by the blue and red nodes in \cref{fig:dg_scheme}, respectively. Similar strong form formulations can be obtained through a correction function-type approach~\citep{Huynh2007}, which use corresponding solution and flux nodes akin to volume and surface quadrature nodes. Along each element interface $\partial \Omega_k$, two solutions exist, the solution within the element of interest (denoted by the $\mathbf{w}^-$) and the solution of its face-adjacent neighbor (denoted by $\mathbf{w}^+$). The flux function here is replaced by a numerical interface flux $\overline{\mathbf{F}}(\mathbf{w}^-, \mathbf{w}^+, \mathbf{n})$, which depends on these two solutions as well as the outward facing normal vector $\mathbf{n}$. The numerical flux is generally computed using exact~\citep{Toro1997_4} or approximate Riemann solvers (e.g., \citet{Rusanov1962}).

While the solution of high-order DG schemes generally does not abide by certain key properties of the underlying hyperbolic systems (e.g., maximum principles, invariant regions, etc.), the scheme possesses a useful feature in that the evolution of the element-wise mean of the solution
$\overline{\mathbf{w}}$, defined as 
\begin{equation}
    \overline{\mathbf{w}} = \frac{\int_{\Omega_k}\mathbf{w} (\mathbf{x})\ \mathrm{d}\mathbf{x}}{\int_{\Omega_k} \mathrm{d}\mathbf{x}},
\end{equation}
is dependent on only the interface flux, i.e.,
\begin{equation}
    \partial_t \overline{\mathbf{w}} = - \int_{\partial \Omega_k} \overline{\mathbf{F}}(\mathbf{w}^-, \mathbf{w}^+, \mathbf{n}) \ \mathrm{d}\mathbf{s},
\end{equation}
which allows it to be represented as a convex combination of temporal updates of first-order Godunov schemes~\citep{Godunov1959}. Under some relatively minor assumptions on the numerical scheme (e.g., Courant-–Friedrichs–-Lewy condition, appropriate Riemann solver, explicit strong stability preserving time integration, appropriate quadrature, etc.), this ensures that the element-wise mean preserves convex invariants of the system~\citep{Zhang2010}.

\subsection{Maximum principle preserving limiting}\label{ssec:lmp}
This property of the element-wise mean of high-order DG schemes allows it to be used in many numerical stabilization approaches, some examples of which are referenced as follows ~\citep{Zhang2010,Zhang2011b,Zhang2017,Dzanic2022,Ching2023,Dzanic2023b}.
In the simplest example, the element-wise mean can be used to enforce a maximum principle on the advection equation, given as 
\begin{equation}
    \partial_t u(\mathbf{x},t) + \mathbf{c}{\cdot}\nabla u = 0,
\end{equation}
where $\mathbf{c}$ is the velocity. It is known that the solution of scalar hyperbolic conservation laws abides by a maximum principle in the sense that the extrema are bounded by the initial data $u(\mathbf{x}, 0)$. In a discrete setting, this can be formulated as a discrete maximum principle as 
\begin{equation}
    m \leq u(\mathbf{x}_i,t) \leq M\ \forall\ i \in V, t \geq 0,
\end{equation}
where $m$ and $M$ are the lower and upper bounds of the solution, respectively, and $V$ is some set of nodes over which the maximum principle should be enforced. In general, $V$ corresponds to the discrete locations which the solution needs to be evaluated (e.g., volume/surface quadrature nodes for DG schemes, solution/flux nodes for flux reconstruction schemes, etc.). 

One can apply a global maximum principle, where the bounds are computed over the domain as
\begin{equation}
    m = \underset{\mathbf{x} \in \Omega}{\min}\ u(\mathbf{x}, 0), \quad  M = \underset{\mathbf{x} \in \Omega}{\max}\ u(\mathbf{x}, 0),
\end{equation}
which results in a single set of bounds for the entire solution. The discrete solution within each element at some time step $t^n$ can be limited in an \emph{a priori} (or \emph{a posteriori}) way using the approach of \citet{Zhang2010} as 
\begin{equation}
    \widetilde u^n(\mathbf{x}) = \overline u^n + \alpha \left(u^n(\mathbf{x}) -  \overline u^n\right),
\end{equation}
where $\widetilde u^n(\mathbf{x})$ is the limited solution, $\overline u^n $ is the element-wise mean, and $\alpha$ is an element-wise constant limiting factor. It can be shown that if $\alpha$ is computed as 
\begin{equation}
    \alpha = \min \left[1, 
    \left |\frac{M - \overline{u}}{\underset{i \in V}{\max}\ u(\mathbf{x}_i) - \overline{u}} \right|, 
    \left |\frac{m - \overline{u}}{\underset{i \in V}{\min}\ u(\mathbf{x}_i) - \overline{u}} \right|\right],
\end{equation}
the limited solution is guaranteed to preserve the maximum principle, i.e., $m \leq \widetilde u (\mathbf{x}_i) \leq M \ \forall \ i \in V$.

In practice, this is often not restrictive enough since local extrema can exist in the solution which may not be adequately stabilized by global bounds. It is generally more advisable to enforce a discrete \emph{local} maximum principle, for which the bounds are computed at each time step over some local stencil $S$ for each element $\Omega_k$ as
\begin{equation}
    m_k^{n+1} = \underset{i \in S}{\min}\ u(\mathbf{x}_i, t^n), \quad  M_k^{n+1} = \underset{i \in S}{\max}\ u(\mathbf{x}_i, t^n).
\end{equation}
The local maximum preserving limiting approach requires the user to choose this local stencil $S$ for some corresponding nodal set $V$ to apply the bounds over. Various approaches for this exist, some examples of which are shown in \cref{fig:lmp_stencil}. A common approach is to apply bounds uniformly across an element, with each volume/surface quadrature node having the same bounds, which may be computed over a local domain of influence such as the face/Voronoi neighbors for that element (\cref{fig:lmp_stencil1}) or the numerical stencil for that element (\cref{fig:lmp_stencil2}). More localized/restrictive approaches exist as well, where unique bounds can be computed over a local subset of nodes (e.g., a first-order stencil) within an element and enforced for each single node (\cref{fig:lmp_stencil3}). These different approaches broadly have their own advantages and disadvantages, with less restrictive bounds (e.g., face neighbors) often resulting in less numerical dissipation around smooth regions of the solution and more restrictive bounds (e.g., nodal neighbors) often resulting in better behaved solutions around discontinuities (particularly for very high-order approximations) and vice versa. The numerical intricacies of these approaches as well as their unique drawbacks are described in more detail in \citet{Anderson2017}, although the examples shown here and in the referenced work for local maximum principle limiters of DG approximations of linear transport are by no means exhaustive.

   \begin{figure}[tbhp]
        \centering
        \subfloat[Face neighbors]{\adjustbox{width=0.33\linewidth, valign=b}{\label{fig:lmp_stencil1}     \begin{tikzpicture}[spy using outlines={rectangle, height=3cm,width=2.3cm, magnification=3, connect spies}]
		\begin{axis}[name=plot1,
		    axis line style={draw=none},
		    tick style={draw=none},
		    axis x line=left,
            axis y line=left,
            axis equal image,
            clip mode=individual,
    		xmin=-2.3,
    		xmax=2.3,
    		xticklabels={,,},
    		ymin=-1.1,
    		ymax=1.2,
    		yticklabels={,,},
    		style={font=\Large},
    		scale = 1]

    	\addplot[black,thick]
        coordinates {(-1.0, -1.0) (-0.895, -0.95) (-0.789, -0.906) (-0.684, -0.867) (-0.579, -0.834) (-0.474, -0.806) (-0.368, -0.784) (-0.263, -0.767) (-0.158, -0.756) (-0.053, -0.751) (0.053, -0.751) (0.158, -0.756) (0.263, -0.767) (0.368, -0.784) (0.474, -0.806) (0.579, -0.834) (0.684, -0.867) (0.789, -0.906) (0.895, -0.95) (1.0, -1.0) (1.0, -1.0) (0.947, -0.832) (0.895, -0.672) (0.842, -0.518) (0.789, -0.371) (0.737, -0.231) (0.684, -0.098) (0.632, 0.028) (0.579, 0.147) (0.526, 0.259) (0.474, 0.364) (0.421, 0.463) (0.368, 0.554) (0.316, 0.639) (0.263, 0.716) (0.211, 0.787) (0.158, 0.85) (0.105, 0.907) (0.053, 0.957) (0.0, 1.0) (0.0, 1.0) (-0.053, 0.847) (-0.105, 0.7) (-0.158, 0.558) (-0.211, 0.421) (-0.263, 0.289) (-0.316, 0.163) (-0.368, 0.042) (-0.421, -0.074) (-0.474, -0.185) (-0.526, -0.29) (-0.579, -0.39) (-0.632, -0.485) (-0.684, -0.574) (-0.737, -0.658) (-0.789, -0.737) (-0.842, -0.811) (-0.895, -0.879) (-0.947, -0.942) (-1.0, -1.0) };

    	\addplot[black,thick]
        coordinates {(1.2, -1.0) (1.147, -0.832) (1.095, -0.672) (1.042, -0.518) (0.989, -0.371) (0.937, -0.231) (0.884, -0.098) (0.832, 0.028) (0.779, 0.147) (0.726, 0.259) (0.674, 0.364) (0.621, 0.463) (0.568, 0.554) (0.516, 0.639) (0.463, 0.716) (0.411, 0.787) (0.358, 0.85) (0.305, 0.907) (0.253, 0.957) (0.2, 1.0) (0.2, 1.0) (0.305, 0.98) (0.411, 0.962) (0.516, 0.947) (0.621, 0.934) (0.726, 0.922) (0.832, 0.914) (0.937, 0.907) (1.042, 0.902) (1.147, 0.9) (1.253, 0.9) (1.358, 0.902) (1.463, 0.907) (1.568, 0.914) (1.674, 0.922) (1.779, 0.934) (1.884, 0.947) (1.989, 0.962) (2.095, 0.98) (2.2, 1.0) (2.2, 1.0) (2.147, 0.935) (2.095, 0.865) (2.042, 0.791) (1.989, 0.712) (1.937, 0.629) (1.884, 0.541) (1.832, 0.449) (1.779, 0.353) (1.726, 0.252) (1.674, 0.147) (1.621, 0.037) (1.568, -0.077) (1.516, -0.196) (1.463, -0.319) (1.411, -0.446) (1.358, -0.578) (1.305, -0.714) (1.253, -0.855) (1.2, -1.0) };

    	\addplot[black,thick]
        coordinates {(-0.2, 1.0) (-0.253, 0.847) (-0.305, 0.7) (-0.358, 0.558) (-0.411, 0.421) (-0.463, 0.289) (-0.516, 0.163) (-0.568, 0.042) (-0.621, -0.074) (-0.674, -0.185) (-0.726, -0.29) (-0.779, -0.39) (-0.832, -0.485) (-0.884, -0.574) (-0.937, -0.658) (-0.989, -0.737) (-1.042, -0.811) (-1.095, -0.879) (-1.147, -0.942) (-1.2, -1.0) (-1.2, -1.0) (-1.253, -0.915) (-1.305, -0.828) (-1.358, -0.738) (-1.411, -0.646) (-1.463, -0.552) (-1.516, -0.456) (-1.568, -0.357) (-1.621, -0.256) (-1.674, -0.153) (-1.726, -0.048) (-1.779, 0.059) (-1.832, 0.169) (-1.884, 0.281) (-1.937, 0.395) (-1.989, 0.512) (-2.042, 0.63) (-2.095, 0.751) (-2.147, 0.875) (-2.2, 1.0) (-2.2, 1.0) (-2.095, 1.02) (-1.989, 1.038) (-1.884, 1.053) (-1.779, 1.066) (-1.674, 1.078) (-1.568, 1.086) (-1.463, 1.093) (-1.358, 1.098) (-1.253, 1.1) (-1.147, 1.1) (-1.042, 1.098) (-0.937, 1.093) (-0.832, 1.086) (-0.726, 1.078) (-0.621, 1.066) (-0.516, 1.053) (-0.411, 1.038) (-0.305, 1.02) (-0.2, 1.0) };

            \addplot [color=black, style={thick}, only marks, mark=*, mark options={scale=1, fill=red!80}] 
            coordinates { (1.042, -0.518) (0.832, 0.028) (0.516, 0.639) (0.621, 0.934) (1.147, 0.9) (1.779, 0.934) (1.937, 0.629) (1.621, 0.037) (1.358, -0.578) (-0.358, 0.558) (-0.621, -0.074) (-0.937, -0.658) (-1.411, -0.646) (-1.726, -0.048) (-2.042, 0.63) (-1.779, 1.066) (-1.253, 1.1) (-0.621, 1.066) 
            };

            \addplot [color=black, style={thick}, only marks, mark=*, mark options={scale=1, fill=blue!80}] 
            coordinates {(-0.56, -0.66) (1.213, -0.44) (0.973, 0.22) (0.613, 0.76) (1.213, 0.7) (1.813, 0.76) (1.513, 0.28) (1.213, -0.44) (-0.58, 0.773)  (-0.76, 0.193) (-1.18, -0.427) (-1.51, 0.173) (-1.78, 0.773) (-1.18, 0.833) (-0.58, 0.773) 
            };

            \addplot [color=green!70!black, style={very thick}, only marks, mark=*, mark options={scale=1, fill=red!80}] 
            coordinates {(-0.579, -0.834) (-0.053, -0.751) (0.579, -0.834) (0.842, -0.518) (0.632, 0.028) (0.316, 0.639) (-0.158, 0.558) (-0.421, -0.074) (-0.737, -0.658)
            };
            
            \addplot [color=green!80!black, style={very thick}, only marks, mark=*, mark options={scale=1, fill=blue!80}] 
            coordinates {(-0.56, -0.66) (0.04, -0.51) (0.64, -0.66) (0.4, 0.0) (0.04, 0.54) (-0.14, 0.0) 
            };

            
		\end{axis}

	\end{tikzpicture}}}
        \subfloat[Numerical stencil]{\adjustbox{width=0.33\linewidth, valign=b}{\label{fig:lmp_stencil2}     \begin{tikzpicture}[spy using outlines={rectangle, height=3cm,width=2.3cm, magnification=3, connect spies}]
		\begin{axis}[name=plot1,
		    axis line style={draw=none},
		    tick style={draw=none},
		    axis x line=left,
            axis y line=left,
            axis equal image,
            clip mode=individual,
    		xmin=-2.3,
    		xmax=2.3,
    		xticklabels={,,},
    		ymin=-1.1,
    		ymax=1.2,
    		yticklabels={,,},
    		style={font=\Large},
    		scale = 1]

    	\addplot[black,thick]
        coordinates {(-1.0, -1.0) (-0.895, -0.95) (-0.789, -0.906) (-0.684, -0.867) (-0.579, -0.834) (-0.474, -0.806) (-0.368, -0.784) (-0.263, -0.767) (-0.158, -0.756) (-0.053, -0.751) (0.053, -0.751) (0.158, -0.756) (0.263, -0.767) (0.368, -0.784) (0.474, -0.806) (0.579, -0.834) (0.684, -0.867) (0.789, -0.906) (0.895, -0.95) (1.0, -1.0) (1.0, -1.0) (0.947, -0.832) (0.895, -0.672) (0.842, -0.518) (0.789, -0.371) (0.737, -0.231) (0.684, -0.098) (0.632, 0.028) (0.579, 0.147) (0.526, 0.259) (0.474, 0.364) (0.421, 0.463) (0.368, 0.554) (0.316, 0.639) (0.263, 0.716) (0.211, 0.787) (0.158, 0.85) (0.105, 0.907) (0.053, 0.957) (0.0, 1.0) (0.0, 1.0) (-0.053, 0.847) (-0.105, 0.7) (-0.158, 0.558) (-0.211, 0.421) (-0.263, 0.289) (-0.316, 0.163) (-0.368, 0.042) (-0.421, -0.074) (-0.474, -0.185) (-0.526, -0.29) (-0.579, -0.39) (-0.632, -0.485) (-0.684, -0.574) (-0.737, -0.658) (-0.789, -0.737) (-0.842, -0.811) (-0.895, -0.879) (-0.947, -0.942) (-1.0, -1.0) };

    	\addplot[black,thick]
        coordinates {(1.2, -1.0) (1.147, -0.832) (1.095, -0.672) (1.042, -0.518) (0.989, -0.371) (0.937, -0.231) (0.884, -0.098) (0.832, 0.028) (0.779, 0.147) (0.726, 0.259) (0.674, 0.364) (0.621, 0.463) (0.568, 0.554) (0.516, 0.639) (0.463, 0.716) (0.411, 0.787) (0.358, 0.85) (0.305, 0.907) (0.253, 0.957) (0.2, 1.0) (0.2, 1.0) (0.305, 0.98) (0.411, 0.962) (0.516, 0.947) (0.621, 0.934) (0.726, 0.922) (0.832, 0.914) (0.937, 0.907) (1.042, 0.902) (1.147, 0.9) (1.253, 0.9) (1.358, 0.902) (1.463, 0.907) (1.568, 0.914) (1.674, 0.922) (1.779, 0.934) (1.884, 0.947) (1.989, 0.962) (2.095, 0.98) (2.2, 1.0) (2.2, 1.0) (2.147, 0.935) (2.095, 0.865) (2.042, 0.791) (1.989, 0.712) (1.937, 0.629) (1.884, 0.541) (1.832, 0.449) (1.779, 0.353) (1.726, 0.252) (1.674, 0.147) (1.621, 0.037) (1.568, -0.077) (1.516, -0.196) (1.463, -0.319) (1.411, -0.446) (1.358, -0.578) (1.305, -0.714) (1.253, -0.855) (1.2, -1.0) };

    	\addplot[black,thick]
        coordinates {(-0.2, 1.0) (-0.253, 0.847) (-0.305, 0.7) (-0.358, 0.558) (-0.411, 0.421) (-0.463, 0.289) (-0.516, 0.163) (-0.568, 0.042) (-0.621, -0.074) (-0.674, -0.185) (-0.726, -0.29) (-0.779, -0.39) (-0.832, -0.485) (-0.884, -0.574) (-0.937, -0.658) (-0.989, -0.737) (-1.042, -0.811) (-1.095, -0.879) (-1.147, -0.942) (-1.2, -1.0) (-1.2, -1.0) (-1.253, -0.915) (-1.305, -0.828) (-1.358, -0.738) (-1.411, -0.646) (-1.463, -0.552) (-1.516, -0.456) (-1.568, -0.357) (-1.621, -0.256) (-1.674, -0.153) (-1.726, -0.048) (-1.779, 0.059) (-1.832, 0.169) (-1.884, 0.281) (-1.937, 0.395) (-1.989, 0.512) (-2.042, 0.63) (-2.095, 0.751) (-2.147, 0.875) (-2.2, 1.0) (-2.2, 1.0) (-2.095, 1.02) (-1.989, 1.038) (-1.884, 1.053) (-1.779, 1.066) (-1.674, 1.078) (-1.568, 1.086) (-1.463, 1.093) (-1.358, 1.098) (-1.253, 1.1) (-1.147, 1.1) (-1.042, 1.098) (-0.937, 1.093) (-0.832, 1.086) (-0.726, 1.078) (-0.621, 1.066) (-0.516, 1.053) (-0.411, 1.038) (-0.305, 1.02) (-0.2, 1.0) };

            \addplot [color=black, style={thick}, only marks, mark=*, mark options={scale=1, fill=red!80}] 
            coordinates { 
            (1.042, -0.518) (0.832, 0.028) (0.516, 0.639)(-0.358, 0.558) (-0.621, -0.074) (-0.937, -0.658) 
            };

            \addplot [color=black, style={thick}, only marks, mark=*, mark options={scale=1, fill=red!80}, fill opacity=0.15, draw opacity=0.15] 
            coordinates {  (0.621, 0.934) (1.147, 0.9) (1.779, 0.934) (1.937, 0.629) (1.621, 0.037) (1.358, -0.578) (-1.411, -0.646) (-1.726, -0.048) (-2.042, 0.63) (-1.779, 1.066) (-1.253, 1.1) (-0.621, 1.066) 
            };

            \addplot [color=black, style={thick}, only marks, mark=*, mark options={scale=1, fill=blue!80}, fill opacity=0.15, draw opacity=0.15] 
            coordinates {(-0.56, -0.66) (1.213, -0.44) (0.973, 0.22) (0.613, 0.76) (1.213, 0.7) (1.813, 0.76) (1.513, 0.28) (1.213, -0.44) (-0.58, 0.773)  (-0.76, 0.193) (-1.18, -0.427) (-1.51, 0.173) (-1.78, 0.773) (-1.18, 0.833) (-0.58, 0.773) 
            };
            
            \addplot [color=green!80!black, style={very thick}, only marks, mark=*, mark options={scale=1, fill=blue!80}] 
            coordinates {(-0.56, -0.66) (0.04, -0.51) (0.64, -0.66) (0.4, 0.0) (0.04, 0.54) (-0.14, 0.0) 
            };

            \addplot [color=green!70!black, style={thick}, only marks, mark=*, mark options={scale=1, fill=red!80}] 
            coordinates {(-0.579, -0.834) (-0.053, -0.751) (0.579, -0.834) (0.842, -0.518) (0.632, 0.028) (0.316, 0.639) (-0.158, 0.558) (-0.421, -0.074) (-0.737, -0.658)
            };

            
		\end{axis}

	\end{tikzpicture}}}
        \subfloat[Nodal neighbors]{\adjustbox{width=0.33\linewidth, valign=b}{\label{fig:lmp_stencil3}     \begin{tikzpicture}[spy using outlines={rectangle, height=3cm,width=2.3cm, magnification=3, connect spies}]
		\begin{axis}[name=plot1,
		    axis line style={draw=none},
		    tick style={draw=none},
		    axis x line=left,
            axis y line=left,
            axis equal image,
            clip mode=individual,
    		xmin=-2.3,
    		xmax=2.3,
    		xticklabels={,,},
    		ymin=-1.1,
    		ymax=1.2,
    		yticklabels={,,},
    		style={font=\Large},
    		scale = 1]

    	\addplot[black,thick]
        coordinates {(-1.0, -1.0) (-0.895, -0.95) (-0.789, -0.906) (-0.684, -0.867) (-0.579, -0.834) (-0.474, -0.806) (-0.368, -0.784) (-0.263, -0.767) (-0.158, -0.756) (-0.053, -0.751) (0.053, -0.751) (0.158, -0.756) (0.263, -0.767) (0.368, -0.784) (0.474, -0.806) (0.579, -0.834) (0.684, -0.867) (0.789, -0.906) (0.895, -0.95) (1.0, -1.0) (1.0, -1.0) (0.947, -0.832) (0.895, -0.672) (0.842, -0.518) (0.789, -0.371) (0.737, -0.231) (0.684, -0.098) (0.632, 0.028) (0.579, 0.147) (0.526, 0.259) (0.474, 0.364) (0.421, 0.463) (0.368, 0.554) (0.316, 0.639) (0.263, 0.716) (0.211, 0.787) (0.158, 0.85) (0.105, 0.907) (0.053, 0.957) (0.0, 1.0) (0.0, 1.0) (-0.053, 0.847) (-0.105, 0.7) (-0.158, 0.558) (-0.211, 0.421) (-0.263, 0.289) (-0.316, 0.163) (-0.368, 0.042) (-0.421, -0.074) (-0.474, -0.185) (-0.526, -0.29) (-0.579, -0.39) (-0.632, -0.485) (-0.684, -0.574) (-0.737, -0.658) (-0.789, -0.737) (-0.842, -0.811) (-0.895, -0.879) (-0.947, -0.942) (-1.0, -1.0) };

    	\addplot[black,thick]
        coordinates {(1.2, -1.0) (1.147, -0.832) (1.095, -0.672) (1.042, -0.518) (0.989, -0.371) (0.937, -0.231) (0.884, -0.098) (0.832, 0.028) (0.779, 0.147) (0.726, 0.259) (0.674, 0.364) (0.621, 0.463) (0.568, 0.554) (0.516, 0.639) (0.463, 0.716) (0.411, 0.787) (0.358, 0.85) (0.305, 0.907) (0.253, 0.957) (0.2, 1.0) (0.2, 1.0) (0.305, 0.98) (0.411, 0.962) (0.516, 0.947) (0.621, 0.934) (0.726, 0.922) (0.832, 0.914) (0.937, 0.907) (1.042, 0.902) (1.147, 0.9) (1.253, 0.9) (1.358, 0.902) (1.463, 0.907) (1.568, 0.914) (1.674, 0.922) (1.779, 0.934) (1.884, 0.947) (1.989, 0.962) (2.095, 0.98) (2.2, 1.0) (2.2, 1.0) (2.147, 0.935) (2.095, 0.865) (2.042, 0.791) (1.989, 0.712) (1.937, 0.629) (1.884, 0.541) (1.832, 0.449) (1.779, 0.353) (1.726, 0.252) (1.674, 0.147) (1.621, 0.037) (1.568, -0.077) (1.516, -0.196) (1.463, -0.319) (1.411, -0.446) (1.358, -0.578) (1.305, -0.714) (1.253, -0.855) (1.2, -1.0) };

    	\addplot[black,thick]
        coordinates {(-0.2, 1.0) (-0.253, 0.847) (-0.305, 0.7) (-0.358, 0.558) (-0.411, 0.421) (-0.463, 0.289) (-0.516, 0.163) (-0.568, 0.042) (-0.621, -0.074) (-0.674, -0.185) (-0.726, -0.29) (-0.779, -0.39) (-0.832, -0.485) (-0.884, -0.574) (-0.937, -0.658) (-0.989, -0.737) (-1.042, -0.811) (-1.095, -0.879) (-1.147, -0.942) (-1.2, -1.0) (-1.2, -1.0) (-1.253, -0.915) (-1.305, -0.828) (-1.358, -0.738) (-1.411, -0.646) (-1.463, -0.552) (-1.516, -0.456) (-1.568, -0.357) (-1.621, -0.256) (-1.674, -0.153) (-1.726, -0.048) (-1.779, 0.059) (-1.832, 0.169) (-1.884, 0.281) (-1.937, 0.395) (-1.989, 0.512) (-2.042, 0.63) (-2.095, 0.751) (-2.147, 0.875) (-2.2, 1.0) (-2.2, 1.0) (-2.095, 1.02) (-1.989, 1.038) (-1.884, 1.053) (-1.779, 1.066) (-1.674, 1.078) (-1.568, 1.086) (-1.463, 1.093) (-1.358, 1.098) (-1.253, 1.1) (-1.147, 1.1) (-1.042, 1.098) (-0.937, 1.093) (-0.832, 1.086) (-0.726, 1.078) (-0.621, 1.066) (-0.516, 1.053) (-0.411, 1.038) (-0.305, 1.02) (-0.2, 1.0) };

            \addplot [color=black, style={thick}, only marks, mark=*, mark options={scale=1, fill=red!80}, fill opacity=0.15, draw opacity=0.15] 
            coordinates {(-0.579, -0.834) (-0.053, -0.751) (0.579, -0.834) (0.842, -0.518) (0.632, 0.028) (0.316, 0.639) (-0.158, 0.558) (-0.737, -0.658) (1.042, -0.518) (0.832, 0.028) (0.516, 0.639) (0.621, 0.934) (1.147, 0.9) (1.779, 0.934) (1.937, 0.629) (1.621, 0.037) (1.358, -0.578) (-0.358, 0.558) (-0.621, -0.074) (-0.937, -0.658) (-1.411, -0.646) (-1.726, -0.048) (-2.042, 0.63) (-1.779, 1.066) (-1.253, 1.1) (-0.621, 1.066) 
            };

            \addplot [color=black, style={thick}, only marks, mark=*, mark options={scale=1, fill=blue!80}, fill opacity=0.15, draw opacity=0.15] 
            coordinates {(-0.56, -0.66)  (0.64, -0.66)   (-0.56, -0.66) (1.213, -0.44) (0.973, 0.22) (0.613, 0.76) (1.213, 0.7) (1.813, 0.76) (1.513, 0.28) (1.213, -0.44) (-0.58, 0.773)  (-0.76, 0.193) (-1.18, -0.427) (-1.51, 0.173) (-1.78, 0.773) (-1.18, 0.833) (-0.58, 0.773) 
            };

            \addplot [color=green!70!black, style={very thick}, only marks, mark=*, mark options={scale=1, fill=blue!80}] 
            coordinates {(-0.14, 0.0)
            };
            
            \addplot [color=black, style={thick}, only marks, mark=*, mark options={scale=1, fill=blue!80}] 
            coordinates {(0.04, -0.51) (0.04, 0.54) (0.4, 0.0)
            };
            
            \addplot [color=black, style={thick}, only marks, mark=*, mark options={scale=1, fill=red!80}] 
            coordinates {(-0.421, -0.074) 
            };

            
		\end{axis}

	\end{tikzpicture}}}
        \newline
        \caption{\label{fig:lmp_stencil} Various examples of choices for the stencil to compute bounds over ($S$, shown by solid red/blue nodes) as well as the nodes to enforce bounds over ($V$, shown by nodes with green outlines), including face/Voronoi neighbors (left), numerical stencil (middle), and nodal neighbors (right). 
        }
    \end{figure}
\section{Methodology}\label{sec:methodology}
The overarching goal of this work is to use insights from high-order limiting schemes for the Boltzmann equation to guide the development of limiting schemes for the Euler equations. Essentially, this work aims to find a novel approach to computing and enforcing bounds on the solution of Euler equations through its connection with an underlying kinetic transport equation, which allows for the use of limiting techniques developed for linear advection. Due to the scalar linear nature of the advection equation, these limiting techniques have become significantly more established over the decades and are simpler and more robust than methods required for nonlinear hyperbolic systems of equations, with many options of varying degrees of algorithmic simplicity, robustness, and accuracy. Most importantly, through this approach, one instead recovers local bounds for the \emph{conserved variables} of the solution, i.e.,
\begin{subequations}
\begin{align}
    \rho (\mathbf{x}_i) &\in [\rho_{\min}, \rho_{\max}],\\
    \mathbf{m}_j (\mathbf{x}_i) &\in [\mathbf{m}_{j, \min}, \mathbf{m}_{j, \max}] \quad \text{for} \quad j \in \{1, \ldots, d\},\\
    E (\mathbf{x}_i) &\in [E_{\min}, E_{\max}],
\end{align}
\end{subequations}
for some arbitrary point $\mathbf{x}_i$, which allows for the use of simple scalar limiting techniques and differs from many limiting approaches which enforce bounds on linear/nonlinear functionals of the solution (e.g., density, pressure, entropy, etc.).

\subsection{Boltzmann representation of the Euler equations}
The Euler equations can be represented in terms of the zero Knudsen number limit of the Boltzmann equation. Consider an \emph{explicit} approximation of the Boltzmann equation at some time step $t^n$ (denoted by the superscript $n$), written as
\begin{equation}
    \partial_t f^n (\mathbf{x}, \mathbf{u}) = - \mathbf{u} {\cdot} \nabla f^n + \frac{1}{Kn}\mathcal C(f^n, f^{n'}).
\end{equation}
For zero Knudsen number, the distribution function is in equilibrium, such that we obtain
\begin{equation}
    f^n (\mathbf{x}, \mathbf{u}) = g^n (\mathbf{x}, \mathbf{u}) = g (\mathbf{w}^n (\mathbf{x})),
\end{equation}
i.e., the distribution function is the Maxwell--Boltzmann distribution in \cref{eq:maxwell} for the corresponding macroscopic solution state $\mathbf{w}$. Furthermore, the collision term in the inviscid limit \emph{for a distribution function in equilibrium}, while indeterminate, can be approximated by taking the limit as $Kn \to 0$ as
\begin{equation}
    \underset{Kn \to 0}{\lim}\  \frac{1}{Kn}\mathcal C(g, g') = 0.
\end{equation}
Therefore, a single \emph{explicit} temporal integration step for the Boltzmann equation with a solution initially at equilibrium ($f^n = g^n$) yields the collisionless Boltzmann equation, written as
\begin{equation}
    \partial_t f^n (\mathbf{x}, \mathbf{u}) = - \mathbf{u} {\cdot} \nabla f^n,
\end{equation}
which can be integrated across the velocity domain to yield the Euler equations as 
\begin{equation}
    \int_{\mathbb R^d} \partial_t f^n (\mathbf{x}, \mathbf{u}) \boldsymbol{\psi} (\mathbf{u})\ \mathrm{d}\mathbf{u} = \partial_t \mathbf{w}^n (\mathbf{x}) \quad \text{and} \quad \int_{\mathbb R^d} - \mathbf{u} {\cdot} \nabla f^n (\mathbf{x}, \mathbf{u}) \boldsymbol{\psi} (\mathbf{u})\ \mathrm{d}\mathbf{u} = - \boldsymbol{\nabla}{\cdot}\mathbf{F}\left(\mathbf{w}^n (\mathbf{x})\right).
\end{equation}
Note that this observation is only valid for a single explicit step for an equilibrium state, as any deviations from equilibrium (driven by the transport term) are then immediately reversed by the (infinitely large) collision term at the next time step.

\subsection{Maximum principle bounds for the Boltzmann equation}
The collisionless Boltzmann equation can be considered as a high-dimensional linear transport equation, and for any arbitrary point in the velocity domain $\mathbf{u} = \mathbf{u}_0$, $f_0 = f(\mathbf{u}_0)$, this simply yields a corresponding linear advection equation,
\begin{equation}
    \partial_t f_0(\mathbf{x}, t) +  \mathbf{u}_0 {\cdot} \nabla f_0 = 0. 
\end{equation}
As such, we can apply standard stabilization techniques developed for high-order approximations of linear advection such as a local maximum principle preserving limiter. For example, one can apply element-wise forward-in-time maximum principle bounds as presented in \cref{ssec:lmp}, which, for an arbitrary element $\Omega_k$, enforces bounds on the nodes $j \in V$ computed over the stencil $i \in S$ as 
\begin{equation}
    \underset{i \in S}{\min}\ f_0^n(\mathbf{x}_i) \leq f_0^{n+1}(\mathbf{x}_j) \leq \underset{i \in S}{\max}\ f_0^n(\mathbf{x}_i) \quad \forall \ j \in V.
\end{equation}
We present this in terms a forward Euler approximation, but one can straightforwardly generalize this to higher-order time stepping schemes. Since $\mathbf{u}_0$ is arbitrary, this local maximum principle can be applied uniformly across the velocity domain, i.e.,
\begin{equation}\label{eq:boltzmann_lmp}
    f_{\min}^{n+1}(\mathbf{u})\leq f^{n+1}(\mathbf{x}_j, \mathbf{u}) \leq f_{\max}^{n+1}(\mathbf{u})\quad \forall \ j \in V \ , \ \forall \ \mathbf{u}  \in \mathbb R^d.
\end{equation}
where
\begin{equation}
    f_{\min}^{n+1}(\mathbf{u}) = \underset{i \in S}{\min}\ f^n(\mathbf{x}_i, \mathbf{u}) 
    \quad
    \text{and}
    \quad 
    f_{\max}^{n+1}(\mathbf{u}) = \underset{i \in S}{\max}\ f^n(\mathbf{x}_i, \mathbf{u}).
\end{equation}
We note here that limiting methods of this form are not necessarily just constrained to collisionless Boltzmann approximations, and various techniques that include the effects of collision can also be used (e.g., operator splitting, retaining only positivity bounds~\citep{Jaiswal2022,Dzanic2023b}, etc.). However, this is not necessary in the proposed approach as we are dealing only with kinetic approximations in equilibrium.

\subsection{From kinetic to macroscopic}
Since discrete velocity approximations of the Boltzmann equation are effectively identical to the advection equation, the goal of this approach is to extend limiting techniques (such as the aforementioned maximum principle preserving limiter) from the advection equation to the Euler equations through the Boltzmann equation. The connection between the local maximum principle bounds in \cref{eq:boltzmann_lmp} and the Euler equations can be made through the Maxwell--Boltzmann distribution in \cref{eq:maxwell}, which analytically defines the distribution function as
\begin{equation}
    f^n(\mathbf{x}_i, \mathbf{u}) = g(\mathbf{w}^n(\mathbf{x}_i)).
\end{equation}
An example of this is shown in the schematic in \cref{fig:1d_example}, which showcases an example Euler solution within a DG element (which aims to local maximum principle bounds on the discrete numerical stencil) and the respective distribution function for the discrete solution nodes across that stencil. 

   \begin{figure}[tbhp]
        \centering
        \adjustbox{width=0.54\linewidth, valign=b}{\begin{tikzpicture}[every node/.style={font=\small},scale=2]
    
    \begin{scope}
        \draw[thick,black] (-1,0) -- (1,0);
        \draw[thick,black] (-1,-0.05) -- (-1, 0.05);
        \draw[thick,black] ( 1,-0.05) -- ( 1, 0.05);
        
        \draw[thick,black,densely dashed] (-1.1,-0.05) -- (-1.1, 1.2);
        \draw[thick,black,densely dashed] ( 1.1,-0.05) -- ( 1.1, 1.2);
        
        \draw[thick,black] (-1.2,0) -- (-1.9,0);
        \draw[thick,black,densely dashed] (-1.9,0) -- (-2.05,0);
        \draw[thick,black] ( 1.2,0) -- ( 1.9,0);
        \draw[thick,black,densely dashed] ( 1.9,0) -- ( 2.05,0);
        \draw[thick,black] (-1.2,-0.05) -- (-1.2, 0.05);
        \draw[thick,black] ( 1.2,-0.05) -- ( 1.2, 0.05);
    \end{scope}

    \begin{scope}[on behind layer]
        \draw[thick,black,densely dotted] (-1.200,0.3) -- (-1.200,0);
        \draw[thick,black,densely dotted] (-1.000,0.3) -- (-1.000,0);
        \draw[thick,black,densely dotted] (-0.775,0.3) -- (-0.775,0);
        \draw[thick,black,densely dotted] ( 0.000,0.2) -- ( 0.000,0);
        \draw[thick,black,densely dotted] ( 0.775,0.25) -- ( 0.775,0);
        \draw[thick,black,densely dotted] ( 1.000,0.3) -- ( 1.000,0);
        \draw[thick,black,densely dotted] ( 1.200,0.2) -- ( 1.200,0);

        \draw[thick,Set1-D] plot[smooth,domain=-1:1] (\x,{0.5*(0.24974*(\x)^2 - 0.0645161*\x + 0.4)});    
        \draw[thick,Set1-D] plot[samples at={-0.775, 0.0, 0.775}, color=black, style={thick}, only marks, mark=*, mark options={scale=0.6, fill=blue!80}] (\x,{0.5*(0.24974*(\x)^2 - 0.0645161*\x + 0.4)});
        \draw[thick,Set1-D] plot[samples at={-1, 1}, color=black, style={thick}, only marks, mark=*, mark options={scale=0.6, fill=red!80}] (\x,{0.5*(0.24974*(\x)^2 - 0.0645161*\x + 0.4)});

        \draw[thick,Set1-C] plot[smooth,domain=-1:1] (\x,{-0.1*(\x)^3 - 0.294038*(\x)^2 - 0.2*\x + 1.0});
        \draw[thick,Set1-C] plot[samples at={-0.775, 0.0, 0.775}, color=black, style={thick}, only marks, mark=*, mark options={scale=0.6, fill=blue!80}] (\x,{-0.1*(\x)^3 - 0.294038*(\x)^2 - 0.2*\x + 1.0});
        \draw[thick,Set1-C] plot[samples at={-1, 1}, color=black, style={thick}, only marks, mark=*, mark options={scale=0.6, fill=red!80}] (\x,{-0.1*(\x)^3 - 0.294038*(\x)^2 - 0.2*\x + 1.0});

        \draw[thick,Set1-E!80!black] plot[smooth,domain=-1:1] (\x,{-0.0244652*(\x)^3 + 0.294038*(\x)^2 - 0.125535*\x + 0.455962});
        \draw[thick,Set1-E!80!black] plot[samples at={-0.775, 0.0, 0.775}, color=black, style={thick}, only marks, mark=*, mark options={scale=0.6, fill=blue!80}] (\x,{-0.0244652*(\x)^3 + 0.294038*(\x)^2 - 0.125535*\x + 0.455962});
        \draw[thick,Set1-E!80!black] plot[samples at={-1, 1}, color=black, style={thick}, only marks, mark=*, mark options={scale=0.6, fill=red!80}] (\x,{-0.0244652*(\x)^3 + 0.294038*(\x)^2 - 0.125535*\x + 0.455962});
    \end{scope}
    
    \begin{scope}[on behind layer]
        \draw[thick,Set1-D] plot[smooth,domain=0.3:1] ({\x-2.2},{0.270552*(\x)^2 - 0.0483871*\x + 0.1});
        \draw[thick,Set1-D] plot[samples at={0.775}, color=black, style={thick}, only marks, mark=*, mark options={scale=0.6, fill=blue!80, fill opacity=0.15, draw opacity=0.15}] ({\x-2.2},{0.270552*(\x)^2 - 0.0483871*\x + 0.1});
        \draw[thick,Set1-D] plot[samples at={1}, color=black, style={thick}, only marks, mark=*, mark options={scale=0.6, fill=red!80}] ({\x-2.2},{0.270552*(\x)^2 - 0.0483871*\x + 0.1});
        \draw[thick,Set1-D,densely dashed] plot[smooth,domain=0.15:0.3] ({\x-2.2},{0.270552*(\x)^2 - 0.0483871*\x + 0.1});
        
        \draw[thick,Set1-C] plot[smooth,domain=0.3:1] ({\x-2.2},{-0.122554*(\x)^3 + 0.147019*(\x)^2 + 0.147554*\x + 0.727981});
        \draw[thick,Set1-C] plot[samples at={0.775}, color=black, style={thick}, only marks, mark=*, mark options={scale=0.6, fill=blue!80, fill opacity=0.15, draw opacity=0.15}]({\x-2.2},{-0.122554*(\x)^3 + 0.147019*(\x)^2 + 0.147554*\x + 0.727981});
        \draw[thick,Set1-C] plot[samples at={1}, color=black, style={thick}, only marks, mark=*, mark options={scale=0.6, fill=red!80}]({\x-2.2},{-0.122554*(\x)^3 + 0.147019*(\x)^2 + 0.147554*\x + 0.727981});
        \draw[thick,Set1-C,densely dashed] plot[smooth,domain=0.15:0.3] ({\x-2.2},{-0.122554*(\x)^3 + 0.147019*(\x)^2 + 0.147554*\x + 0.727981});
        
        \draw[thick,Set1-E!80!black] plot[smooth,domain=0.3:1] ({\x-2.2},{0.147019*(\x)^2 + 0.347554*\x + 0.027981});
        \draw[thick,Set1-E!80!black] plot[samples at={0.775}, color=black, style={thick}, only marks, mark=*, mark options={scale=0.6, fill=blue!80, fill opacity=0.15, draw opacity=0.15}] ({\x-2.2},{0.147019*(\x)^2 + 0.347554*\x + 0.027981});
        \draw[thick,Set1-E!80!black] plot[samples at={1}, color=black, style={thick}, only marks, mark=*, mark options={scale=0.6, fill=red!80}] ({\x-2.2},{0.147019*(\x)^2 + 0.347554*\x + 0.027981});
        \draw[thick,Set1-E!80!black,densely dashed] plot[smooth,domain=0.15:0.3] ({\x-2.2},{0.147019*(\x)^2 + 0.347554*\x + 0.027981});

        \draw[thick,Set1-D] plot[smooth,domain=-1:-0.3] ({\x+2.2},{0.145682*(\x)^2 + 0.0483871*\x + 0.1});
        \draw[thick,Set1-D] plot[samples at={-0.775}, color=black, style={thick}, only marks, mark=*, mark options={scale=0.6, fill=blue!80, fill opacity=0.15, draw opacity=0.15}] ({\x+2.2},{0.145682*(\x)^2 + 0.0483871*\x + 0.1});
        \draw[thick,Set1-D] plot[samples at={-1}, color=black, style={thick}, only marks, mark=*, mark options={scale=0.6, fill=red!80}] ({\x+2.2},{0.145682*(\x)^2 + 0.0483871*\x + 0.1});
        \draw[thick,Set1-D,densely dashed] plot[smooth,domain=-0.3:-0.15] ({\x+2.2},{0.145682*(\x)^2 + 0.0483871*\x + 0.1});
        
        \draw[thick,Set1-C] plot[smooth,domain=-1:-0.3] ({\x+2.2},{-0.303915*(\x)^3 + 0.117615*(\x)^2 + 0.303915*\x + 0.482385});
        \draw[thick,Set1-C] plot[samples at={-0.775}, color=black, style={thick}, only marks, mark=*, mark options={scale=0.6, fill=blue!80, fill opacity=0.15, draw opacity=0.15}] ({\x+2.2},{-0.303915*(\x)^3 + 0.117615*(\x)^2 + 0.303915*\x + 0.482385});
        \draw[thick,Set1-C] plot[samples at={-1}, color=black, style={thick}, only marks, mark=*, mark options={scale=0.6, fill=red!80}] ({\x+2.2},{-0.303915*(\x)^3 + 0.117615*(\x)^2 + 0.303915*\x + 0.482385});
        \draw[thick,Set1-C,densely dashed] plot[smooth,domain=-0.3:-0.15] ({\x+2.2},{-0.303915*(\x)^3 + 0.117615*(\x)^2 + 0.303915*\x + 0.482385});
        
        \draw[thick,Set1-E!80!black] plot[smooth,domain=-1:-0.3] ({\x+2.2},{0.517615*(\x)^2 + 0.303915*\x + 0.282385});
        \draw[thick,Set1-E!80!black] plot[samples at={-0.775}, color=black, style={thick}, only marks, mark=*, mark options={scale=0.6, fill=blue!80, fill opacity=0.15, draw opacity=0.15}] ({\x+2.2},{0.517615*(\x)^2 + 0.303915*\x + 0.282385});
        \draw[thick,Set1-E!80!black] plot[samples at={-1}, color=black, style={thick}, only marks, mark=*, mark options={scale=0.6, fill=red!80}] ({\x+2.2},{0.517615*(\x)^2 + 0.303915*\x + 0.282385});
        \draw[thick,Set1-E!80!black,densely dashed] plot[smooth,domain=-0.3:-0.15] ({\x+2.2},{0.517615*(\x)^2 + 0.303915*\x + 0.282385});
    \end{scope}
        
    \node[black] at (-1.200,-0.2) {$\mathbf{w}_1$};
    \node[black] at (-1.000,-0.2) {$\mathbf{w}_2$};
    \node[black] at (-0.775,-0.2) {$\mathbf{w}_3$};
    \node[black] at ( 0.000,-0.2) {$\mathbf{w}_4$};
    \node[black] at ( 0.775,-0.2) {$\mathbf{w}_5$};
    \node[black] at ( 1.000,-0.2) {$\mathbf{w}_6$};
    \node[black] at ( 1.200,-0.2) {$\mathbf{w}_7$};

    \node[black] at (0.9,1.1) {$\Omega_k$};
    \node[black] at (-1.3,1.1) {$\Omega_{k-1}$};
    \node[black] at (1.9,1.1) {$\Omega_{k+1}$};

    \node[Set1-C] at (-0.4,1.15) {\footnotesize $E(x)$};
    \node[Set1-E!80!black] at (-0.4,0.7) {\footnotesize $m(x)$};
    \node[Set1-D] at (-0.4,0.35) {\footnotesize $\rho(x)$};
    
        
        
    
\end{tikzpicture}}
        \adjustbox{width=0.45\linewidth, valign=b}{\begin{tikzpicture}[spy using outlines={rectangle, height=3cm,width=2.3cm, magnification=3, connect spies}]
	\begin{axis}[name=plot1,
		axis line style={latex-latex},
	    axis x line=left,
        axis y line=left,
        width=\axisdefaultwidth,
        height=0.7*\axisdefaultheight,
        scale only axis=true,
		xlabel={$u$},
    	xmin=-4, xmax=6,
    	ylabel={$f$},
    	ymin=0,ymax=0.12,
    	ytick={0, 0.02, 0.04, 0.06, 0.08, 0.10, 0.12},
        y tick label style={
            /pgf/number format/.cd,
            fixed,
            precision=2,
            /tikz/.cd},
        clip mode=individual,
    	legend style={at={(0.97, 0.97)},anchor=north east,font=\small},
    	legend cell align={left},
    	style={font=\normalsize}]
    	
        \addplot[color=red!90!black, style={semithick}]
        table[x=u, y=fa, col sep=comma]{./figs/data/f_example.csv};
        \addlegendentry{$\mathbf{w}_1$}
        \addplot[color=blue!90!black, style={semithick}]
        table[x=u, y=fb, col sep=comma]{./figs/data/f_example.csv};
        \addlegendentry{$\mathbf{w}_2$}
        \addplot[color=green!70!black, style={semithick}]
        table[x=u, y=fc, col sep=comma]{./figs/data/f_example.csv};
        \addlegendentry{$\mathbf{w}_3$}
        \addplot[color=yellow!90!black, style={semithick}]
        table[x=u, y=fd, col sep=comma]{./figs/data/f_example.csv};
        \addlegendentry{$\mathbf{w}_4$}
        \addplot[color=cyan!90!black, style={semithick}]
        table[x=u, y=fe, col sep=comma]{./figs/data/f_example.csv};
        \addlegendentry{$\mathbf{w}_5$}
        \addplot[color=black, style={semithick}]
        table[x=u, y=ff, col sep=comma]{./figs/data/f_example.csv};
        \addlegendentry{$\mathbf{w}_6$}
        \addplot[color=Set1-G!90!black, style={semithick}]
        table[x=u, y=fg, col sep=comma]{./figs/data/f_example.csv};
        \addlegendentry{$\mathbf{w}_7$}

	\end{axis}
\end{tikzpicture}}
        \caption{\label{fig:1d_example} Example schematic for a one-dimensional Euler solution within an element $\Omega_k$, including volume (solid blue circles) and surface (solid red circles) quadrature nodes for the local discrete stencil (denoted by $\{\mathbf{w}_1, \ldots, \mathbf{w}_7 \} \in S$). Corresponding equilibrium distribution function (in the velocity domain) for each solution node shown on the right. 
        }
    \end{figure}

Given the analytic definition of the distribution function for each solution node in the stencil $S$, the bound for the distribution function at the next time step can be computed as 
\begin{equation}\label{eq:fminmax}
     f_{\min}^{n+1}(\mathbf{x}_i, \mathbf{u}) = \underset{i \in S}{\min}\ g(\mathbf{w}^n(\mathbf{x}_i)) \quad \text{and} \quad 
     f_{\max}^{n+1}(\mathbf{x}_i, \mathbf{u}) = \underset{i \in S}{\max}\ g(\mathbf{w}^n(\mathbf{x}_i)),
\end{equation}
which are piecewise defined over the velocity domain and are effectively the lower and upper bounds of the corresponding distribution functions of each solution node in the stencil $S$. A schematic of these distribution function bounds for the corresponding solution and stencil in \cref{fig:1d_example} is shown in \cref{fig:lmp_dist}. This essentially yields an expression for the region of ``admissible'' distribution function values over the velocity domain, bounded from below and above by $f_{\min}$ and $f_{\max}$, respectively. 

   \begin{figure}[tbhp]
        \centering
        \adjustbox{width=0.45\linewidth, valign=b}{\begin{tikzpicture}[spy using outlines={rectangle, height=3cm,width=2.3cm, magnification=3, connect spies}]
	\begin{axis}[name=plot1,
		axis line style={latex-latex},
	    axis x line=left,
        axis y line=left,
        unbounded coords=jump,
        scale only axis=true,
		xlabel={$u$},
    	xmin=-4, xmax=6,
    	ylabel={$f_{\min}$, $f_{\max}$},
    	ymin=0,ymax=0.12,
    	ytick={0, 0.02, 0.04, 0.06, 0.08, 0.10, 0.12},
        y tick label style={
            /pgf/number format/.cd,
            fixed,
            precision=2,
            /tikz/.cd},
        clip mode=individual,
    	legend style={at={(0.97, 0.97)},anchor=north east,font=\small},
    	legend cell align={left},
    	style={font=\normalsize}]

        \addlegendimage{color=red, style={semithick}}
        \addlegendentry{$g(\mathbf{w}_1)$}
        \addlegendimage{color=blue, style={semithick}}
        \addlegendentry{$g(\mathbf{w}_2)$}
        \addlegendimage{color=green, style={semithick}}
        \addlegendentry{$g(\mathbf{w}_3)$}
        \addlegendimage{color=yellow, style={semithick}}
        \addlegendentry{$g(\mathbf{w}_4)$}
        \addlegendimage{color=cyan, style={semithick}}
        \addlegendentry{$g(\mathbf{w}_5)$}
        \addlegendimage{color=black, style={semithick}}
        \addlegendentry{$g(\mathbf{w}_6)$}
        \addlegendimage{color=Set1-G, style={semithick}}
        \addlegendentry{$g(\mathbf{w}_7)$}

        \node[black] at (1, 0.11) {$f_{\max}$};
        \node[black] at (1, 0.025) {$f_{\min}$};
       
        \addplot[color=red!60, style={semithick, densely dashed}]
        table[x=u, y=fa, col sep=comma]{./figs/data/f_example.csv};
        \addplot[color=blue!60, style={semithick, densely dashed}]
        table[x=u, y=fb, col sep=comma]{./figs/data/f_example.csv};
        \addplot[color=green!60, style={semithick, densely dashed}]
        table[x=u, y=fc, col sep=comma]{./figs/data/f_example.csv};
        \addplot[color=yellow!60, style={semithick, densely dashed}]
        table[x=u, y=fd, col sep=comma]{./figs/data/f_example.csv};
        \addplot[color=cyan!60, style={semithick, densely dashed}]
        table[x=u, y=fe, col sep=comma]{./figs/data/f_example.csv};
        \addplot[color=black!60, style={semithick, densely dashed}]
        table[x=u, y=ff, col sep=comma]{./figs/data/f_example.csv};
        \addplot[color=Set1-G!60, style={semithick, densely dashed}]
        table[x=u, y=fg, col sep=comma]{./figs/data/f_example.csv};

        \addplot[color=red!90!black, style={ultra thick}]
        table[x=u, y=fa, col sep=comma]{./figs/data/f_example_max.csv};
        \addplot[color=blue!90!black, style={ultra thick}]
        table[x=u, y=fb, col sep=comma]{./figs/data/f_example_max.csv};
        \addplot[color=green!70!black, style={ultra thick}]
        table[x=u, y=fc, col sep=comma]{./figs/data/f_example_max.csv};
        \addplot[color=yellow!90!black, style={ultra thick}]
        table[x=u, y=fd, col sep=comma]{./figs/data/f_example_max.csv};
        \addplot[color=cyan!90!black, style={ultra thick}]
        table[x=u, y=fe, col sep=comma]{./figs/data/f_example_max.csv};
        \addplot[color=black, style={ultra thick}]
        table[x=u, y=ff, col sep=comma]{./figs/data/f_example_max.csv};
        \addplot[color=Set1-G!90!black, style={ultra thick}]
        table[x=u, y=fg, col sep=comma]{./figs/data/f_example_max.csv};

        \addplot[color=red!90!black, style={ultra thick}]
        table[x=u, y=fa, col sep=comma]{./figs/data/f_example_min.csv};
        \addplot[color=blue!90!black, style={ultra thick}]
        table[x=u, y=fb, col sep=comma]{./figs/data/f_example_min.csv};
        \addplot[color=green!70!black, style={ultra thick}]
        table[x=u, y=fc, col sep=comma]{./figs/data/f_example_min.csv};
        \addplot[color=yellow!90!black, style={ultra thick}]
        table[x=u, y=fd, col sep=comma]{./figs/data/f_example_min.csv};
        \addplot[color=cyan!90!black, style={ultra thick}]
        table[x=u, y=fe, col sep=comma]{./figs/data/f_example_min.csv};
        \addplot[color=black, style={ultra thick}]
        table[x=u, y=ff, col sep=comma]{./figs/data/f_example_min.csv};
        \addplot[color=Set1-G!90!black, style={ultra thick}]
        table[x=u, y=fg, col sep=comma]{./figs/data/f_example_min.csv};

        \addplot+[name path=A, black, opacity=0, mark=none] 
        table[x=u, y=fmin, col sep=comma]{./figs/data/f_example_min.csv};
        \addplot+[name path=B, black, opacity=0, mark=none] 
        table[x=u, y=fmax, col sep=comma]{./figs/data/f_example_max.csv};

        \addplot[black!5] fill between[of=A and B];
	\end{axis}
\end{tikzpicture}}
        \caption{\label{fig:lmp_dist} Example schematic for minimum and maximum bounds for the distribution function ($f_{\min}$ and $f_{\max}$, respectively) for an element $\Omega_k$ as computed by the local discrete stencil of the macroscopic solution $\{\mathbf{w}_1, \ldots, \mathbf{w}_7 \} \in S$ from \cref{fig:1d_example}. Distribution function portions highlighted by the corresponding equilibrium distribution functions. Shaded region denotes the admissible set of the distribution function values. 
        }
    \end{figure}

Once the admissible distribution function values are defined, it is possible to then recover the corresponding admissible macroscopic states. In particular, we are interested in finding the minimum and maximum values of each component of the macroscopic solution $[\rho, \mathbf{m}, E]^T$ that can be recovered from a distribution function $f$ which resides in the admissible region $f_{\min} \leq f \leq f_{\max}$. For the density and energy, where the respective collision invariants are strictly non-negative (i.e., a non-negative measure), we can simply integrate the moments of the minimum/maximum distribution function to yield the minimum/maximum density and energy bounds, i.e.,
    \begin{equation}\label{eq:rhobounds}
    \rho_{\min}^{n+1} = \int_{\mathbb R^d} f_{\min}^{n+1}\ \mathrm{d}\mathbf{u}  \quad \text{and} \quad \rho_{\max}^{n+1} = \int_{\mathbb R^d} f_{\max}^{n+1}\ \mathrm{d}\mathbf{u},
    \end{equation}
    \begin{equation}\label{eq:Ebounds}
    E_{\min}^{n+1} = \int_{\mathbb R^d} f_{\min}^{n+1}\ \left(\frac{1}{2}\mathbf{u}\cdot\mathbf{u}\right) \  \mathrm{d}\mathbf{u}  \quad \text{and} \quad E_{\max}^{n+1} = \int_{\mathbb R^d} f_{\max}^{n+1}\ \left(\frac{1}{2}\mathbf{u}\cdot\mathbf{u}\right)  \ \mathrm{d}\mathbf{u}.
    \end{equation}
For momentum, where the collision invariants are not strictly non-negative, the process is slightly more involved, requiring splitting the integration over the velocity domain as 
    \begin{equation}\label{eq:mmbounds}
        \mathbf{m}_{i,\min}^{n+1} = \int_{\underset{\mathbf{u}_i < 0}{\mathbb R^d}} f_{\max}^{n+1}\ \mathbf{u}_i \  \mathrm{d}\mathbf{u} + \int_{\underset{\mathbf{u}_i \geq 0}{\mathbb R^d}} f_{\min}^{n+1}\ \mathbf{u}_i \  \mathrm{d}\mathbf{u}
    \end{equation}
and
    \begin{equation}\label{eq:mpbounds}
        \mathbf{m}_{i,\max}^{n+1} = \int_{\underset{\mathbf{u}_i < 0}{\mathbb R^d}} f_{\min}^{n+1}\ \mathbf{u}_i \  \mathrm{d}\mathbf{u} + \int_{\underset{\mathbf{u}_i \geq 0}{\mathbb R^d}} f_{\max}^{n+1}\ \mathbf{u}_i \  \mathrm{d}\mathbf{u},
    \end{equation}
for each momentum component $i$. This yields a set of bounds for the conserved variables at the next time step which can be used for the limiting method of choice, which we denote by $\mathbf{w}_{i,\min}$ and $\mathbf{w}_{i,\max}$, respectively, for each solution component $1 \leq i \leq d+2$. The technique for approximating these integrals is later presented in \cref{ssec:int}.

    \subsection{Limiting}
As the bounds for the proposed approach are given in terms of the conserved variables instead of functionals thereof (e.g., pressure, entropy, etc.), the limiting can be performed by straightforwardly treating the individual solution components as scalar local maximum principle limiting tasks. In particular, the limiting method in this work is the squeeze limiter of \citet{Zhang2010} presented in \cref{ssec:lmp}, which uniformly contracts the high-order solution to the element-wise mean. For the proposed approach, the limiting factor $\alpha$ for each solution component $i$ can be computed as
\begin{equation}\label{eq:alpha}
    \alpha_i = \min \left[1, 
    \left |\frac{\mathbf{w}_{i, \max} - \overline{\mathbf{w}}_i}{\underset{j \in V}{\max}\ \mathbf{w}_{i}(\mathbf{x}_j) - \overline{\mathbf{w}}_i} \right|, 
    \left |\frac{\mathbf{w}_{i, \min} - \overline{\mathbf{w}}_i}{\underset{j \in V}{\min}\ \mathbf{w}_{i}(\mathbf{x}_j) - \overline{\mathbf{w}}_i} \right|\right].
\end{equation}
We set the limiting factor explicitly as $\alpha_i = 0$ if $|\underset{j \in V}{\min}\ \mathbf{w}_{i}(\mathbf{x}_j) - \overline{\mathbf{w}}_i| < \epsilon$ or $|\underset{j \in V}{\max}\ \mathbf{w}_{i}(\mathbf{x}_j) - \overline{\mathbf{w}}_i| < \epsilon$, where $\epsilon = 10^{-13}$. This limiting factor is computed separately for each element in the domain. 

As the bounds are decoupled by solution component, the limiting can be performed independently (where each solution component is limited by its respective limiting factor) or coupled (where each solution component is limited by the minimum of all limiting factors). For this work, we use the latter approach, such that the limited solution $\widetilde{\mathbf{w}}(\mathbf{x})$ is obtained as
\begin{equation}\label{eq:limiter}
    \widetilde{\mathbf{w}}(\mathbf{x}) = \overline{\mathbf{w}} + \alpha_{\min} (\mathbf{w}(\mathbf{x}) - \overline{\mathbf{w}}),
\end{equation}
where 
\begin{equation}
    \alpha_{\min} = \underset{i}{\min} \ \alpha_i,
\end{equation}
although one may obtain a marginally less dissipative approach by using the former. It will later be shown in \cref{sec:analysis} that a limiting factor $0 \leq \alpha_{\min} \leq 1$ is guaranteed to exist for which the limited solution discretely satisfies the imposed bounds, i.e.,
\begin{equation*}
    \underset{j \in V}{\min}\ \widetilde{\mathbf{w}}_{i}(\mathbf{x}_j) \geq \mathbf{w}_{i, \min} \quad \text{and} \quad \underset{j \in V}{\max}\ \widetilde{\mathbf{w}}_{i}(\mathbf{x}_j) \leq \mathbf{w}_{i, \max} \quad \text{for} \quad 1 \leq i \leq d+2.
\end{equation*}

\subsection{Relaxing the bounds}
While strictly enforcing a local maximum principle on the high-order scheme is often enough to yield a well-behaved solution around discontinuities (although in some cases more restrictive conditions are needed, see \citet{Anderson2017} for example), it is often too restrictive of a condition in smooth regions of the flow
and will often degrade the accuracy of the scheme around smooth extrema, sometimes referred to as the \emph{peak clipping} effect. This effect is seen in various limiting approaches for high-order approximations of scalar transport and gas dynamics~\citep{Kuzmin2013}, and one remedy for this effect is to marginally relax the bounds to 
recover high-order accuracy in smooth regions of the flow.

However, as will be shown in the numerical experiments, \emph{the proposed approach does not appear to require relaxing the bounds to recover high-order accuracy around smooth extrema whilst simultaneously enforcing sufficient stabilization to ensure convergence around discontinuities.} Nevertheless, it can be beneficial to relax the bounds, both to allow for a slight numerical tolerance in terms of limiting and to compensate for any minor numerical integration errors in computing the bounds from the distribution function. As such, we relax the bounds in a relative manner as 
\begin{equation}\label{eq:relaxed}
     f_{\min}^{n+1}(\mathbf{x}_i, \mathbf{u}) = \underset{i \in S}{\min}\ (1 - r_{\text{fac}}) g(\mathbf{w}^n(\mathbf{x}_i)) \quad \text{and} \quad 
     f_{\max}^{n+1}(\mathbf{x}_i, \mathbf{u}) = \underset{i \in S}{\max}\ (1 + r_{\text{fac}}) g(\mathbf{w}^n(\mathbf{x}_i)),
\end{equation}
where the relaxation factor is set as $r_{\text{fac}} = 10^{-3}$, which is identical to relaxing $f_{\min}^{n+1}$ and $f_{\max}^{n+1}$ themselves. Similar relaxation approaches are taken in other limiting methods (see, for example, \citet{Dumbser2016} and \citet{RuedaRamrez2022}). The effects of enforcing bounds with and without relaxation for smooth and discontinuous solutions will later be shown in the numerical results, but, unless otherwise stated, the relaxed bounds are used in this work.
    
    \subsection{Modeling internal energy}
On its own, the Maxwell--Boltzmann distribution models a monatomic particle which has $d$ translation degrees of freedom. As such, this would constrain the specific heat ratio to $\gamma = 1 + 2/d$. To recover any arbitrary value of $\gamma \in (1, 1 + 2/d]$ and effectively model polyatomic molecules, it is necessary to include the effects of internal degrees of freedom in the distribution function. This is typically achieved through modeling an additional dimension for the internal energy domain, $\zeta \in [0, \infty)$, with the macroscopic state then computed as
\begin{equation}
    \mathbf{w}(\mathbf{x}, t) = 
    \int_{0}^{\infty} \int_{\mathbb R^d} f (\mathbf{x}, \mathbf{u}, \zeta, t)\ \boldsymbol{\psi} (\mathbf{u}, \zeta) \ \mathrm{d}\mathbf{u} \ \mathrm{d}\zeta,
\end{equation}
where the updated collision invariants are defined as $\boldsymbol{\psi} (\mathbf{u}, \zeta) \coloneqq [1, \mathbf{u}, (\mathbf{u}\cdot\mathbf{u})/2 + \zeta]^T$ (see \citet{Baranger2020}, for example). However, naively implementing an internal energy model would add another dimension for evaluating/integrating the distribution function, which can drastically increase the computational cost. Since we are only interested in the effects of the internal energy on the macroscopic state, we can instead take a reduced distribution approach~\citep{Baranger2020} by analytically integrating over the internal energy domain to compute the macroscopic state as
\begin{equation}
    \mathbf{w}(\mathbf{x}, t) = \int_{\mathbb R^d} f_T (\mathbf{x}, \mathbf{u}, \zeta, t)\ [1, \mathbf{u}, (\mathbf{u}\cdot\mathbf{u})/2]^T \ \mathrm{d}\mathbf{u} + \int_{\mathbb R^d} f_I (\mathbf{x}, \mathbf{u}, \zeta, t)\ [0, 0, 1]^T \ \mathrm{d}\mathbf{u},
\end{equation}
where $f_T$ and $f_I$ represent the translational (monatomic) and internal energy components of the distribution function, respectively, defined as
\begin{equation}
    f_T (\mathbf{x}, \mathbf{u}, t) = \int_{0}^{\infty} f (\mathbf{x}, \mathbf{u}, \zeta, t) \ \mathrm{d}\zeta \quad \text{and} \quad f_I (\mathbf{x}, \mathbf{u}, t) = \int_{0}^{\infty} \zeta f (\mathbf{x}, \mathbf{u}, \zeta, t) \ \mathrm{d}\zeta.
\end{equation} This is further simplified through the relation that for an equilibrium distribution function, 
\begin{equation}
    g_I = \frac{\theta \delta}{2}g_T,
\end{equation}
where $\delta$ is the number of internal degrees of freedom and $\theta$ is the scaled temperature of the macroscopic state. By setting $\delta = \frac{2}{\gamma -1} - d$, one can recover the desired specific heat ratio $\gamma$. From a computational perspective, computing the bounds for the distribution function for a particle with internal energy effects (as opposed to the monatomic case) simply requires the addition of the $\frac{\theta \delta}{2}$ term in the collision invariant for the total energy.

    \subsection{Integrating the distribution function}\label{ssec:int}
Computing the bounds requires integrating the distribution functions $f_{\min}$ and $f_{\max}$. Although it is possible in theory to perform this integration analytically, it is not tractable to do it in any computationally feasible manner outside of simple one-dimensional cases. As such, this integration must be performed numerically via quadrature. The quadrature techniques for this approach broadly follow the techniques used in numerical approximations of the Boltzmann equation, where similar integration over the velocity domain is necessary to recover the moments of the distribution function. However, since the distribution functions to be integrated are no longer the Maxwell--Boltzmann distribution (but a piecewise combination thereof), certain techniques which rely on its properties do not directly apply (e.g., Gauss--Hermite quadrature). 

First, the infinite velocity domain can be truncated to a reduced velocity domain $\Omega_u \subset \mathbb R^d$. Since the Maxwell--Boltzmann distribution is a Gaussian distribution centered around $\mathbf{U}$ with a standard deviation of $\sqrt{\theta}$, one can quantify the maximum error of truncating the domain to $k$ standard deviations around the mean (i.e., $\mathbf{U} \pm k \sqrt{\theta}$). For example, $k =4$ and $k =8$ bound the relative error of
\begin{equation*}
    \int_{\mathbb R^d} g \ \mathrm{d}\mathbf{u} - \int_{\Omega_u} g \ \mathrm{d}\mathbf{u}
\end{equation*}
by approximately $6\cdot 10^{-5}$ and $1 \cdot 10^{-15}$, respectively. We use a similar approach to truncate the velocity domain, where the largest extent over the stencil (along each direction) is used, i.e.,
    \begin{equation}\label{eq:truncdomain}
        \Omega_u = \left [ \underset{i \in S}{\min} \ \left(\mathbf{U(\mathbf{x}_i)} - k\sqrt{\theta(\mathbf{x}_i)} \right), \underset{i \in S}{\max} \ \left(\mathbf{U(\mathbf{x}_i)} - k\sqrt{\theta(\mathbf{x}_i)}\right) \right]
    \end{equation}
We assume that $k$ is chosen appropriately to achieve the desired tolerance, which should be related to the relaxation factor used (i.e., $k$ should ostensibly be chosen such that the relative error is $\lesssim r_{\text{fac}}$). For $r_{\text{fac}} = 10^{-3}$ in this work, this implies that $k \geq 3.5$, whereas without relaxation, one should integrate the bounds to machine precision (i.e, $k=8$). We note here however that exact integration (even over the truncated domain) isn't actually necessary, and any numerical integration errors simply result in minor under/overprediction in the limiting factors without affecting the positivity-preserving properties of the approach. This positivity-preserving property of the bounds, which will be proven in \cref{sec:analysis}, holds even under discrete (inexact) integration and is separate from the explicit positivity-preservation enforced, which will be described in \cref{ssec:pp}. 

With the truncated domain defined, one can introduce a quadrature operation as
    \begin{equation}
        \sum_{i = 1}^{N_u} M_i f(\mathbf{u}_i) \psi(\mathbf{u}_i) \approx \int_{\mathbb R^d} f(\mathbf{u}) \psi(\mathbf{u})\ \mathrm{d}\mathbf{u},
    \end{equation}
where $M_i$ is a quadrature weight and $\mathbf{u}_i \in \Omega_u$ correspond to a set of $N_u$ quadrature nodes in the velocity domain. A wide variety of numerical integration techniques can be used here. In this work, we particularly rely on trapezoidal rule integration over an equispaced Cartesian mesh which is commonly employed for discrete velocity approximations of the Boltzman equation. Higher-order (albeit less robust) quadrature techniques can also be used, although the piecewise nature of the underlying distribution function tends to mitigate benefits of many higher-order techniques which have an underlying assumption of the smoothness of the integrand. 

    \subsection{Bounding density and pressure}\label{ssec:pp}
As will be shown in \cref{sec:analysis}, the proposed bounds guarantee that the density and pressure/internal energy will be positive for any $\mathbf{w} \in [\mathbf{w}_{\min}, \mathbf{w}_{\max}]$. However, it is often desirable to enforce some minimum threshold on the density/pressure since just strictly enforcing positivity can lead to numerical difficulties for very small values of density/pressure (e.g., ill-behaved Riemann solver, unfeasible time step restrictions). As such, we enforce a floor value of $10^{-12}$ on the density and pressure. For density, this entails simply replacing the minimum density bound $\rho_{\min}$ with $\max (\rho_{\min}, 10^{-12})$. For pressure, this process is slightly more involved since the proposed bounds are given in terms of the conserved variables. If the discrete limited pressure, computed as 
\begin{equation}\label{eq:pmin}
    P_{\min} = \underset{i \in V}{\min}\  P\left(\widetilde{\mathbf{w}}(\mathbf{x}_i)\right),
\end{equation}
was below the threshold, we apply another limiting step as 
\begin{equation}\label{eq:doublelimit}
    \widetilde{\widetilde{\mathbf{w}}}(\mathbf{x}) = \overline{\mathbf{w}} + \alpha_{P} (\widetilde{\mathbf{w}}(\mathbf{x}) - \overline{\mathbf{w}}),
\end{equation}
with the limiting factor computed as
\begin{equation}
    \alpha_P =
    \left |\frac{10^{-12} - P(\overline{\mathbf{w}})}{P_{\min} - P(\overline{\mathbf{w}})} \right|,
\end{equation}
The quasiconcavity of the pressure functional ensures that $\underset{i \in V}{\min}\  P\left(\widetilde{\widetilde{\mathbf{w}}}(\mathbf{x}_i)\right) \geq 10^{-12}$.

    \subsection{Overview}

With the details of the proposed approach presented, we summarize here an overview of the method applied to high-order approximations of the Euler equations. This overview is shown in terms of a single temporal update $t^n \to t^{n+1}$ for an arbitrary element $\Omega_k$, but this can directly be generalized to every element and each substage of higher-order (SSP) time integration schemes. The necessary code for the implementation of this approach is provided in the electronic supplementary material. 

\begin{enumerate}
    \item Collect the macroscopic solution $\mathbf{w}^n(\mathbf{x}_i)$ at each node in the chosen stencil $S$. 
    \item Evaluate the distribution function bounds $f^{n+1}_{\min}$ and $f^{n+1}_{\max}$ (per \cref{eq:fminmax} and \cref{eq:relaxed}) and the collision invariants $\boldsymbol{\psi} = [1, \mathbf{u}, \mathbf{u}{\cdot}\mathbf{u}/2 + \zeta]$ at each discrete velocity node $\mathbf{u}_i$. Perform a quadrature summation over the discrete velocity nodes ($1 \leq i \leq N_u$) to recover the macroscopic bounds per \cref{eq:rhobounds,eq:Ebounds,eq:mmbounds,eq:mpbounds}.
    \item Apply a minimum threshold on the density bound. 
    \item Advance solution in time $\mathbf{w}^n(\mathbf{x}) \to \mathbf{w}^{n+1}(\mathbf{x})$.
    \item Compute the limited solution per \cref{eq:limiter}.
    \item Compute the minimum pressure of the limited solution per \cref{eq:pmin}.
    \item If minimum pressure is below the minimum threshold, apply a second limiting step per \cref{eq:doublelimit}. 
    \item Replace the high-order DG solution $\mathbf{w}^{n+1}(\mathbf{x})$ with the limited solution $\widetilde{\mathbf{w}}^{n+1}(\mathbf{x})$) (or $\widetilde{\widetilde{\mathbf{w}}}^{n+1}(\mathbf{x})$). 
\end{enumerate}

\section{Analysis}
\label{sec:analysis}
With this formulation, we now move on to state and prove the numerical properties of the proposed approach. The subsequent proofs rely on some (or all) of the following assumptions. 
\begin{enumerate}[(i)]
    \item Let there be a stencil $S$ over which the bounds are computed for an arbitrary element $\Omega_k$, and let $G$ be the set of solutions with positive density $\rho > 0$ and positive pressure/internal energy $P, \rho e > 0$. The initial solution $\mathbf{w}_i \in G \ \forall \ i \in S$ at some time $t^n$.
    \item The distribution functions are integrated via a discrete integration operator $\mathbf{M}$ consisting of $N_u$ entries with positive weights corresponding to integration over a subdomain $\Omega_u \subset \mathbb R^d$ for which we denote its extent as $u_{\max} = \underset{\mathbf{u} \in \Omega_u}{\max}\ \|\mathbf{u}\|_2$ which is assumed to be at least two thermal/absolute velocity units (i.e., \cref{eq:truncdomain} with $k \geq 2$). It is assumed that the integration operator can integrate moments of an arbitrary distribution function $f$ to some desired numerical tolerance, i.e.,
    \begin{equation}
        \int_{\mathbb R^d} f(\mathbf{u}) \psi(\mathbf{u})\ \mathrm{d}\mathbf{u} \approx \sum_{i = 1}^{N_u} M_i f(\mathbf{u}_i) \psi(\mathbf{u}_i).
    \end{equation}
    We denote this operation in shorthand as $\mathbf{M}\cdot f \otimes \psi(\mathbf{u})$. 
    \item The nodes of the stencil over which to compute bounds (denoted by $S$) and the nodes over which the limiting is performed (denoted by $V \subset S$) are chosen appropriately such that: 
    \begin{enumerate}
        \item There exists a set of positive quadrature weights $c_i$ such that 
        \begin{equation*}
            \overline{\mathbf{w}} = \frac{\int_{\Omega_k} \mathbf{w}(\mathbf{x})\ \mathrm{d}\mathbf{x}}{\int_{\Omega_k} \mathrm{d}\mathbf{x}} = \frac{\sum_{i \in V} c_i \mathbf{w}(\mathbf{x}_i)}{\sum_{i \in V} c_i}.
        \end{equation*}
        \item There exists a subset of nodes $V' \subset V$ (e.g., surface quadrature nodes) with a corresponding set of positive quadrature weights $c_i'$ such that 
        \begin{equation*}
             \int_{\partial \Omega_k} \mathbf{w}(\mathbf{x})\ \mathrm{d}\mathbf{s} \approx \sum_{i \in V'} c_i' \mathbf{w}(\mathbf{x}_i).
        \end{equation*}
        See, for example, \citet{Zhang2011} and \citet{Carlier2023}.
    \end{enumerate}
    \item The interface numerical fluxes are computed using a Rusanov-type approach~\citep{Rusanov1962} as
    \begin{equation*}
        \overline{\mathbf{F}}(\mathbf{w}^+, \mathbf{w}^-, \mathbf{n}) = \frac{1}{2}\left [\mathbf{F}(\mathbf{w}^-) + \mathbf{F}(\mathbf{w}^+) \right]\cdot {\mathbf{n}} - \frac{1}{2}\lambda_{\max}(\mathbf{w}^+, \mathbf{w}^-)\left [\mathbf{w}^+ -\mathbf{w}^- \right],
    \end{equation*}
    where $\lambda_{\max}$ is an upper bound on the local maximum wavespeed of the system. 
\end{enumerate}

\begin{theorem}[Constant state preservation]\label{thm:const}
    The proposed bounds preserve constant states in the sense that if $\mathbf{w}(\mathbf{x}_i, t^n) = \mathbf{w}_0 \ \forall\  i \in S$ for some constant state $\mathbf{w}_0$, then $\mathbf{w}_{\min}^{n+1} = \mathbf{w}_{\max}^{n+1} = \mathbf{w}_0$.
\end{theorem}

\newproof{pot_const}{Proof of \cref{thm:const}}
\begin{pot_const}
The proof follows directly from the definition of $f_{\min}^{n+1}$ and $f_{\max}^{n+1}$ in \cref{eq:fminmax}, which yields 
\begin{equation}
    f_{\min}^{n+1} = f_{\max}^{n+1} = g(\mathbf{w}_0). 
\end{equation}
By integrating the moments, one yields
\begin{equation}
    \mathbf{w}_{\min}^{n+1} = \mathbf{w}_{\max}^{n+1} = \mathbf{M}\cdot g(\mathbf{w}_0) \otimes \psi(\mathbf{u}) \approx \int_{\mathbb R^d} g(\mathbf{w}_0) \ \boldsymbol{\psi} (\mathbf{u}) \ \mathrm{d}\mathbf{u} = \mathbf{w}_0.
\end{equation}
\end{pot_const}

\begin{theorem}[Positivity preservation]\label{thm:pp} The proposed bounds preserve positivity of density and pressure/internal energy for all $\mathbf{w} \in [\mathbf{w}_{\min}^{n+1}, \mathbf{w}_{\max}^{n+1}]$.
\end{theorem}

\newproof{pot_pp}{Proof of \cref{thm:pp}}
\begin{pot_pp}
To show positivity of density, it suffices to show that $\rho_{\min}^{n+1} > 0$, which can be expressed as the condition
\begin{equation}
    \mathbf{M} \cdot f_{\min}^{n+1} > 0.
\end{equation}
From \cref{eq:maxwell}, the distribution function $g(\mathbf{w})$ is strictly positive if $\mathbf{w} \in G$. Therefore, $f_{\min}^{n+1}(\mathbf{u}), f_{\max}^{n+1}(\mathbf{u}) > 0$ for all $\mathbf{u} \in \mathbb R^d$, which means that $\rho_{\min}^{n+1} > 0$.

To show positivity of pressure/internal energy for all $\mathbf{w} \in [\mathbf{w}_{\min}^{n+1}, \mathbf{w}_{\max}^{n+1}]$, we use the relation that any $\mathbf{w}_0 \in [\mathbf{w}_{\min}^{n+1}, \mathbf{w}_{\max}^{n+1}]$ can be expressed as 
\begin{equation}
    \mathbf{w}_0 = \mathbf{M}\cdot f_0^{n+1} \otimes \boldsymbol{\psi}(\mathbf{u})
\end{equation}
for some distribution function $f_{\min}^{n+1}(\mathbf{u}) \leq f_0^{n+1}(\mathbf{u}) \leq f_{\max}^{n+1}(\mathbf{u})$, where $f_0^{n+1}(\mathbf{u}) > 0$. The internal energy can be calculated as 
\begin{equation}
    \rho e = E - \frac{1}{2}\rho \mathbf{U}{\cdot}\mathbf{U} = \mathbf{M} \cdot f_0^{n+1}\left (\frac{1}{2} \mathbf{u}{\cdot}\mathbf{u} -  \frac{1}{2} \mathbf{u}{\cdot}\mathbf{U}\right)
    = \mathbf{M} \cdot \frac{1}{2}f_0^{n+1} \mathbf{c}{\cdot}\mathbf{c},
\end{equation}
where $\mathbf{c} = \mathbf{u} - \mathbf{U}$ is the peculiar velocity. From the positivity of the distribution function and measure $\mathbf{c}{\cdot}\mathbf{c}$, it follows that $\rho e > 0$ and, by extension, $P = (\gamma-1)\rho e > 0$. 
\end{pot_pp}

\begin{theorem}[Riemann-averaged states]\label{thm:ras}
For any arbitrary $\mathbf{w}_i, \mathbf{w}_j\in G$ and corresponding bounds $[f_{\min}, f_{\max}]$ and $[\mathbf{w}_{\min}, \mathbf{w}_{\max}]$, the Riemann-averaged state
\begin{equation}
    \overline{\mathbf{w}}_{ij} = \frac{1}{2}\left(\mathbf{w}_i + \mathbf{w}_j \right) - \frac{1}{2 \lambda_{\max}}\left(\mathbf{F}(\mathbf{w}_j) - \mathbf{F}(\mathbf{w}_i) \right)\cdot{\mathbf{n}_{ij}}
\end{equation}
preserves the bounds in the sense that $\overline{\mathbf{w}}_{ij} \in [\mathbf{w}_{\min}, \mathbf{w}_{\max}]$ for any arbitrary unit vector $\mathbf{n}_{ij}$ if $\lambda_{\max}$ is an upper bound on the local maximum wavespeed of the system. 
\end{theorem}
\newproof{pot_ras}{Proof of \cref{thm:ras}}
\begin{pot_ras}
By using the relation that 
\begin{equation}
    \mathbf{F}(\mathbf{w}) = \int_{\mathbb R^d} \mathbf{u} g(\mathbf{w})\ \boldsymbol{\psi}(\mathbf{u})\ \mathrm{d}\mathbf{u},
\end{equation}
the Riemann-averaged state can be expressed in terms of a distribution function $\overline{g}_{ij}$ as 
\begin{equation}
    \overline{\mathbf{w}}_{ij} = \mathbf{M} \cdot \overline{g}_{ij} \otimes \boldsymbol{\psi} (\mathbf{u}) = \mathbf{M} \cdot \left [ \frac{1}{2}\left(g_i + g_j \right) - \frac{1}{2 \lambda_{\max}}\left(g_j - g_i \right)\mathbf{u}{\cdot}{\mathbf{n}_{ij}} \right] \otimes \boldsymbol{\psi} (\mathbf{u}).
\end{equation}
Let $\lambda_{\max} = u_{\max} \geq \max\ \left(|\mathbf{U}_i{\cdot}\mathbf{n}_{ij}| + \sqrt{\gamma \theta_i}, |\mathbf{U}_j{\cdot}\mathbf{n}_{ij}| + \sqrt{\gamma \theta_j}\right)$. The Riemann-averaged distribution function can be expressed as
\begin{equation}
    \overline{g}_{ij} = g_i \left (\frac{1}{2} + \frac{1}{2}\frac{\mathbf{u}{\cdot}\mathbf{n}}{\lambda_{\max}} \right) + g_j \left (\frac{1}{2} - \frac{1}{2}\frac{\mathbf{u}{\cdot}\mathbf{n}}{\lambda_{\max}} \right).
\end{equation}
Since $\lambda_{\max} \geq |\mathbf{u}{\cdot}\mathbf{n}_{ij}|$, $\overline{g}_{ij}$ is a convex combination of the states $g_i$ and $g_j$. Therefore,
\begin{equation}
    f_{\min} = \min\ (g_i, g_j) \leq \overline{g}_{ij} \leq \max\ (g_i, g_j) = f_{\max},
\end{equation}
from which it directly follows that $\overline{\mathbf{w}}_{ij} \in [\mathbf{w}_{\min}, \mathbf{w}_{\max}]$. 
\end{pot_ras}

\begin{remark}
    The proof of \cref{thm:ras} is based on an upper bound on the local maximum particle speed $u_{\max}$ (on the truncated domain $\Omega_u$) instead of the local maximum macroscopic wavespeed (e.g., the Davis~\citep{Davis1988} wavespeed estimate $\lambda_{\max} = |\mathbf{U}{\cdot}\mathbf{n}| + \sqrt{\gamma \theta}$). At high Mach numbers, where the thermal velocity is small compared to the flow velocity, these estimates do not drastically differ, but at low Mach numbers, the local maximum particle speed may be several times larger than the maximum macroscopic wavespeed --- a factor of $k/\sqrt{\gamma}$, where $k$ is the extent factor of the truncated domain  --- which can affect the maximum admissible time step required to satisfy the bounds. However, we observed in empirical studies that setting $\lambda_{\max}$ as the local maximum macroscopic wavespeed still always ensured that the Riemann-averaged state satisfied the bounds, but we do not have a rigorous mathematical proof for this behavior.
\end{remark}

\begin{theorem}[First-order scheme]\label{thm:fo}
For any arbitrary $\mathbf{w}_0^n, \mathbf{w}_1^n, \mathbf{w}_2^n \in G$ and corresponding bounds $[\mathbf{w}_{\min}^{n+1}, \mathbf{w}_{\max}^{n+1}]$, the first-order Lax--Friedrichs scheme, 
\begin{equation}
    \mathbf{w}_1^{n+1} = \mathbf{w}_1^n - \frac{\Delta t}{\Delta x}\left ( \overline{\mathbf{F}}(\mathbf{w}_1^n, \mathbf{w}_2^n) - \overline{\mathbf{F}}(\mathbf{w}_0^n, \mathbf{w}_1^n) \right) ,
\end{equation}
where 
\begin{equation}
    \overline{\mathbf{F}}(\mathbf{w}_a, \mathbf{w}_b) = \frac{1}{2}\left [\mathbf{F}(\mathbf{w}_a) + \mathbf{F}(\mathbf{w}_b) \right] - \frac{1}{2}\lambda_{\max}(\mathbf{w}_a, \mathbf{w}_b)\left [\mathbf{w}_b -\mathbf{w}_a \right]
\end{equation}
preserves bounds in the sense that $\mathbf{w}_1^{n+1} \in [\mathbf{w}_{\min}^{n+1}, \mathbf{w}_{\max}^{n+1}]$ under the standard Courant-–Friedrichs–-Lewy (CFL) condition $\lambda_{\max}\frac{\Delta t}{\Delta x} \leq 1$.
\end{theorem}
\newproof{pot_fo}{Proof of \cref{thm:fo}}
\begin{pot_fo}

The first-order Lax--Friedrichs scheme can be written in terms of the Riemann-averaged states as 
\begin{equation}
    \mathbf{w}_1^{n+1} = \left(1 - \lambda_{\max}^{0,1} \frac{\Delta t}{\Delta x} - \lambda_{\max}^{1,2} \frac{\Delta t}{\Delta x}\right)\mathbf{w}_1^n + \lambda_{\max}^{0,1} \frac{\Delta t}{\Delta x} \overline{\mathbf{w}}_{01} + \lambda_{\max}^{1,2} \frac{\Delta t}{\Delta x} \overline{\mathbf{w}}_{12},
\end{equation}
where $\lambda_{\max}^{0,1} = \lambda_{\max}(\mathbf{w}_0, \mathbf{w}_1)$, $\lambda_{\max}^{1,2} = \lambda_{\max}(\mathbf{w}_1, \mathbf{w}_2)$, and $\mathbf{n}_{01} = \mathbf{n}_{12} = 1$. If $\Delta t$ is sufficiently small such that $\lambda_{\max}\frac{\Delta t}{\Delta x} \leq 1$, where $\lambda_{\max} = \max (\lambda_{\max}^{0,1}, \lambda_{\max}^{1,2})$, then $\mathbf{w}_1^{n+1}$ is a convex combination of states in $[\mathbf{w}_{\min}^{n+1}, \mathbf{w}_{\max}^{n+1}]$. As a result, $\mathbf{w}_1^{n+1} \in [\mathbf{w}_{\min}^{n+1}, \mathbf{w}_{\max}^{n+1}]$.

\end{pot_fo}

\begin{theorem}[Mean bounds preservation]\label{thm:mbp}
For some arbitrary element $\Omega_k$ and its associated mean $\overline{\mathbf{w}}^n = \overline{\mathbf{w}}(t^n)$, there exists a sufficiently small $\Delta t > 0$ such that $\overline{\mathbf{w}}(t^n + \Delta t) = \overline{\mathbf{w}}^{n+1} \in [\mathbf{w}_{\min}^{n+1}, \mathbf{w}_{\max}^{n+1}]$. 
\end{theorem}
\newproof{pot_mbp}{Proof of \cref{thm:mbp}}
\begin{pot_mbp}

The element-wise mean can be computed as 
\begin{equation}
    \overline{\mathbf{w}}^n = \frac{\sum_{i \in V} c_i \left [ \mathbf{M} \cdot g(\mathbf{w}_i^n) \otimes \boldsymbol{\psi} (\mathbf{u})\right]}{\sum_{i \in V} c_i},
\end{equation}
which, by linearity of the quadrature operators, can be expressed as 
\begin{equation}
    \overline{\mathbf{w}}^n = \mathbf{M} \cdot \left [ \frac{\sum_{i \in V} c_i  g(\mathbf{w}_i^n)}{\sum_{i \in V} c_i} \right ] \otimes \boldsymbol{\psi} (\mathbf{u}),
\end{equation}
Since, by definition, $g(\mathbf{w}_i^n) \in [f_{\min}^{n+1}, f_{\max}^{n+1}] \ \forall \ \mathbf{u} \in \mathbb R^d \ \forall \ i \in V$ and $[f_{\min}^{n+1}, f_{\max}^{n+1}]$ is a convex set, the convex combination $\sum_{i \in V} c_i  g(\mathbf{w}_i^n)/\sum_{i \in V} c_i \in [f_{\min}^{n+1}, f_{\max}^{n+1}] \ \forall \ \mathbf{u} \in \mathbb R^d \ \forall \ i \in V$. Therefore, $\overline{\mathbf{w}}^n \in [\mathbf{w}_{\min}^{n+1}, \mathbf{w}_{\max}^{n+1}]$.

It now remains to be shown that $\overline{\mathbf{w}}^{n+1} \in [\mathbf{w}_{\min}^{n+1}, \mathbf{w}_{\max}^{n+1}]$. The temporal update of the element-wise mean can be represented as
\begin{equation}
    \overline{\mathbf{w}}^{n+1} = \overline{\mathbf{w}}^n - \Delta t \sum_{i \in V'} c_i' \overline{\mathbf{F}}(\mathbf{w}^{n, -}_i, \mathbf{w}^{n, +}_i, \mathbf{n}_i) =  \frac{\sum_{i \in V} c_i \mathbf{w}^{n}_i}{\sum_{i \in V} c_i} - \Delta t \sum_{i \in V'} c_i' \overline{\mathbf{F}}(\mathbf{w}^{n, -}_i, \mathbf{w}^{n, +}_i, \mathbf{n}_i) ,
\end{equation}
where $\mathbf{w}^{n, -}_i, \mathbf{w}^{n, +}_i \in [\mathbf{w}_{\min}^{n+1}, \mathbf{w}_{\max}^{n+1}]$ for all $i \in V'$. We use the shorthand $\mathbf{w}_i$ to denote $\mathbf{w}(\mathbf{x}_i)$. $\overline{\mathbf{w}}^{n+1}$ can be equivalently expressed as a convex combination as 
\begin{equation}
    \overline{\mathbf{w}}^{n+1} = \frac{\sum_{i \in V\setminus V'} c_i \mathbf{w}^{n}_i + \sum_{i \in  V'} c_i \mathbf{w}^{*}_i}{\sum_{i \in V} c_i },
\end{equation}
where
\begin{equation}
    \mathbf{w}^{*}_i = \mathbf{w}^{n, -}_i - \Delta t \frac{c'_i}{c_i}\overline{\mathbf{F}}(\mathbf{w}^{n, -}_i, \mathbf{w}^{n, +}_i, \mathbf{n}_i).
\end{equation}
Thus, it is sufficient to show that the intermediate states $\mathbf{w}^{*}_i \in [\mathbf{w}_{\min}^{n+1}, \mathbf{w}_{\max}^{n+1}]$ for all $i \in V'$. These states can be represented in terms of the Riemann-averaged states as
\begin{equation}
     \mathbf{w}^{*}_i =  \left [1 - C \right]\mathbf{w}^{n, -}_i + C \overline{\mathbf{w}}_{i^- i^+},
\end{equation}
where 
\begin{equation}
    C = \Delta t \lambda_{\max}(\mathbf{w}^{n, -}_i, \mathbf{w}^{n, +}_i)\frac{c'_i}{c_i}.
\end{equation}
It can be seen that $\mathbf{w}^{*}_i \in [\mathbf{w}_{\min}^{n+1}, \mathbf{w}_{\max}^{n+1}]$ for any $0 \leq C \leq 1$. Note here that $c_i'/c_i$ is $\mathcal O(1/h)$, where $h$ is some local mesh scale, as it is a surface quadrature weight normalized by a volume quadrature weight. Therefore, the above condition is met under a CFL-like constraint $\Delta t \lambda_{\max}(\mathbf{w}^{n, -}_i, \mathbf{w}^{n, +}_i)\frac{c'_i}{c_i} \leq 1$, which then guarantees that $\overline{\mathbf{w}}^{n+1} \in [\mathbf{w}_{\min}^{n+1}, \mathbf{w}_{\max}^{n+1}]$.
\end{pot_mbp}

\begin{corollary}[Existence of solution]\label{thm:eos}
If $\overline{\mathbf{w}}(t^{n+1}) \in [\mathbf{w}_{\min}^{n+1}, \mathbf{w}_{\max}^{n+1}]$, there exists a limiting factor $\alpha \in [0,1]$ such that $\widetilde{\mathbf{w}}(\mathbf{x}_i, t^{n+1}) \in [\mathbf{w}_{\min}^{n+1}, \mathbf{w}_{\max}^{n+1}]\ \forall\ i \in V$.
\end{corollary}
\section{Implementation}\label{sec:implementation}
The proposed limiting approach was implemented using the flux reconstruction (FR) scheme of \citet{Huynh2007}, a strong-form generalization of the nodal DG method, within PyFR~\citep{Witherden2014, Witherden2024}, a high-order unstructured FR solver that can efficiently target massively-parallel CPU and GPU computing architectures. To showcase generalizability with respect to solution node sets, we utilize open solution/flux node sets (i.e., where the interface flux nodes are not a subset of the interior solution nodes), namely the Gauss--Legendre/Gauss--Legendre solution/flux nodes for line/quadrilateral elements and Williams--Shunn/Gauss--Legendre solution/flux nodes for triangular elements~\citep{Williams2014}. Temporal integration was performed using a three-stage, third-order strong stability preserving Runge--Kutta scheme (SSP-RK3~\citep{Gottlieb2001}), with limiting performed after each substage of the scheme. Common interface numerical fluxes were computed using the HLLC Riemann solver of \citet{Toro1994}. We note here that while the numerical analysis in this work only presents some theoretical guarantees for Rusanov-type Riemann solvers, it was found empirically in the numerical experiments that many other positivity-preserving entropy-stable Riemann solvers always preserved the bounds on the element-wise mean (under the standard assumptions), such that less dissipative approaches like HLLC were found to be more performant. 

The bounds were computed over the numerical stencil $S$, consisting of the interior solution and flux nodes as well as the exterior (face-adjacent) flux nodes, and limiting was performed over the element stencil $V$, consisting of just the interior solution and flux nodes, as shown in \cref{fig:lmp_stencil2}. The distribution functions were integrated over the velocity domain using uniform trapezoidal rule quadrature with the domain extent set to four thermal/absolute velocity units (i.e., $k=4$). Similar accuracy was observed with more advanced quadrature techniques such as Romberg integration. For the one- and two-dimensional experiments to be presented, the velocity domain resolution was fixed separately. For one-dimensional problems, we purposely use a large highly-resolved velocity domain ($k = 8$, $N_u  = 2048$) which can integrate the bounds to essentially near machine precision levels so that we can effectively compare the accuracy of the approach using strict and relaxed bounds without introducing any implicit relaxation in the bounds due to numerical under-integration of the distribution functions. For two-dimensional problems, we use a more appropriate resolution of $k = 4$, $N_u = 32^2$ for every problem, although we note that similar results could be obtained with a lower resolution of $N_u = 16^2$. In theory, one can use more efficient adaptive quadrature techniques to yield a more performant algorithm since the piecewise nature of the distribution function gives \emph{a priori} information about the required local resolution, but this is not explored in this work as its primary purpose is to evaluate the accuracy of the proposed approach. 
\section{Results}\label{sec:results}

\subsection{Sod shock tube}

\begin{figure}[htbp!] 
    \centering
    \begin{tabular}{|c | c | c | c | c | c | c | c | c |}
    \hline
    $N$ & $\mathbb P_2$ & Rate & $\mathbb P_3$ & Rate &  $\mathbb P_4$ & Rate &  $\mathbb P_5$ & Rate \\ 
    \hline
    25 & $1.258 \times 10^{-2}$ & - & $9.788 \times 10^{-3}$ & - & $1.222 \times 10^{-2}$ & - & $1.196 \times 10^{-2}$ & - \\
    50 & $4.513 \times 10^{-3}$ & 1.479 & $4.392 \times 10^{-3}$ & 1.156 & $4.621 \times 10^{-3}$ & 1.402 & $4.921 \times 10^{-3}$ & 1.282 \\
    100 & $2.424 \times 10^{-3}$ & 0.897 & $2.178 \times 10^{-3}$ & 1.012 & $2.443 \times 10^{-3}$ & 0.920 & $2.503 \times 10^{-3}$ & 0.975 \\
    200 & $1.259 \times 10^{-3}$ & 0.945 & $1.181 \times 10^{-3}$ & 0.884 & $1.359 \times 10^{-3}$ & 0.846 & $1.576 \times 10^{-3}$ & 0.667 \\
    400 & $6.227 \times 10^{-4}$ & 1.016 & $5.747 \times 10^{-4}$ & 1.039 & $7.750 \times 10^{-4}$ & 0.810 & $8.608 \times 10^{-4}$ & 0.873 \\
    800 & $3.287 \times 10^{-4}$ & 0.922 & $2.886 \times 10^{-4}$ & 0.994 & $4.877 \times 10^{-4}$ & 0.668 & $4.786 \times 10^{-4}$ & 0.847 \\
    1600 & $1.731 \times 10^{-4}$ & 0.925 & $1.469 \times 10^{-4}$ & 0.974 & $2.961 \times 10^{-4}$ & 0.720 & $2.911 \times 10^{-4}$ & 0.717 \\
    3200 & $9.132 \times 10^{-5}$ & 0.923 & $7.815 \times 10^{-5}$ & 0.911 & $1.819 \times 10^{-4}$ & 0.703 & $1.677 \times 10^{-4}$ & 0.795 \\
    \hline
    & & 0.987 & & 0.989 & & 0.834 & & 0.851 \\
    \hline
    \end{tabular}
    \captionof{table}{\label{tab:sod_error_rel} 
     Convergence of the $L^1$ norm of the density error for the Sod shock tube problem at $t = 0.2$ computed using varying approximation orders and mesh resolution \emph{with relaxed bounds}. Average rate of convergence shown on bottom.}
\end{figure}

\begin{figure}[htbp!] 
    \centering
    \begin{tabular}{|c | c | c | c | c | c | c | c | c |}
    \hline
    $N$ & $\mathbb P_2$ & Rate & $\mathbb P_3$ & Rate &  $\mathbb P_4$ & Rate &  $\mathbb P_5$ & Rate \\ 
    \hline    
    25 & $1.273 \times 10^{-2}$ & - & $1.030 \times 10^{-2}$ & - & $1.231 \times 10^{-2}$ & - & $1.309 \times 10^{-2}$ & - \\
    50 & $4.602 \times 10^{-3}$ & 1.468 & $4.879 \times 10^{-3}$ & 1.078 & $5.212 \times 10^{-3}$ & 1.240 & $5.602 \times 10^{-3}$ & 1.224 \\
    100 & $2.496 \times 10^{-3}$ & 0.882 & $2.475 \times 10^{-3}$ & 0.979 & $2.902 \times 10^{-3}$ & 0.845 & $3.104 \times 10^{-3}$ & 0.852 \\
    200 & $1.314 \times 10^{-3}$ & 0.926 & $1.403 \times 10^{-3}$ & 0.819 & $1.715 \times 10^{-3}$ & 0.759 & $1.979 \times 10^{-3}$ & 0.650 \\
    400 & $6.586 \times 10^{-4}$ & 0.996 & $7.438 \times 10^{-4}$ & 0.916 & $9.820 \times 10^{-4}$ & 0.804 & $1.130 \times 10^{-3}$ & 0.808 \\
    800 & $3.529 \times 10^{-4}$ & 0.900 & $4.100 \times 10^{-4}$ & 0.859 & $5.689 \times 10^{-4}$ & 0.788 & $6.749 \times 10^{-4}$ & 0.744 \\
    1600 & $1.887 \times 10^{-4}$ & 0.903 & $2.349 \times 10^{-4}$ & 0.804 & $3.292 \times 10^{-4}$ & 0.789 & $4.132 \times 10^{-4}$ & 0.708 \\
    3200 & $1.001 \times 10^{-4}$ & 0.914 & $1.424 \times 10^{-4}$ & 0.722 & $1.936 \times 10^{-4}$ & 0.766 & $2.593 \times 10^{-4}$ & 0.673 \\
    \hline
    & & 0.969 & & 0.879 & & 0.830 & & 0.784 \\
    \hline
    \end{tabular}
    \captionof{table}{\label{tab:sod_error_norel} 
     Convergence of the $L^1$ norm of the density error for the Sod shock tube problem at $t = 0.2$ computed using varying approximation orders and mesh resolution \emph{without relaxed bounds}. Average rate of convergence shown on bottom. }
\end{figure}

An initial evaluation of the shock capturing and discontinuity resolving capabilities of the proposed approach was performed with the Sod shock tube~\citep{Sod1978}. This canonical gas dynamics test case consists of the three main features of the Riemann problem: shock waves, rarefaction waves, and contact discontinuities. The domain was set as $\Omega=[0,1]$ with Dirichlet boundary conditions, and the initial conditions were set as

    \begin{equation}
        \mathbf{q}(x,0) = \begin{cases}
            \mathbf{q}_l, &\mbox{if } x\leqslant 0.5,\\
            \mathbf{q}_r, &\mbox{else},
        \end{cases} \quad \mathrm{where} \quad \mathbf{q}_l = \begin{bmatrix}
            1 \\ 0 \\ 1
        \end{bmatrix}, \quad \mathbf{q}_r = \begin{bmatrix}
            0.125 \\ 0 \\ 0.1
        \end{bmatrix}.
    \end{equation}

The problem was solved using varying approximation orders and mesh resolution, and the convergence of the density error with respect to the exact solution was computed, calculated as the point-mean $L^1$ norm on the solution nodes. The errors and rates are shown in \cref{tab:sod_error_rel} as computed with the \emph{relaxed} bounds. It can be seen that the optimal rate of unity was nearly achieved at all approximation orders. A marginally higher convergence rate was observed at lower approximation orders ($\mathbb P_2$ and $\mathbb P_3$) than at higher approximation orders ($\mathbb P_4$ and $\mathbb P_5$), which is consistent with observations that element-wise local maximum principle limiters are better behaved at lower approximation orders~\citep{Anderson2017}. The experiments were also repeated with \emph{strict} bounds, the results of which are tabulated in \cref{tab:sod_error_norel}. It can be seen that very similar error magnitudes were obtained with and without relaxing the bounds, with the relaxed bounds even achieving marginally higher convergence rates.

\subsection{Smooth convecting density pulse}

    \begin{figure}[htbp!] 
        \centering
        \begin{tabular}{|c | c | c | c | c | c | c | c | c |}
        \hline
        $N$ & $\mathbb P_2$ & Rate & $\mathbb P_3$ & Rate &  $\mathbb P_4$ & Rate &  $\mathbb P_5$ & Rate \\ 
        \hline
        10 & $1.663 \times 10^{-1}$ & - & $3.622 \times 10^{-2}$ & - & $5.097 \times 10^{-3}$ & - & $6.811 \times 10^{-4}$ & - \\
        15 & $6.048 \times 10^{-2}$ & 2.495 & $4.934 \times 10^{-3}$ & 4.917 & $4.187 \times 10^{-4}$ & 6.164 & $4.657 \times 10^{-5}$ & 6.617 \\
        20 & $2.605 \times 10^{-2}$ & 2.927 & $1.346 \times 10^{-3}$ & 4.515 & $8.759 \times 10^{-5}$ & 5.439 & $4.800 \times 10^{-6}$ & 7.898 \\
        25 & $1.014 \times 10^{-2}$ & 4.227 & $4.821 \times 10^{-4}$ & 4.601 & $2.517 \times 10^{-5}$ & 5.588 & $1.958 \times 10^{-6}$ & 4.020 \\
        30 & $5.177 \times 10^{-3}$ & 3.689 & $1.925 \times 10^{-4}$ & 5.036 & $8.523 \times 10^{-6}$ & 5.940 & $6.527 \times 10^{-7}$ & 6.025 \\
        35 & $2.551 \times 10^{-3}$ & 4.590 & $1.008 \times 10^{-4}$ & 4.196 & $4.250 \times 10^{-6}$ & 4.514 & $2.559 \times 10^{-7}$ & 6.075 \\
        40 & $1.435 \times 10^{-3}$ & 4.311 & $5.225 \times 10^{-5}$ & 4.920 & $2.331 \times 10^{-6}$ & 4.499 & $1.021 \times 10^{-7}$ & 6.877 \\
        \hline
        & & 3.469 & & 4.696 & & 5.570 & & 6.248 \\
        \hline
        \end{tabular}
        \captionof{table}{\label{tab:pulse_error_rel} 
          Convergence of the $L^\infty$ norm of the density error for the convecting density pulse problem at $t = 1$ computed using varying approximation orders and mesh resolution \emph{with relaxed bounds}. Average rate of convergence shown on bottom.}
    \end{figure}

    \begin{figure}[htbp!] 
        \centering
        \begin{tabular}{|c | c | c | c | c | c | c | c | c |}
        \hline
        $N$ & $\mathbb P_2$ & Rate & $\mathbb P_3$ & Rate &  $\mathbb P_4$ & Rate &  $\mathbb P_5$ & Rate \\ 
        \hline
        10 & $1.663 \times 10^{-1}$ & - & $3.622 \times 10^{-2}$ & - & $5.098 \times 10^{-3}$ & - & $6.812 \times 10^{-4}$ & - \\
        15 & $6.048 \times 10^{-2}$ & 2.495 & $4.934 \times 10^{-3}$ & 4.917 & $4.187 \times 10^{-4}$ & 6.164 & $4.650 \times 10^{-5}$ & 6.621 \\
        20 & $2.605 \times 10^{-2}$ & 2.927 & $1.346 \times 10^{-3}$ & 4.515 & $8.773 \times 10^{-5}$ & 5.433 & $4.752 \times 10^{-6}$ & 7.928 \\
        25 & $1.014 \times 10^{-2}$ & 4.227 & $4.821 \times 10^{-4}$ & 4.602 & $2.503 \times 10^{-5}$ & 5.620 & $1.958 \times 10^{-6}$ & 3.974 \\
        30 & $5.177 \times 10^{-3}$ & 3.689 & $1.926 \times 10^{-4}$ & 5.033 & $8.401 \times 10^{-6}$ & 5.988 & $5.215 \times 10^{-7}$ & 7.255 \\
        35 & $2.551 \times 10^{-3}$ & 4.590 & $1.007 \times 10^{-4}$ & 4.202 & $4.155 \times 10^{-6}$ & 4.568 & $2.559 \times 10^{-7}$ & 4.620 \\
        40 & $1.435 \times 10^{-3}$ & 4.311 & $5.231 \times 10^{-5}$ & 4.908 & $2.183 \times 10^{-6}$ & 4.821 & $1.021 \times 10^{-7}$ & 6.877 \\
        \hline
        & & 3.469 & & 4.695 & & 5.605 & & 6.290 \\
        \hline
        \end{tabular}
        \captionof{table}{\label{tab:pulse_error_norel} 
          Convergence of the $L^\infty$ norm of the density error for the convecting density pulse problem at $t = 1$ computed using varying approximation orders and mesh resolution \emph{without relaxed bounds}. Average rate of convergence shown on bottom.}
    \end{figure}
    
To verify that the proposed approach retains the high-order accuracy of the underlying DG scheme even at smooth extrema, a smooth convecting density pulse was simulated. The problem was evaluated on the periodic domain $\Omega=[-0.5, 0.5]$, and the initial conditions are given as
    \begin{subequations}
        \begin{align}
            \rho & = 1 + \exp( - \sigma x^2)  \\
            u &= 1, \\
            P &= 1,
        \end{align}
    \end{subequations}
where $\sigma = 500$ is the strength of the Gaussian wave. This problem consists of a density perturbation convecting through the domain at constant velocity. To evaluate the behavior of the limiting approach at smooth extrema, the $L^\infty$ norm of the density error was calculated at the solution nodes after one flow-through of the domain ($t = 1$). The error and convergence rates at varying approximation orders and mesh resolution are shown in \cref{tab:pulse_error_rel} and \cref{tab:pulse_error_norel} for the limiting with and without relaxing the bounds, respectively. It can be seen that the approach recovers the expected high-order (i.e., approximately $p+1$) convergence rates for all approximation orders. Most importantly, it can also be seen that strictly enforcing the bounds does not degrade the accuracy of the scheme at smooth extrema, indicating that the proposed approach is not overly restrictive in smooth regions of the flow. This is in contrast to many limiting techniques which require ad hoc corrections to recover high-order accuracy at smooth extrema.

\subsection{Double expansion wave}

   \begin{figure}[tbhp]
        \centering
        \subfloat[Density]{\adjustbox{width=0.4\linewidth, valign=b}{\begin{tikzpicture}[spy using outlines={rectangle, height=3cm,width=2.3cm, magnification=3, connect spies}]
	\begin{axis}[name=plot1,
		axis line style={latex-latex},
	    axis x line=left,
        axis y line=left,
		xlabel={$x$},
    	xmin=0, xmax=1,
    	xtick={0, 0.2, 0.4, 0.6, 0.8, 1.0},
    	ylabel={$\rho$},
    	ymin=0,ymax=1.05,
    	ytick={0, 0.2, 0.4, 0.6, 0.8, 1.0},
        clip mode=individual,
    	ylabel style={rotate=-90},
    	legend style={at={(0.5, 0.95)},anchor=north,font=\small},
    	legend cell align={left},
    	style={font=\normalsize}]
    	
        \addplot[color=black!80, style={thin}, only marks, mark=o, mark options={scale=0.5}]
        table[x=x, y=r, col sep=comma]{./figs/data/123_exact.csv};
        
        \addplot[color=blue!80!black, 
                style={thin, dashed}]
        table[x=x, y=r, col sep=comma]{./figs/data/123_p2_n50.csv};
        
        \addplot[color=red!80!black, 
                style={semithick, densely dotted}]
        table[x=x, y=r, col sep=comma]{./figs/data/123_p2_n100.csv};
        
        \addplot[color=green!50!black, 
                style={thin, dash dot}]
        table[x=x, y=r, col sep=comma]{./figs/data/123_p2_n200.csv};
        
        \addplot[color=black, 
                style={thin}]
        table[x=x, y=r, col sep=comma]{./figs/data/123_p2_n400.csv};

	\end{axis}
\end{tikzpicture}}}
        \subfloat[Velocity]{\adjustbox{width=0.4\linewidth, valign=b}{\begin{tikzpicture}[spy using outlines={rectangle, height=3cm,width=2.3cm, magnification=3, connect spies}]
	\begin{axis}[name=plot1,
		axis line style={latex-latex},
	    axis x line=left,
        axis y line=left,
		xlabel={$x$},
    	xmin=0, xmax=1,
    	xtick={0, 0.2, 0.4, 0.6, 0.8, 1.0},
    	ylabel={$U$},
    	ymin=-2.05,ymax=2.05,
        clip mode=individual,
    	ylabel style={rotate=-90},
    	legend style={at={(0.03, 0.97)},anchor=north west,font=\small},
    	legend cell align={left},
    	style={font=\normalsize}]
    	
        \addplot[color=black!80, style={thin}, only marks, mark=o, mark options={scale=0.5}]
        table[x=x, y=u, col sep=comma]{./figs/data/123_exact.csv};
        \addlegendentry{Exact};
        
        \addplot[color=blue!80!black, 
                style={thin, dashed}]
        table[x=x, y=u, col sep=comma]{./figs/data/123_p2_n50.csv};
        \addlegendentry{$N = 50$};
        
        \addplot[color=red!80!black, 
                style={semithick, densely dotted}]
        table[x=x, y=u, col sep=comma]{./figs/data/123_p2_n100.csv};
        \addlegendentry{$N = 100$};
        
        \addplot[color=green!50!black, 
                style={thin, dash dot}]
        table[x=x, y=u, col sep=comma]{./figs/data/123_p2_n200.csv};
        \addlegendentry{$N = 200$};
        
        \addplot[color=black, 
                style={thin}]
        table[x=x, y=u, col sep=comma]{./figs/data/123_p2_n400.csv};
        \addlegendentry{$N = 400$};

	\end{axis}
\end{tikzpicture}}}
        \newline
        \subfloat[Pressure]{\adjustbox{width=0.4\linewidth, valign=b}{\begin{tikzpicture}[spy using outlines={rectangle, height=3cm,width=2.3cm, magnification=3, connect spies}]
	\begin{axis}[name=plot1,
		axis line style={latex-latex},
	    axis x line=left,
        axis y line=left,
		xlabel={$x$},
    	xmin=0, xmax=1,
    	xtick={0, 0.2, 0.4, 0.6, 0.8, 1.0},
    	ylabel={$P$},
    	ymin=0,ymax=0.41,
        clip mode=individual,
    	ylabel style={rotate=-90},
    	legend style={at={(0.5, 0.95)},anchor=north,font=\small},
    	legend cell align={left},
    	style={font=\normalsize}]
    	
        \addplot[color=black!80, style={thin}, only marks, mark=o, mark options={scale=0.5}]
        table[x=x, y=p, col sep=comma]{./figs/data/123_exact.csv};
        
        \addplot[color=blue!80!black, 
                style={thin, dashed}]
        table[x=x, y=p, col sep=comma]{./figs/data/123_p2_n50.csv};
        
        \addplot[color=red!80!black, 
                style={semithick, densely dotted}]
        table[x=x, y=p, col sep=comma]{./figs/data/123_p2_n100.csv};
        
        \addplot[color=green!50!black, 
                style={thin, dash dot}]
        table[x=x, y=p, col sep=comma]{./figs/data/123_p2_n200.csv};
        
        \addplot[color=black, 
                style={thin}]
        table[x=x, y=p, col sep=comma]{./figs/data/123_p2_n400.csv};

	\end{axis}
\end{tikzpicture}}}
        \subfloat[Specific internal energy]{\adjustbox{width=0.4\linewidth, valign=b}{\begin{tikzpicture}[spy using outlines={rectangle, height=3cm,width=2.3cm, magnification=3, connect spies}]
	\begin{axis}[name=plot1,
		axis line style={latex-latex},
	    axis x line=left,
        axis y line=left,
		xlabel={$x$},
    	xmin=0,xmax=1,
    	xtick={0, 0.2, 0.4, 0.6, 0.8, 1.0},
    	ylabel={$e$},
    	ymin=0,ymax=1.05,
        clip mode=individual,
    	ylabel style={rotate=-90},
    	legend style={at={(1.0, 0.03)},anchor=south east,font=\small},
    	legend cell align={left},
    	style={font=\normalsize}]
    	
        \addplot[color=black!80, style={thin}, only marks, mark=o, mark options={scale=0.5}]
        table[x=x, y=e, col sep=comma]{./figs/data/123_exact.csv};
        
        \addplot[color=blue!80!black, 
                style={thin, dashed}]
        table[x=x, y=e, col sep=comma]{./figs/data/123_p2_n50.csv};
        
        \addplot[color=red!80!black, 
                style={semithick, densely dotted}]
        table[x=x, y=e, col sep=comma]{./figs/data/123_p2_n100.csv};
        
        \addplot[color=green!50!black, 
                style={thin, dash dot}]
        table[x=x, y=e, col sep=comma]{./figs/data/123_p2_n200.csv};
        
        \addplot[color=black, 
                style={thin}]
        table[x=x, y=e, col sep=comma]{./figs/data/123_p2_n400.csv};

	\end{axis}
\end{tikzpicture}}}
        \newline
        \caption{\label{fig:123} 
        Density (top left), velocity (top right), pressure (bottom left), and specific internal energy (bottom right) profiles for the double expansion wave problem at $t = 0.15$ computed using a $\mathbb P_2$ approximation with varying mesh resolution.
        }
    \end{figure}
    
The positivity-preserving properties of the proposed approach at near-vacuum conditions were then evaluated through the double expansion wave problem. This problem, known as the \textit{123 problem} (or Test 2) in \citet{Toro1997_4}, was solved on the domain $\Omega = [0,1]$ with Dirichlet boundary conditions, and the initial conditions were set as 
\begin{equation}\label{eq:expansion}
    \mathbf{q}(x, 0) = 
    \begin{cases}
    \mathbf{q}_L, \quad \mathrm{if } x \leq 0.5,\\
    \mathbf{q}_R, \quad \mathrm{else},
    \end{cases}
    \quad \mathrm{where} \quad 
    \mathbf{q}_L = 
    \begin{bmatrix}
           1 \\
           -2 \\
           0.4
    \end{bmatrix},
    \quad
    \mathbf{q}_R = 
    \begin{bmatrix}
           1 \\
           2 \\
           0.4
    \end{bmatrix}.
\end{equation}
This Riemann problem develops into two outrunning expansion waves, causing a low density and pressure region in the center. The resulting symmetric wave structure and the near-vacuum center region are notably difficult to resolve, with even robust first-order schemes often performing poorly or failing altogether.

The problem was solved using a $\mathbb P_2$ approximation with varying mesh resolution between $N = 50$ and $N = 400$ elements. The profiles of density, velocity, pressure, and specific internal energy at $t = 0.15$ are shown in \cref{fig:123}. The predicted solution profiles showed clear qualitative convergence to the exact solution with increasing resolution, even in the near-vacuum center region. The specific internal energy profile showed a distinct spike in the center, a very commonly encountered numerical artifact in this problem indicative of spurious physical entropy generation in the center~\citep{Toro1994}. However, as expected, this spike decreased with increasing resolution, also converging toward the exact solution. 

\subsection{Shu--Osher problem}
The Shu--Osher problem~\citep{Shu1988} was then simulated, which consists of more complex flow behavior including shock waves and oscillatory small-scale features. The domain was set as $\Omega = [0,10]$ with Dirichlet boundary conditions, and the initial conditions were set as 
    \begin{equation}
         \mathbf{q}(x,0) =  \begin{cases}
            \mathbf{q}_l, &\mbox{if } x\leqslant 1, \\
            \mathbf{q}_r, &\mbox{else},
        \end{cases} \quad \mathrm{where} \quad
        \mathbf{q}_l = \begin{bmatrix}
            3.857143 \\ 2.629369 \\ 10.333333
        \end{bmatrix}, \quad
        \mathbf{q}_r = \begin{bmatrix}
            1 + 0.2\sin{5x} \\ 0\\ 1
        \end{bmatrix}.
    \end{equation}
The problem consists of a shock wave interacting with a sinusoidally-perturbed density field which induces instabilities in the flow field that may be erroneously damped by overly dissipative shock capturing schemes. 
    
\begin{figure}[tbhp]
    \centering
    \subfloat[$\mathbb P_1$]{ 
    \adjustbox{width=0.49\linewidth, valign=b}{     \begin{tikzpicture}[spy using outlines={rectangle, height=3.5cm,width=3cm, magnification=2, connect spies}]
		\begin{axis}[name=plot1,
		    axis line style={latex-latex},
		    axis x line=left,
            axis y line=left,
		    clip mode=individual,
		    xlabel={$x$},
    		xmin=0,
    		xmax=10,
    		x tick label style={
        		/pgf/number format/.cd,
            	fixed,
            	precision=1,
        	    /tikz/.cd},
    		ylabel={$\rho$},
    		ylabel style={rotate=-90},
    		ytick={0,1,2,3,4,5},
    		ymin=0,
    		ymax=5,
    		y tick label style={
        		/pgf/number format/.cd,
            	fixed,
            	fixed zerofill,
            	precision=0,
        	    /tikz/.cd},
    		legend style={at={(1.02,1.0)},anchor=north east,font=\tiny},
    		legend cell align={left},
    		legend columns=1,
    		style={font=\normalsize}]
    		
            \addplot[color=black!80, style={thin}, only marks, mark=o, mark options={scale=0.5}, mark repeat = 7, mark phase = 0,
            x filter/.expression={x < 5.5 ? \pgfmathresult : NaN)}
            ]
				table[x=x,y=r,col sep=comma,unbounded coords=jump]{./figs/data/shu-osher_ref.csv};
    		\addlegendentry{Ref.}
    		
            \addplot[color=blue!80!black, 
                    style={thin, densely dashed}]
            table[x=x, y=r, col sep=comma]{./figs/data/shu-osher_p1_n100.csv};
            \addlegendentry{$N = 100$};
            
            \addplot[color=red!80!black, 
                    style={semithick, densely dotted}]
            table[x=x, y=r, col sep=comma]{./figs/data/shu-osher_p1_n200.csv};
            \addlegendentry{$N = 200$};
            
            \addplot[color=black, 
                    style={thin}]
            table[x=x, y=r, col sep=comma]{./figs/data/shu-osher_p1_n400.csv};
            \addlegendentry{$N = 400$};
    		
            \addplot[color=black!80, style={thin}, only marks, mark=o, mark options={scale=0.5}, mark repeat = 2, mark phase = 0,
            x filter/.expression={abs(x - 6.5) < 1 ? \pgfmathresult : NaN)}
            ]
				table[x=x,y=r,col sep=comma,unbounded coords=jump]{./figs/data/shu-osher_ref.csv};
    
            \addplot[color=black!80, style={thin}, only marks, mark=o, mark options={scale=0.5}, mark repeat = 7, mark phase = 0,
            x filter/.expression={x > 7.5 ? \pgfmathresult : NaN)}
            ]
				table[x=x,y=r,col sep=comma,unbounded coords=jump]{./figs/data/shu-osher_ref.csv};

			\coordinate (spypoint) at (axis cs:6.5,3.85);
            \coordinate (magnifyglass) at (axis cs:2.35,1.65);

		\end{axis} 		
	
        \spy [black] on (spypoint) in node[fill=white] at (magnifyglass);
		
	\end{tikzpicture}}}
    \subfloat[$\mathbb P_3$]{
    \adjustbox{width=0.49\linewidth, valign=b}{     \begin{tikzpicture}[spy using outlines={rectangle, height=3.5cm,width=3cm, magnification=2, connect spies}]
		\begin{axis}[name=plot1,
		    axis line style={latex-latex},
		    axis x line=left,
            axis y line=left,
		    clip mode=individual,
		    xlabel={$x$},
    		xmin=0,
    		xmax=10,
    		x tick label style={
        		/pgf/number format/.cd,
            	fixed,
            	precision=1,
        	    /tikz/.cd},
    		ylabel={$\rho$},
    		ylabel style={rotate=-90},
    		ytick={0,1,2,3,4,5},
    		ymin=0,
    		ymax=5,
    		y tick label style={
        		/pgf/number format/.cd,
            	fixed,
            	fixed zerofill,
            	precision=0,
        	    /tikz/.cd},
    		legend style={at={(1.02,1.0)},anchor=north east,font=\tiny},
    		legend cell align={left},
    		style={font=\normalsize}]
    		
            \addplot[color=black!80, style={thin}, only marks, mark=o, mark options={scale=0.5}, mark repeat = 7, mark phase = 0,
            x filter/.expression={x < 5.5 ? \pgfmathresult : NaN)}
            ]
				table[x=x,y=r,col sep=comma,unbounded coords=jump]{./figs/data/shu-osher_ref.csv};
    		\addlegendentry{Ref.}
    		
            \addplot[color=blue!80!black, 
                    style={thin, densely dashed}]
            table[x=x, y=r, col sep=comma]{./figs/data/shu-osher_p3_n50.csv};
            \addlegendentry{$N = 50$};
            
            \addplot[color=red!80!black, 
                    style={semithick, densely dotted}]
            table[x=x, y=r, col sep=comma]{./figs/data/shu-osher_p3_n100.csv};
            \addlegendentry{$N = 100$};
            
            \addplot[color=black, 
                    style={thin}]
            table[x=x, y=r, col sep=comma]{./figs/data/shu-osher_p3_n200.csv};
            \addlegendentry{$N = 200$};
    		
            \addplot[color=black!80, style={thin}, only marks, mark=o, mark options={scale=0.5}, mark repeat = 2, mark phase = 0,
            x filter/.expression={abs(x - 6.5) < 1 ? \pgfmathresult : NaN)}
            ]
				table[x=x,y=r,col sep=comma,unbounded coords=jump]{./figs/data/shu-osher_ref.csv};
    
            \addplot[color=black!80, style={thin}, only marks, mark=o, mark options={scale=0.5}, mark repeat = 7, mark phase = 0,
            x filter/.expression={x > 7.5 ? \pgfmathresult : NaN)}
            ]
				table[x=x,y=r,col sep=comma,unbounded coords=jump]{./figs/data/shu-osher_ref.csv};

			\coordinate (spypoint) at (axis cs:6.5,3.85);
            \coordinate (magnifyglass) at (axis cs:2.35,1.65);

		\end{axis} 		
	
        \spy [black] on (spypoint) in node[fill=white] at (magnifyglass);

	\end{tikzpicture}}}
    \newline
    \caption{\label{fig:shuosher} Density profile for the Shu--Osher problem at $t = 1.8$ computed using a $\mathbb P_1$ (left) and $\mathbb P_3$ (right) approximation with varying mesh resolution.}
\end{figure}

To explore the efficacy of the proposed approach at both low and high orders, the problem was simulated with both a $\mathbb P_1$ and $\mathbb P_3$ approximation. A reference solution was computed using a highly-resolved exact Godunov-type solver~\citep{Godunov1959}. The density profiles at $t=1.8$ as computed by a $\mathbb P_1$ and $\mathbb P_3$ approximation are shown in \cref{fig:shuosher} in comparison to the reference simulation. It can be seen that the solution was well-behaved around the shock at both approximation orders and various mesh resolutions, with minimal spurious oscillations. Furthermore, good resolution of the oscillatory density field was obtained, with clear qualitative convergence to the exact solution with increasing resolution. As expected, the $\mathbb P_1$ was better behaved, with fewer spurious oscillations. However, the $\mathbb P_3$ approximation showed better resolution of the oscillatory flow features at lower mesh resolution, indicating that the high-order accuracy of the underlying scheme was sufficiently retained in smooth regions. 

\subsection{Leblanc shock tube}
As a final one-dimensional numerical experiment, the Leblanc shock tube was simulated to validate the approach for extreme shocks and flow conditions. This notoriously difficult test case consists of pressure ratios of $\mathcal O(10^9)$ and density ratios of $\mathcal (10^3)$ and is extremely challenging for numerical schemes to accurately and robustly resolve~\citep{Liu2009}. The problem was solved on the domain $\Omega = [0,9]$ with Dirichlet boundary conditions, and the initial conditions were set as

    \begin{equation}
         \mathbf{q}(x,0) =  \begin{cases}
            \mathbf{q}_l, &\mbox{if } x\leqslant 3, \\
            \mathbf{q}_r, &\mbox{else},
        \end{cases} \quad \mathrm{where} \quad
        \mathbf{q}_l = 
        \begin{bmatrix}
            1\\
            0\\
            (\gamma -1) 10^{-1}
        \end{bmatrix}, \quad
        \mathbf{q}_r =  
        \begin{bmatrix}
            10^{-3}\\
            0\\
            (\gamma -1) 10^{-10}
        \end{bmatrix}.
    \end{equation}
This problem consists of a right-running strong shock and contact discontinuity as well as a left-running expansion wave. Due to the large density jump at the contact discontinuity, which can be smeared by numerical dissipation, it is common for numerical approximations to overpredict the specific internal energy at the contact, which, in turn, overpredicts the shock speed. As such, fine mesh resolution is often required to accurately predict the shock location~\citep{Liu2009}.

   \begin{figure}[tbhp]
        \centering
        \subfloat[Density]{\adjustbox{width=0.45\linewidth, valign=b}{    \begin{tikzpicture}[spy using outlines={rectangle, height=4.5cm,width=2.3cm, magnification=1.75, connect spies}]
		\begin{semilogyaxis}[name=plot1,
		    axis x line=left,
            axis y line=left,
		    xlabel={$x$},
		    xtick={0,1,2,3,4,5,6,7,8,9},
    		xmin=0,
    		xmax=9,
    		x tick label style={
        		/pgf/number format/.cd,
            	fixed,
            	precision=1,
        	    /tikz/.cd},
    		ylabel={$\rho$},
    		ylabel style={rotate=-90},
		    ytick={1e-3, 1e-2, 1e-1, 1e-0},
    		ymin=8e-4,
    		ymax=1.02,
    		legend style={at={(0.03,0.03)},anchor=south west,font=\small, column sep=0.2cm},
    		legend cell align={left},
    		style={font=\normalsize}
      ]

            \addplot [color=black!80, style={thin}, only marks, mark=o, mark options={scale=0.5}, mark repeat = 3, mark phase = 0] table[x=x,y=r,col sep=comma,unbounded coords=jump]{./figs/data/leblanc_exact.csv};


        \addplot[color=blue!80!black, 
                style={semithick, densely dotted}]
				table[x=x,y=r,col sep=comma,unbounded coords=jump]{./figs/data/leblanc_p3_n100.csv};

        \addplot[color=red!80!black, 
                style={thin, dashed}]
				table[x=x,y=r,col sep=comma,unbounded coords=jump]{./figs/data/leblanc_p3_n400.csv};

        \addplot[color=black, 
                style={semithick}]
				table[x=x,y=r,col sep=comma,unbounded coords=jump]{./figs/data/leblanc_p3_n1600.csv};
			\coordinate (spypoint) at (axis cs:6.5,0.015);
            \coordinate (magnifyglass) at (axis cs:1.75,0.015);

		\end{semilogyaxis} 		    
		
        \spy [black] on (spypoint) in node[fill=white] at (magnifyglass);
	\end{tikzpicture}}}
        \subfloat[Velocity]{\adjustbox{width=0.45\linewidth, valign=b}{    \begin{tikzpicture}[spy using outlines={rectangle, height=6cm,width=2.3cm, magnification=3, connect spies}]
		\begin{axis}[name=plot1,
		    axis x line=left,
            axis y line=left,
		    xlabel={$x$},
		    xtick={0,1,2,3,4,5,6,7,8,9},
    		xmin=0,
    		xmax=9,
    		x tick label style={
        		/pgf/number format/.cd,
            	fixed,
            	precision=1,
        	    /tikz/.cd},
    		ylabel={$U$},
    		ylabel style={rotate=-90},
		    ytick={0,0.2,0.4,0.6,0.8},
    		ymin=-0.05,
    		ymax=0.7,
    		y tick label style={
        		/pgf/number format/.cd,
            	fixed,
            	fixed zerofill,
            	precision=1,
        	    /tikz/.cd},
    		legend style={at={(0.03,1)},anchor=north west ,font=\small, column sep=0.2cm},
    		legend cell align={left},
    		style={font=\normalsize}]

            \addplot [color=black!80, style={thin}, only marks, mark=o, mark options={scale=0.5}, mark repeat = 3, mark phase = 0] table[x=x,y=u,col sep=comma,unbounded coords=jump]{./figs/data/leblanc_exact.csv};
        \addlegendentry{Exact};


        \addplot[color=blue!80!black, 
                style={semithick, densely dotted}]
				table[x=x,y=u,col sep=comma,unbounded coords=jump]{./figs/data/leblanc_p3_n100.csv};
        \addlegendentry{$N = 100$};

        \addplot[color=red!80!black, 
                style={thin, dashed}]
				table[x=x,y=u,col sep=comma,unbounded coords=jump]{./figs/data/leblanc_p3_n400.csv};
        \addlegendentry{$N = 400$};

        \addplot[color=black, 
                style={semithick}]
				table[x=x,y=u,col sep=comma,unbounded coords=jump]{./figs/data/leblanc_p3_n1600.csv};
        \addlegendentry{$N = 1600$};

		\end{axis} 		    
		
	\end{tikzpicture}}}
        \newline
        \subfloat[Pressure]{\adjustbox{width=0.45\linewidth, valign=b}{    \begin{tikzpicture}[spy using outlines={rectangle, height=6cm,width=2.3cm, magnification=3, connect spies}]
		\begin{semilogyaxis}[name=plot1,
		    axis x line=left,
            axis y line=left,
		    xlabel={$x$},
		    xtick={0,1,2,3,4,5,6,7,8,9},
    		xmin=0,
    		xmax=9,
    		x tick label style={
        		/pgf/number format/.cd,
            	fixed,
            	precision=1,
        	    /tikz/.cd},
    		ylabel={$P$},
    		ylabel style={rotate=-90},
		    ytick={1e-12, 1e-10, 1e-8, 1e-6, 1e-4, 1e-2, 1e-0},
    		ymin=1e-12,
    		ymax=1,
    		legend style={at={(1,1)},anchor=north east ,font=\small, column sep=0.2cm},
    		legend cell align={left},
    		style={font=\normalsize}]

            \addplot [color=black!80, style={thin}, only marks, mark=o, mark options={scale=0.5}, mark repeat = 3, mark phase = 0] table[x=x,y=p,col sep=comma,unbounded coords=jump]{./figs/data/leblanc_exact.csv};

        \addplot[color=blue!80!black, 
                style={semithick, densely dotted}]
				table[x=x,y=p,col sep=comma,unbounded coords=jump]{./figs/data/leblanc_p3_n100.csv};
    
        \addplot[color=red!80!black, 
                style={thin, dashed}]
				table[x=x,y=p,col sep=comma,unbounded coords=jump]{./figs/data/leblanc_p3_n400.csv};
    
        \addplot[color=black, 
                style={semithick}]
				table[x=x,y=p,col sep=comma,unbounded coords=jump]{./figs/data/leblanc_p3_n1600.csv};

		\end{semilogyaxis} 		    
		
	\end{tikzpicture}}}
        \subfloat[Specific internal energy]{\adjustbox{width=0.45\linewidth, valign=b}{    \begin{tikzpicture}[spy using outlines={rectangle, height=3cm,width=2.3cm, magnification=3, connect spies}]
		\begin{axis}[name=plot1,
		    axis x line=left,
            axis y line=left,
		    xlabel={$x$},
		    xtick={0,1,2,3,4,5,6,7,8,9},
    		xmin=0,
    		xmax=9,
    		x tick label style={
        		/pgf/number format/.cd,
            	fixed,
            	precision=1,
        	    /tikz/.cd},
    		ylabel={$e$},
    		ylabel style={rotate=-90},
		    ytick={0, 0.05, 0.10, 0.15, 0.20},
    		ymin=-0.02,
    		ymax=0.23,
    		y tick label style={
        		/pgf/number format/.cd,
            	fixed,
            	precision=2,
        	    /tikz/.cd},
    		legend style={at={(1,1)},anchor=north east ,font=\small, column sep=0.2cm},
    		legend cell align={left},
    		style={font=\normalsize}]

            \addplot [color=black!80, style={thin}, only marks, mark=o, mark options={scale=0.5}, mark repeat = 3, mark phase = 0] table[x=x,y=e,col sep=comma,unbounded coords=jump]{./figs/data/leblanc_exact.csv};


        \addplot[color=blue!80!black, 
                style={semithick, densely dotted}]
				table[x=x,y=e,col sep=comma,unbounded coords=jump]{./figs/data/leblanc_p3_n100.csv};

        \addplot[color=red!80!black, 
                style={thin, dashed}]
				table[x=x,y=e,col sep=comma,unbounded coords=jump]{./figs/data/leblanc_p3_n400.csv};

        \addplot[color=black, 
                style={semithick}]
				table[x=x,y=e,col sep=comma,unbounded coords=jump]{./figs/data/leblanc_p3_n1600.csv};
        		    		
			\coordinate (spypoint) at (axis cs:7,0.205);
            \coordinate (magnifyglass) at (axis cs:4,0.15);

		\end{axis} 		    
		
        \spy [black] on (spypoint) in node[fill=white] at (magnifyglass);
		
	\end{tikzpicture}}}
        \newline
        \caption{\label{fig:leblanc} 
        Density (top left), velocity (top right), pressure (bottom left), and specific internal energy (bottom right) profiles for the Leblanc shock tube problem at $t = 6$ computed using a $\mathbb P_3$ approximation with varying mesh resolution.
        }
    \end{figure}

The problem was solved using a $\mathbb P_3$ approximation with various mesh resolutions ranging from $N = 100$ to $N = 1600$. The predicted density, velocity, pressure, and specific internal energy fields at $t = 6$ are shown in \cref{fig:leblanc}. At all resolutions, the proposed limiting approach effectively stabilized the high-order solution around  very strong discontinuities. At $N = 100$, the shock speed was noticeably overpredicted, and some minor oscillations were observed around the shock and contact. However, at $N=400$, a notably accurate solution was obtained despite the relatively coarse resolution, with only a minor overprediction of the shock speed and specific internal energy at the contact. At $N=1600$, excellent agreement with the exact solution was obtained. 
    
    \begin{figure}[htbp!] 
    \centering
    \begin{tabular}{|c | c | c | c | c | c | c | c | c | }
    \hline
    $N$ & $\rho$ & Rate & $m$ & Rate & $E$ & Rate & $U$ & Rate \\ 
    \hline
    50 & $3.925 \times 10^{-3}$ & - & $1.685 \times 10^{-3}$ & - & $6.329 \times 10^{-4}$ & - & $4.179 \times 10^{-2}$ & - \\
    100 & $2.215 \times 10^{-3}$ & 0.825 & $9.243 \times 10^{-4}$ & 0.866 & $3.595 \times 10^{-4}$ & 0.816 & $2.996 \times 10^{-2}$ & 0.480 \\
    200 & $1.122 \times 10^{-3}$ & 0.982 & $5.062 \times 10^{-4}$ & 0.869 & $1.880 \times 10^{-4}$ & 0.935 & $1.322 \times 10^{-2}$ & 1.180 \\
    400 & $5.471 \times 10^{-4}$ & 1.036 & $2.269 \times 10^{-4}$ & 1.157 & $8.839 \times 10^{-5}$ & 1.088 & $7.492 \times 10^{-3}$ & 0.819 \\
    800 & $2.788 \times 10^{-4}$ & 0.972 & $1.257 \times 10^{-4}$ & 0.853 & $4.683 \times 10^{-5}$ & 0.916 & $3.456 \times 10^{-3}$ & 1.116 \\
    1600 & $1.441 \times 10^{-4}$ & 0.953 & $6.126 \times 10^{-5}$ & 1.036 & $2.354 \times 10^{-5}$ & 0.992 & $1.896 \times 10^{-3}$ & 0.866 \\
    3200 & $7.047 \times 10^{-5}$ & 1.032 & $3.164 \times 10^{-5}$ & 0.953 & $1.180 \times 10^{-5}$ & 0.996 & $8.490 \times 10^{-4}$ & 1.159 \\
    \hline
    \end{tabular}
    \begin{tabular}{|c | c | c | c | c | c | c |}
    \hline
    $N$ & $P$ & Rate & $e$ & Rate & $s$ & Rate \\ 
    \hline
    50 & $1.755 \times 10^{-4}$ & - & $2.009 \times 10^{-2}$ & - & $6.250 \times 10^{-1}$ & - \\
    100 & $1.020 \times 10^{-4}$ & 0.783 & $1.464 \times 10^{-2}$ & 0.456 & $4.595 \times 10^{-1}$ & 0.444 \\
    200 & $4.714 \times 10^{-5}$ & 1.113 & $6.995 \times 10^{-3}$ & 1.066 & $2.204 \times 10^{-1}$ & 1.060 \\
    400 & $2.514 \times 10^{-5}$ & 0.907 & $3.754 \times 10^{-3}$ & 0.898 & $1.188 \times 10^{-1}$ & 0.891 \\
    800 & $1.186 \times 10^{-5}$ & 1.084 & $1.876 \times 10^{-3}$ & 1.001 & $5.928 \times 10^{-2}$ & 1.003 \\
    1600 & $6.309 \times 10^{-6}$ & 0.911 & $1.021 \times 10^{-3}$ & 0.878 & $3.200 \times 10^{-2}$ & 0.889 \\
    3200 & $3.025 \times 10^{-6}$ & 1.060 & $4.759 \times 10^{-4}$ & 1.101 & $1.493 \times 10^{-2}$ & 1.100 \\
    \hline
    \end{tabular} \hspace{0.19\linewidth}
    \captionof{table}{\label{tab:leblanc_error_relaxed} 
     Convergence of the $L^1$ error of the density, momentum, total energy, velocity, pressure, specific internal energy, and physical entropy for the Leblanc shock tube problem at $t = 6$ computed using a $\mathbb P_3$ approximation \emph{with relaxed bounds}.}
    \end{figure}
    \begin{figure}[htbp!] 
    \centering
    \begin{tabular}{|c | c | c | c | c | c | c | c | c | }
    \hline
    $N$ & $\rho$ & Rate & $m$ & Rate & $E$ & Rate & $U$ & Rate \\ 
    \hline
    50 & $3.892 \times 10^{-3}$ & - & $1.679 \times 10^{-3}$ & - & $6.303 \times 10^{-4}$ & - & $4.179 \times 10^{-2}$ & - \\
    100 & $2.198 \times 10^{-3}$ & 0.824 & $9.214 \times 10^{-4}$ & 0.866 & $3.584 \times 10^{-4}$ & 0.815 & $3.000 \times 10^{-2}$ & 0.478 \\
    200 & $1.140 \times 10^{-3}$ & 0.947 & $5.131 \times 10^{-4}$ & 0.845 & $1.916 \times 10^{-4}$ & 0.903 & $1.324 \times 10^{-2}$ & 1.180 \\
    400 & $6.058 \times 10^{-4}$ & 0.912 & $2.478 \times 10^{-4}$ & 1.050 & $9.849 \times 10^{-5}$ & 0.960 & $7.546 \times 10^{-3}$ & 0.811 \\
    800 & $3.276 \times 10^{-4}$ & 0.887 & $1.423 \times 10^{-4}$ & 0.801 & $5.498 \times 10^{-5}$ & 0.841 & $3.490 \times 10^{-3}$ & 1.112 \\
    1600 & $1.844 \times 10^{-4}$ & 0.829 & $7.491 \times 10^{-5}$ & 0.925 & $3.030 \times 10^{-5}$ & 0.860 & $1.938 \times 10^{-3}$ & 0.849 \\
    3200 & $9.658 \times 10^{-5}$ & 0.933 & $4.054 \times 10^{-5}$ & 0.886 & $1.616 \times 10^{-5}$ & 0.907 & $8.740 \times 10^{-4}$ & 1.149 \\
    \hline
    \end{tabular}
    \begin{tabular}{|c | c | c | c | c | c | c |}
    \hline
    $N$ & $P$ & Rate & $e$ & Rate & $s$ & Rate \\ 
    \hline
    50 & $1.733 \times 10^{-4}$ & - & $2.007 \times 10^{-2}$ & - & $6.245 \times 10^{-1}$ & - \\
    100 & $1.012 \times 10^{-4}$ & 0.776 & $1.463 \times 10^{-2}$ & 0.456 & $4.584 \times 10^{-1}$ & 0.446 \\
    200 & $4.963 \times 10^{-5}$ & 1.028 & $6.987 \times 10^{-3}$ & 1.066 & $2.197 \times 10^{-1}$ & 1.061 \\
    400 & $3.144 \times 10^{-5}$ & 0.659 & $3.781 \times 10^{-3}$ & 0.886 & $1.194 \times 10^{-1}$ & 0.880 \\
    800 & $1.724 \times 10^{-5}$ & 0.867 & $1.872 \times 10^{-3}$ & 1.014 & $5.877 \times 10^{-2}$ & 1.022 \\
    1600 & $1.077 \times 10^{-5}$ & 0.678 & $1.015 \times 10^{-3}$ & 0.883 & $3.154 \times 10^{-2}$ & 0.898 \\
    3200 & $5.858 \times 10^{-6}$ & 0.879 & $4.792 \times 10^{-4}$ & 1.083 & $1.490 \times 10^{-2}$ & 1.082 \\
    \hline
    \end{tabular} \hspace{0.19\linewidth}
    \captionof{table}{\label{tab:leblanc_error_strict} 
     Convergence of the $L^1$ error of the density, momentum, total energy, velocity, pressure, specific internal energy, and physical entropy for the Leblanc shock tube problem at $t = 6$ computed using a $\mathbb P_3$ approximation \emph{without relaxed bounds}.}
\end{figure}

Since the problem has an exact solution, the convergence behavior of the proposed approach was evaluated for cases with extreme discontinuities. The error of various solution quantities was computed, calculated as the point-mean $L^1$ norm on the solution nodes normalized by the domain size, including the conserved variables ($\rho$, $m$, $E$) and primitive variables ($u$, $P$) as well as non-quasiconcave functionals ($e$). Furthermore, the physical entropy $s = P\rho^{-\gamma}$ was also computed. The error and convergence rates for these quantities are shown in \cref{tab:leblanc_error_relaxed} for the approach \emph{with relaxed bounds}. It can be seen that the optimal rate of unity was achieved for all quantities. Furthermore, the error and convergence rates for the approach \emph{without relaxed bounds} are shown in \cref{tab:leblanc_error_strict}. Similarly to the Sod shock tube, marginally better error and convergence rates were observed with the relaxed bounds in comparison to the strict bounds. 

\subsection{Sedov blast wave}
The extension to two-dimensional problems was first performed with the Sedov blast wave problem~\citep{Sedov1993}, which consists of an energy point-source in an ambient gas that drives an outward-running radial explosion. The domain was set as $\Omega = [-1.2, 1.2]^2$ with Dirichlet boundary conditions, and the initial conditions were set to consist of an ambient gas with a constant density $\rho_0 = 1$ and velocity $\mathbf{U} = \mathbf{0}$. We take an identical problem setup to \citet{Maire2009} by setting the pressure to an ambient value of $P_a = 10^{-6}$ everywhere except in the element centered at $x = y = 0$, where an overpressure value is set in terms of the internal energy $e_0$ and cell volume $V_0$ as
\begin{equation}
    P_0 = 4(\gamma-1) \rho_0 e_0/V_0.
\end{equation}
Here, $V_0 = \Delta x \Delta y = 2.4^2/N$ is the element volume and $\epsilon_0 = 0.244816$ is the initial specific internal energy. The factor of 4 in the initial pressure is used to account for the fact that the entire domain is being simulated instead of one quadrant.

\begin{figure}[htbp!]
    \centering
    \hfill
    \subfloat[Contours]{\adjustbox{width=0.35\linewidth,valign=b}{\includegraphics{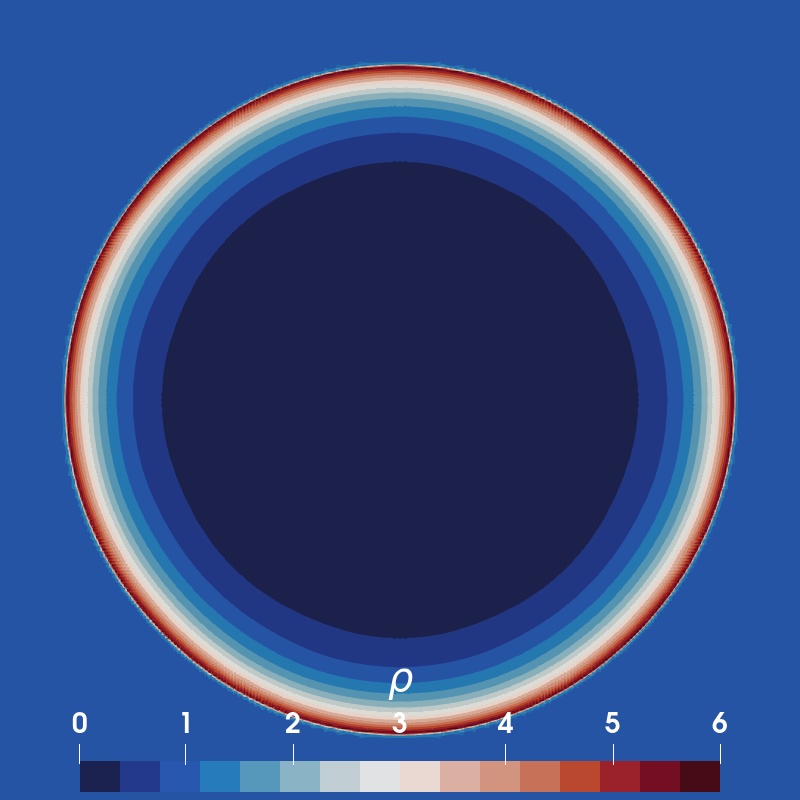}}}
    \hfill
    \subfloat[Profiles]{\adjustbox{width=0.42\linewidth,valign=b}{    \begin{tikzpicture}[spy using outlines={rectangle, height=6cm,width=2.3cm, magnification=3, connect spies}]
		\begin{axis}[name=plot1,
		    axis x line=left,
            axis y line=left,
		    xlabel={$r$},
		    xtick={0,0.2, 0.4, 0.6, 0.8, 1.0, 1.2, 1.4},
    		xmin=0,
    		xmax=1.4,
    		x tick label style={
        		/pgf/number format/.cd,
            	fixed,
            	fixed zerofill,
            	precision=1,
        	    /tikz/.cd},
    		ylabel={$\rho$},
    		ylabel style={rotate=-90},
		    ytick={0,1,2,3,4,5,6},
    		ymin=-0.1,
    		ymax=6,
    		legend style={at={(0.03,0.97)},anchor=north west, font=\small, column sep=0.2cm},
    		legend cell align={left},
    		style={font=\normalsize}
      ]

            \addplot [color=black!80, style={thin}, only marks, mark=o, mark options={scale=0.5}, mark repeat = 1, mark phase = 0] table[x=r,y=d,col sep=comma,unbounded coords=jump]{./figs/data/sedov_exact.csv};
    		\addlegendentry{Exact}


			\addplot[color={black}, style={thin}] table[x=x,y=d,col sep=comma,unbounded coords=jump]{./figs/data/sedov_n257_profiles.csv};
    		\addlegendentry{$N = 257^2$}
      
            \addplot[color=black, style={ultra thin}, only marks, mark=o, mark options={scale=0.5}, samples at={1, 1.02, ..., 1.4}]{1};

		\end{axis} 		    
		
	\end{tikzpicture}}}
    \hfill
    \newline
    \caption{\label{fig:sedov} 
    Contours of density (left) and azimuthally-averaged density profile for the two-dimensional Sedov blast problem at $t = 1$ computed using a $\mathbb P_2$ approximation and $N = 257^2$ elements. }
\end{figure}

The problem was solved using a $\mathbb P_2$ approximation with $N = 257^2$ quadrilateral elements. For this particular case, the Riemann solver was changed to HLL~\citep{Harten1983b} as the HLLC Riemann solver is known to be ill-behaved for this problem. The contours of density as well as the azimuthally-averaged density profile at $t=1$ are shown in \cref{fig:sedov}. We observe the canonical blast wave structure in the solution, with a low density center region and a sharp outward-running shock front. Furthermore, good agreement with the exact solution was obtained, with sharp resolution of the shock and minimal spurious oscillations.

\subsection{Double Mach reflection}
The double mach reflection problem of \citet{Woodward1984} was subsequently used to evaluate the proposed approach in the context of strong shocks with small-scale flow features. The problem consists of a Mach 10 shock impinging on a $30$ degree ramp, which results in multiple strong shock interactions as well as the onset of small-scale Kelvin--Helmholtz instabilities. The problem was solved on the domain $\Omega = [0,4] \times [0,1]$, and the initial conditions were set as 

    \begin{equation*}
         \mathbf{q}(\mathbf{x},0) =  \begin{cases}
            \mathbf{q}_l, &\mbox{if } x < 1/6 + \tan (30^\circ)y,\\
            \mathbf{q}_r, &\mbox{else},
        \end{cases} \quad \mathrm{where} \quad
        \mathbf{q}_l = \begin{bmatrix}
            8 \\ 7.14471\\ -4.125 \\ 116.5
        \end{bmatrix}, \quad
        \mathbf{q}_r = \begin{bmatrix}
            1.4 \\ 0\\ 0 \\ 1
        \end{bmatrix}.
    \end{equation*}
Dirichlet boundary conditions were primarily used, with the post-shock state $\mathbf{q}_l$ at the left boundary and the pre-shock state $\mathbf{q}_r$ for the right boundary. At the bottom boundary, the post-shock state $\mathbf{q}_l$ was set for $x < 1/6$ and a slip adiabatic wall boundary condition was applied for $x \geqslant \frac{1}{6}$. For the top boundary, the exact solution was set, given as
    \begin{equation*}
        \mathbf{q}(\mathbf{x}, t)|_{y=1} = \begin{cases}
                \mathbf{q}_l, &\mbox{if } x\leqslant 1/6  + \tan (30^\circ)y + \frac{10}{\cos(30^\circ)}t,\\
                \mathbf{q}_r, &\mbox{else}.
            \end{cases}
    \end{equation*}

\begin{figure}[htbp!]
    \centering
    \subfloat[$\mathbb P_1$]{{\adjustbox{width=0.8\linewidth,valign=b}{\includegraphics{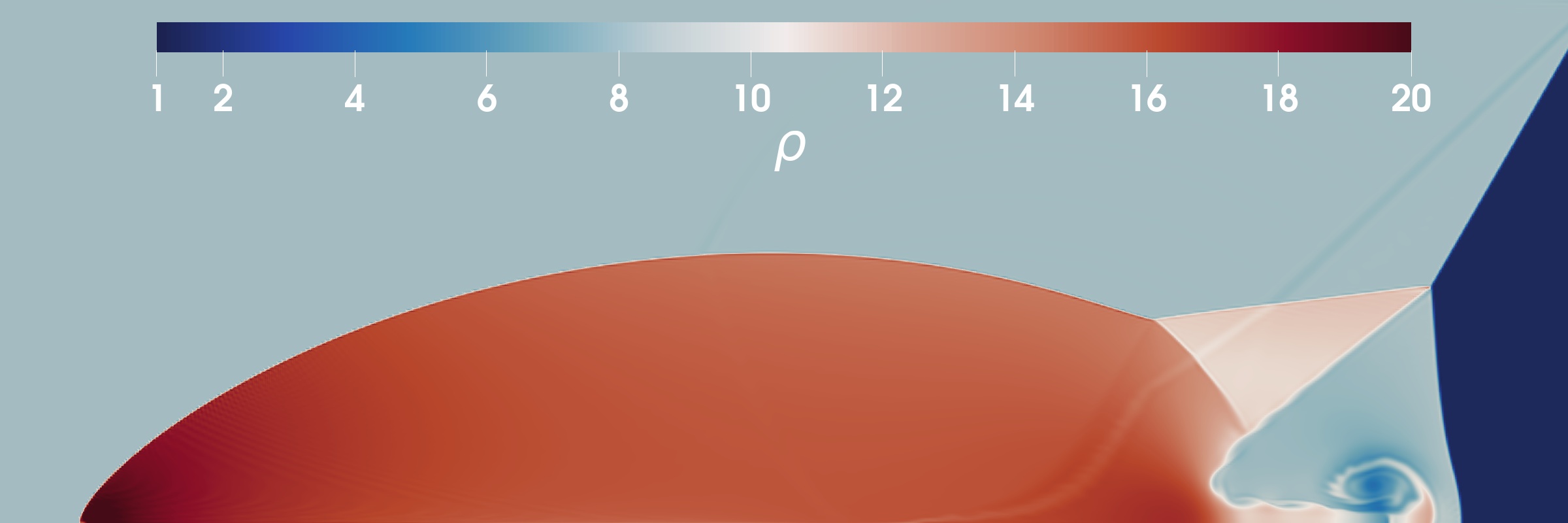}}}}
    \newline
    \subfloat[$\mathbb P_2$]{{\adjustbox{width=0.8\linewidth,valign=b}{\includegraphics{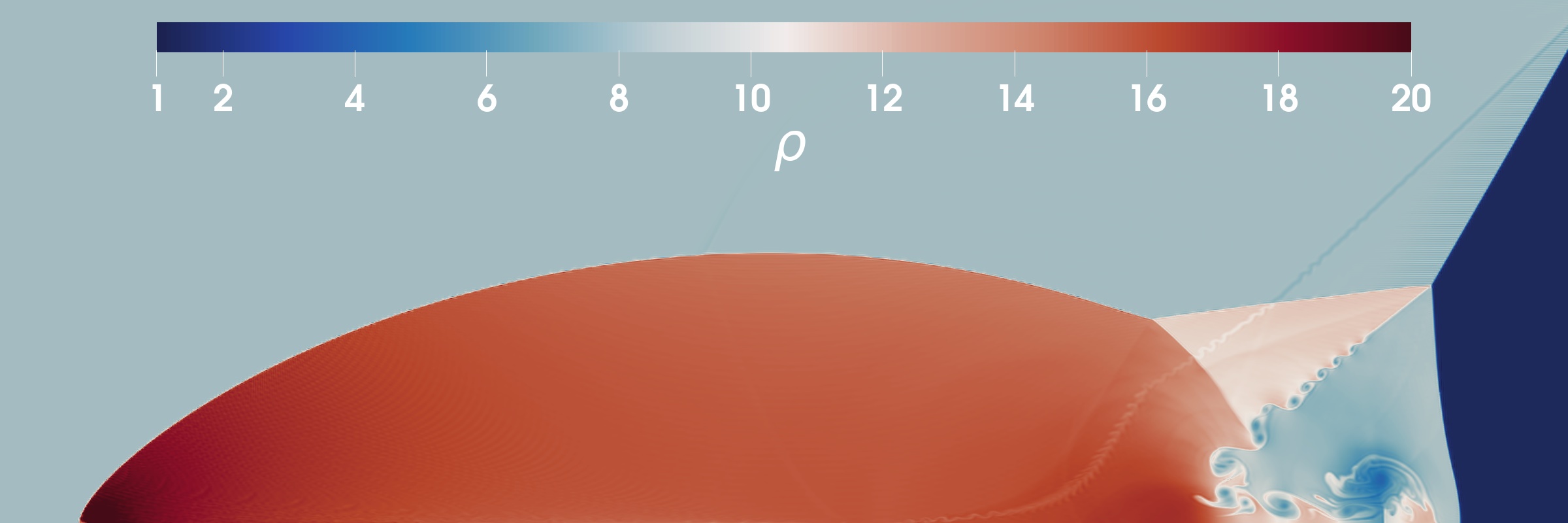}}}}
    \newline
    \subfloat[$\mathbb P_3$]{{\adjustbox{width=0.8\linewidth,valign=b}{\includegraphics{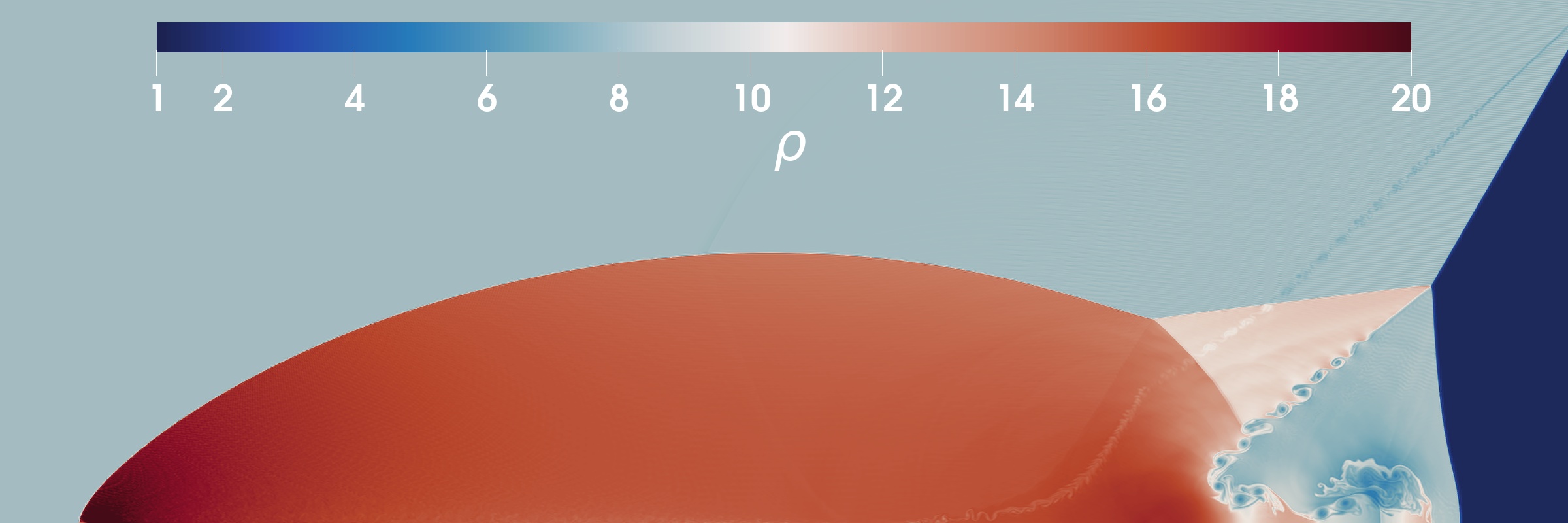}}}}
    \newline
    \caption{\label{fig:doublemach} Contours of density for the Double Mach reflection problem at $t = 0.2$ on the subregion $0 \leq x \leq 3$ computed using a $\mathbb P_1$ to $\mathbb P_3$ approximation with $N = 960 \times 240$ elements. }
\end{figure}

The problem was computed using a $\mathbb P_1$ to $\mathbb P_3$ approximation with $N = 960 \times 240$ quadrilateral elements. The contours of density at $t = 0.2$ are shown in \cref{fig:doublemach}. Excellent resolution of the strong shocks in the flow was observed at all polynomial orders, with effectively sub-element resolution of discontinuities without noticeable spurious oscillations. Furthermore, the onset of the Kelvin--Helmholtz instabilities rolling up off the contact line was well-resolved by the approach, indicating that the limiter does not erroneously dissipate small-scale flow features. It can also be seen that with increasing polynomial order, better resolution of these small-scale flow features was obtained, highlighting the benefits of the approach at higher orders.

\subsection{Richtmyer--Meshkov instability}
\begin{figure}[htbp!]
    \centering
    \subfloat[$N = 10 \times 80$]{{\adjustbox{width=0.195\linewidth,valign=b}{\includegraphics{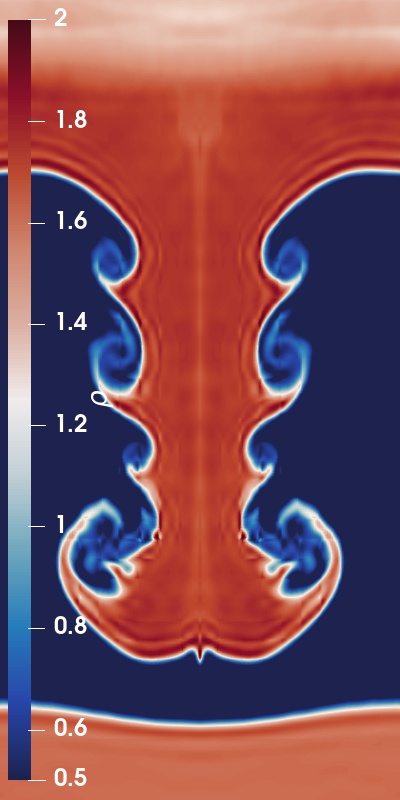}}}}
    \hfill
    \subfloat[$N = 20 \times 160$]{{\adjustbox{width=0.195\linewidth,valign=b}{\includegraphics{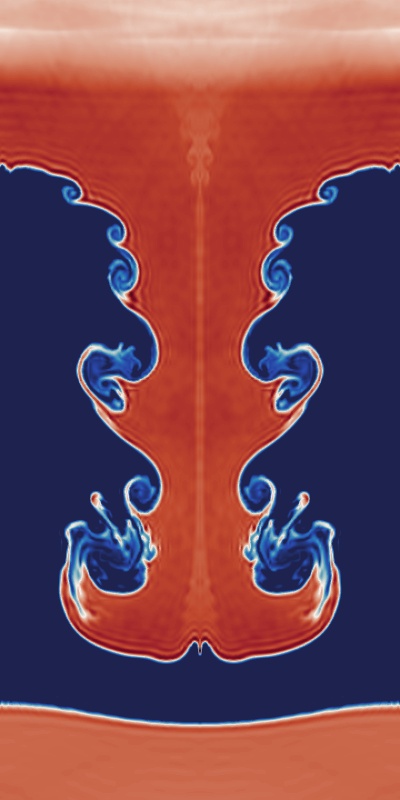}}}}
    \hfill
    \subfloat[$N = 40 \times 320$]{{\adjustbox{width=0.195\linewidth,valign=b}{\includegraphics{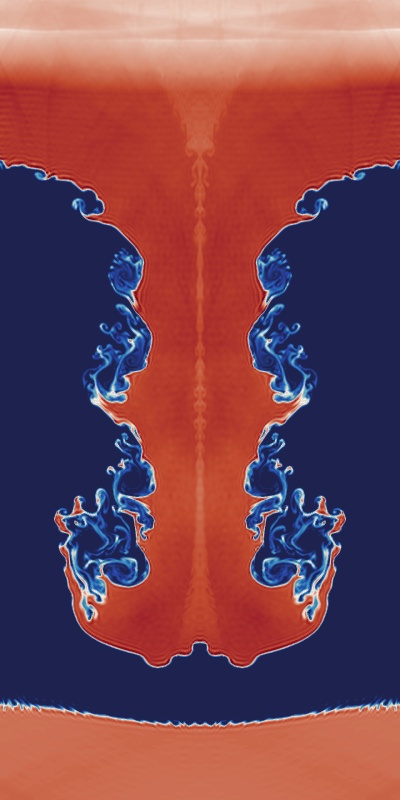}}}}
    \hfill
    \subfloat[$N = 80 \times 640$]{{\adjustbox{width=0.195\linewidth,valign=b}{\includegraphics{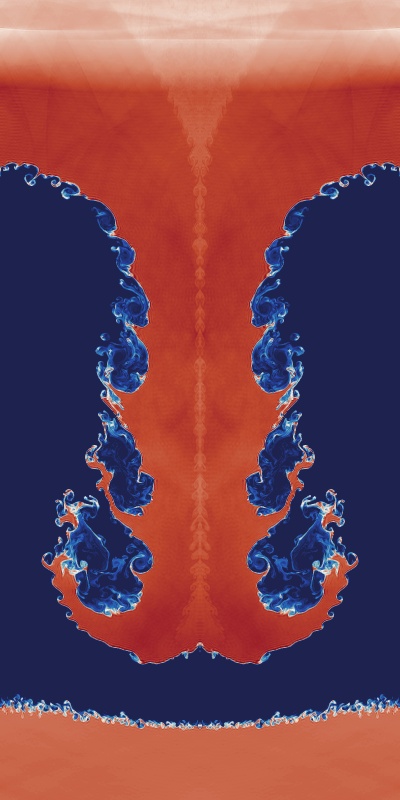}}}}
    \hfill
    \subfloat[$N = 160 \times 1280$]{{\adjustbox{width=0.195\linewidth,valign=b}{\includegraphics{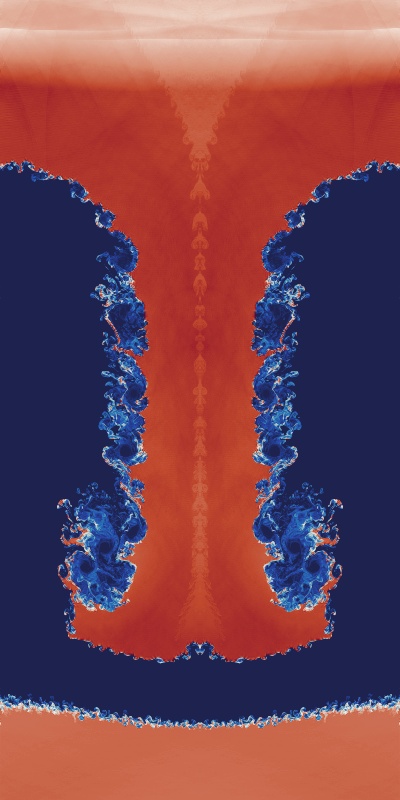}}}}
    \caption{\label{fig:rmi}  Contours of density for the Richtmyer--Meshkov instability problem at $t = 2$ on the subregion $2 \leq y \leq 4$ computed using a $\mathbb P_4$ approximation with varying mesh resolution.  Contours are reflected about the $y$-axis for visualization.}
\end{figure}

To further evaluate the ability of the approach in resolving small-scale vortical flow features in the vicinity of discontinuities, the Richtmyer--Meshkov instability was simulated. This instability, predicted analytically by \citet{Richtmyer1960} and shown experimentally by \citet{Meshkov1972}, stems from an impulse acting upon a contact discontinuity, generally as a result of a propagating shock wave. The problem was solved on the domain $\Omega = [0, 0.5] \times [0, 4]$, and the initial conditions were set as
    \begin{equation*}
         \mathbf{q}(\mathbf{x},0) =  \begin{cases}
            \mathbf{q}_t, &\mbox{if } y \geq 2 - \frac{1}{2}\cos (2 \pi x),\\
            \mathbf{q}_m, &\mbox{if } \frac{1}{2} \leq y < 2 - \frac{1}{2}\cos (2 \pi x) ,\\
            \mathbf{q}_b, &\mbox{else},
        \end{cases} \quad \mathrm{where} \quad
        \mathbf{q}_b = \begin{bmatrix}
            0.25 \\ 0\\ 0 \\ 1
        \end{bmatrix}, \quad
        \mathbf{q}_m = \begin{bmatrix}
            1 \\ 0\\ 0 \\ 1
        \end{bmatrix}, \quad
        \mathbf{q}_t = \begin{bmatrix}
            8/3 \\ 0\\ 0 \\ 4.5
        \end{bmatrix}.
    \end{equation*}
These conditions yield an Atwood number of 0.6. Slip adiabatic wall boundary conditions were enforced along the left and right boundaries (mimicking symmetry) and Dirichlet boundary conditions were enforced along the top and bottom boundaries. 

The problem was computed using a $\mathbb P_4$ approximation with progressively finer mesh resolution, starting with $N = 10 \times 80$ quadrilateral elements and increasing to $N = 160 \times 1280$ quadrilateral elements. The density contours at $t = 2$ as computed on the various meshes are shown in \cref{fig:rmi}. Good resolution of the discontinuities in the flow was observed even with the higher-order approximation, and the rollup of the vortices was effectively captured by the approach with increasingly smaller-scale vortical structures appearing with higher mesh resolution. These results indicate that the low numerical dissipation and scale-resolving capabilities of the underlying high-order DG scheme are not noticeably degraded with the proposed limiter. 

\subsection{Obstructed radial explosion}

\begin{figure}[htbp!]
    \centering
    \subfloat[$t = 1$]{{\adjustbox{width=0.33\linewidth,valign=b}{\includegraphics{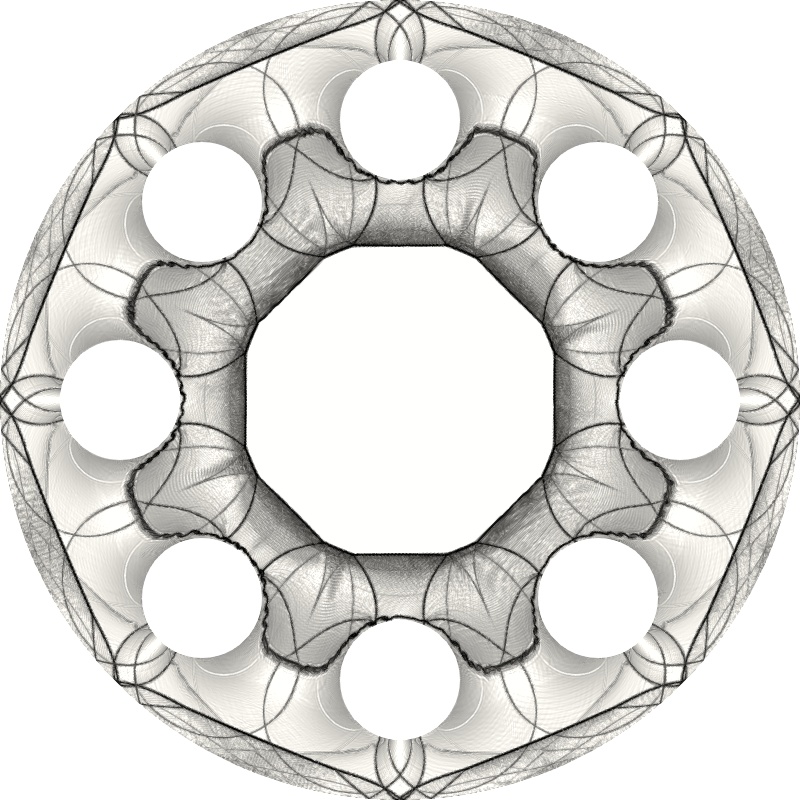}}}}
    \hspace{25pt}
    \subfloat[$t = 2$]{{\adjustbox{width=0.33\linewidth,valign=b}{\includegraphics{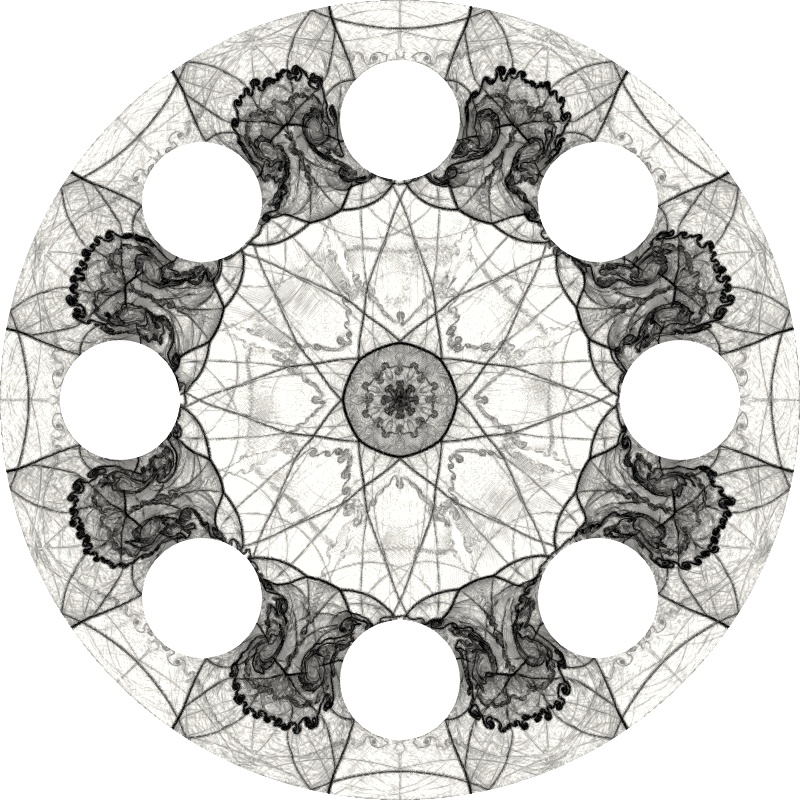}}}}
    \newline
    \subfloat[$t = 3$]{{\adjustbox{width=0.33\linewidth,valign=b}{\includegraphics{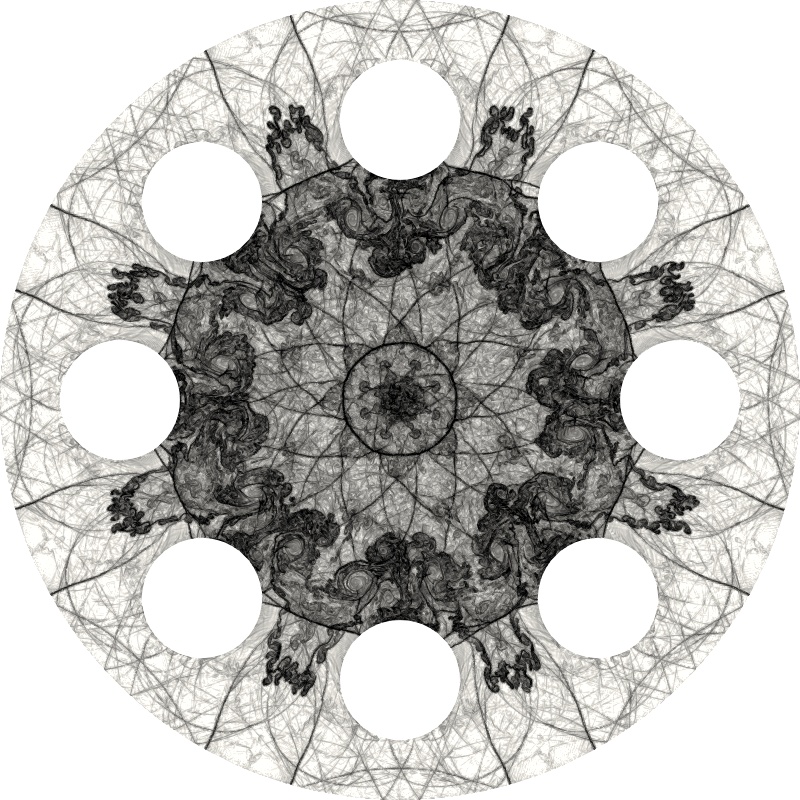}}}}
    \hspace{25pt}
    \subfloat[$t = 4$]{{\adjustbox{width=0.33\linewidth,valign=b}{\includegraphics{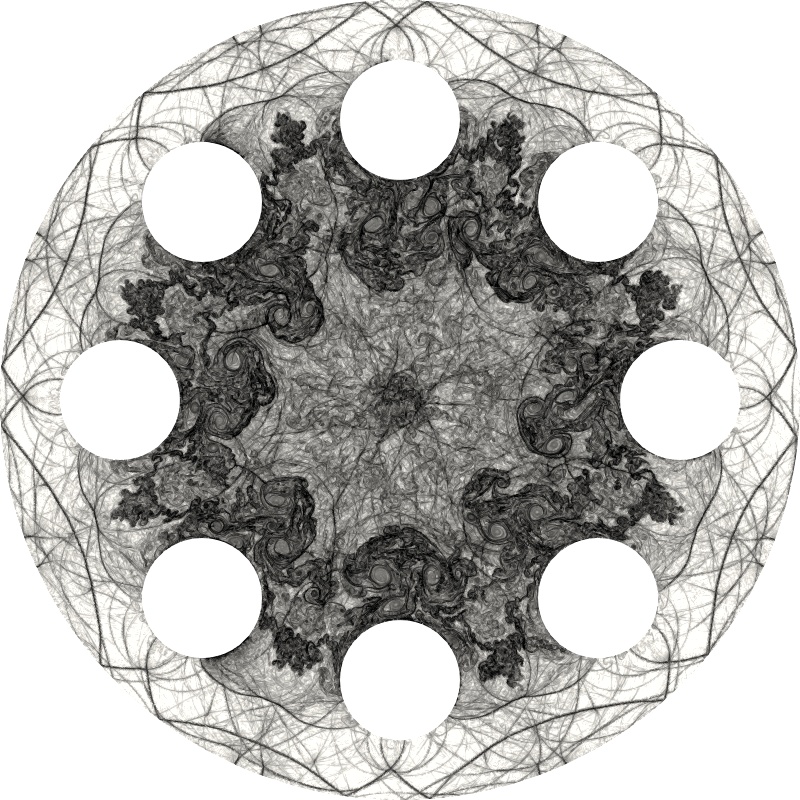}}}}
    \newline
    \caption{\label{fig:exp} Schlieren-type representation of the density gradient norm for the obstructed radial explosion problem at varying times computed using a $\mathbb P_2$ approximation with $9.68 \times 10^5$ unstructured/curved triangular elements.}
\end{figure}

The explosion problem of \citet{Nazarov2017} was then used to showcase the extension of the approach to unstructured/curved grids. The problem consists of a circular domain of radius 2, with eight circular cylinders of radius 0.3 placed azimuthally equispaced such that their centers were located 1.4 units from the origin. An initial overpressure region in the center drives an outward-running shock and contact discontinuity, which interact with the curved surfaces to form complex shock interactions. The initial conditions were set similarly to a radially-symmetric Sod shock tube problem, with an inner gas region of radius $\sqrt{0.4}$ (centered at the origin) set to $\mathbf{q} = [1, 0, 0, 1]$ and the remaining domain set to $\mathbf{q} = [0.125, 0, 0, 0.1]$. Adiabatic slip wall boundary conditions were applied to all surfaces.

An unstructured mesh was generated with approximately $9.68 \times 10^5$ second-order triangular elements. The problem was then simulated on this mesh using a $\mathbb P_2$ approximation. A Schlieren-type representation of the density gradient norm is shown in \cref{fig:exp} at several time instances. The interactions of the shocks with the contact line, which form Richtmyer--Meshkov instabilities and complex vortical flow, were well captured by the approach. Good resolution of the discontinuities in the flow as well as the small-scale flow features was observed, showcasing the applicability of the approach on unstructured meshes.

\subsection{Mach 800 astrophysical jet}
As a final validation study and an evaluation of the approach for very extreme flow conditions, the problem of high-speed astrophysical jets was considered. The test case was introduced  by \citet{Balsara2012} and consists of a Mach 800 jet expanding into an ambient gas. An identical problem setup to \citet{Dzanic2022} was used. The half-domain was set as $\Omega = [0, 0.5] \times [0, 1.5]$, with adiabatic slip wall (symmetry) boundary conditions along the $y$-axis and Dirichlet boundary conditions on the outer boundaries. The initial conditions were set to $\mathbf{q} = [0.1\gamma, 0, 0, 1]$. At the $y=0$ boundary, the inlet region was defined on $x \leqslant 0.05$, with the solution set as $\mathbf{q} = [\gamma, 0, 800, 1]^T$, while the remaining boundary region was set to a Dirichlet boundary conditions.

\begin{figure}[htbp!]
    \centering
    \subfloat[$N = 200 \times 600$]{{\adjustbox{width=0.498\linewidth,valign=b}{\includegraphics{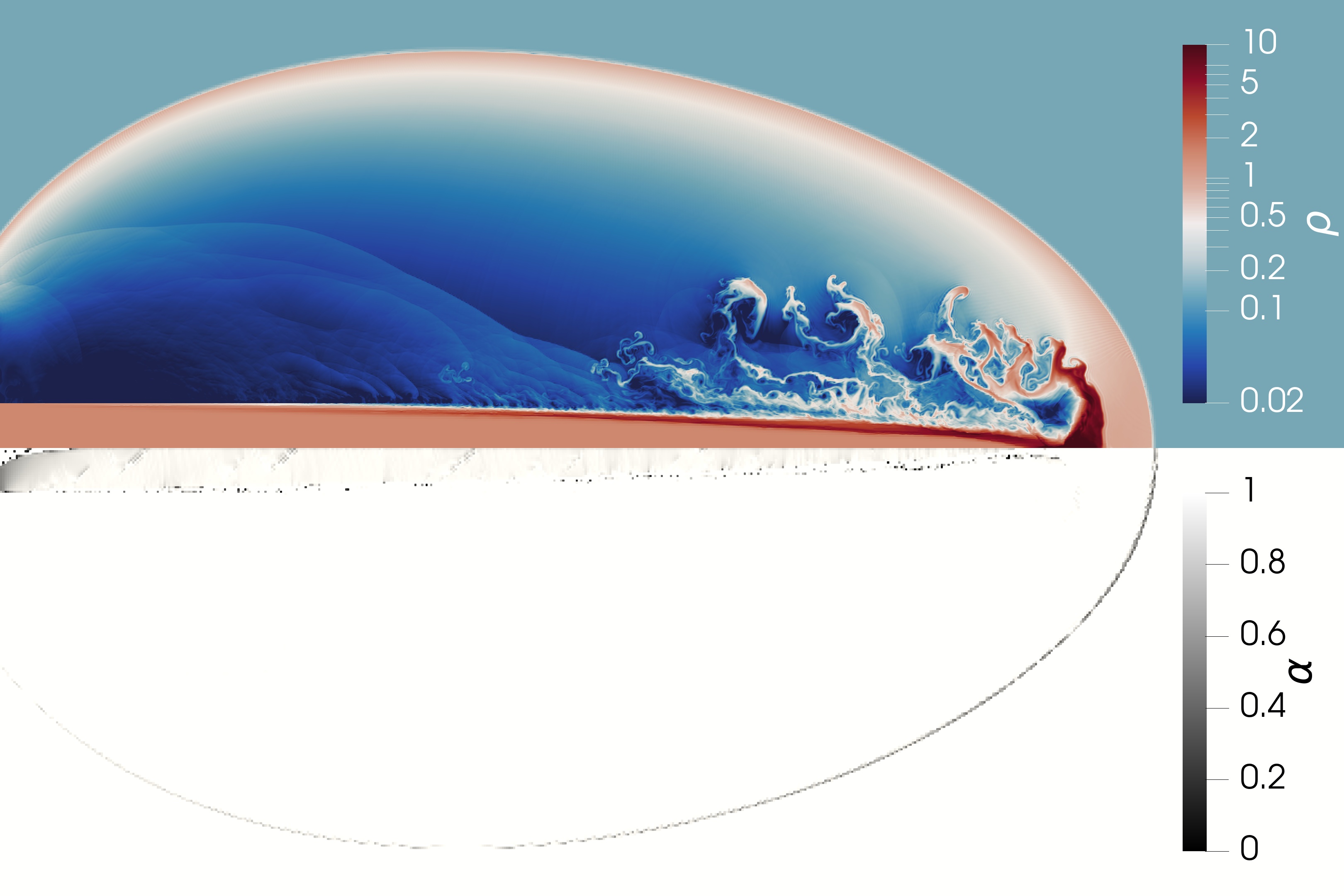}}}}
    \hfill
    \subfloat[$N = 800 \times 2400$]{{\adjustbox{width=0.498\linewidth,valign=b}{\includegraphics{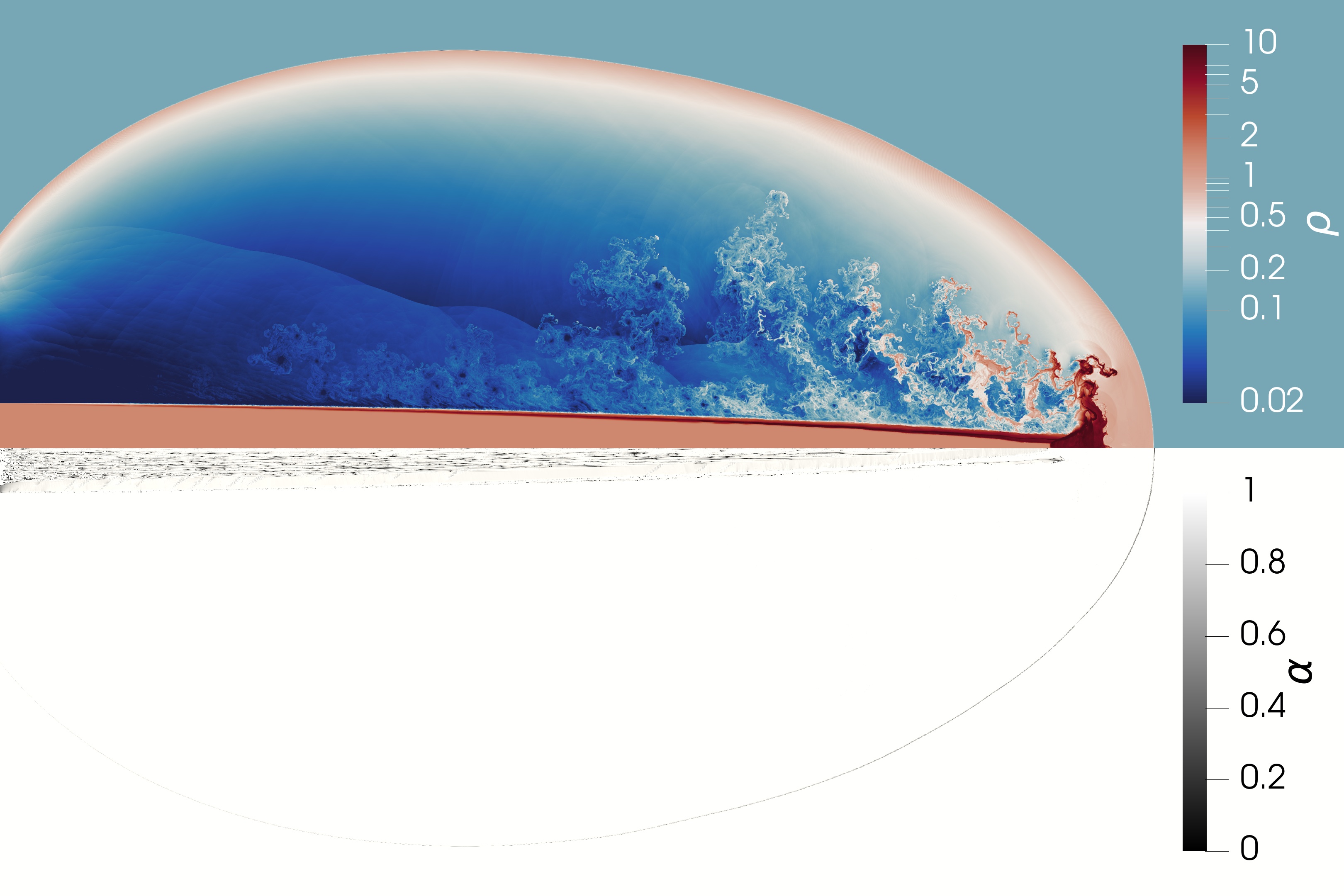}}}}
    \caption{\label{fig:jet} Contours of density (top) and limiting factor distribution (bottom) for the Mach 800 astrophysical jet problem at $t = 0.002$ computed using a $\mathbb P_3$ approximation with $N = 200 \times 600$ (left) and $N = 800 \times 2400$ (right) elements. View is reoriented such that the $+y$ direction is shown left-to-right.}
\end{figure}

The problem was solved using a $\mathbb P_3$ approximation on a coarse and fine mesh, consisting of $N = 200 \times 600$ and $N = 800 \times 2400$  quadrateral elements, respectively. The contours of density and the limiting factor distribution at $t = 0.002$ as computed on these meshes are shown in \cref{fig:jet}. The strong leading shock wave was robustly resolved without noticeable spurious oscillations, and small-scale flow features around the jet showed excellent resolution. With increasing mesh resolution, smaller-scale flow features were observed around the jet. Given the identical case setup, the results can be compared to the entropy filtering approach in \citet{Dzanic2022}. It can be seen that noticeably better resolution of the small-scale flow features was obtained with the proposed approach with similar resolution of the shock structure. Furthermore, the distribution of the limiting factor shows that the limiting is primarily applied along the shock front. No visually discernable limiting was applied within the jet cocoon in the vicinity of small-scale flow features, except within the jet interior where the near-constant solution can spuriously activate the limiter.
\section{Conclusions}\label{sec:conclusion}
In this work, we introduced a novel limiting approach for high-order approximations of the Euler equations. Motivated by the observation that high-order approximations of the Boltzmann equation remain robust around discontinuities and that linear particle transport is notably easier to design limiters for, the proposed approach connects the evolution of the macroscopic solution of the Euler equations to the underlying kinetic representation provided by the Boltzmann equation, which allows for extending limiters designed for linear advection directly to the Euler equations. The technique takes an arbitrary stencil of local solution values and, by applying the corresponding local maximum principle to the analytic expression of linear transport of the Boltzmann particle distribution function in the inviscid limit, yields an expression for the admissible region of particle distribution function values, which may be numerically integrated to yield a set of bounds for the density, momentum, and total energy. Several key properties of these bounds are shown, including preservation of the positivity of density/pressure/internal energy. This approach was applied to an explicit unstructured high-order discontinuous Galerkin/flux reconstruction scheme for a variety of problems in gas dynamics such as extreme shocks, smooth flows, and shock-vortex interactions. It was observed that the approach could robustly resolve strong discontinuities with essentially optimal convergence rates while recovering high-order accuracy in smooth regions without any ad hoc corrections. 

This work presents potential for various extensions and improvements which are in the process of being developed. The approach is not strictly restricted to the limiting method used in this work as it primarily is used to compute solution bounds given a set of discrete solution values. Therefore, it may directly be applied to flux-corrected transport-type techniques which rely on enforcing bounds over a set of Riemann-averaged states~\citep{Guermond2018} since it was shown that the Riemann-averaged states preserve these bounds. Furthermore, the kinetic representation provided by the Boltzmann equation underpins many complex macroscopic governing equations, some for which admissible solution bounds are not necessarily well-understood. As such, the extension of this method to these systems, namely the Navier--Stokes equations and multi-physics systems such as magnetohydrodynamics, can result in more physically consistent limiting techniques. Additionally, algorithmic improvements can be made to increase performance, such as applying adaptive quadrature-based integration methods for the distribution functions to mitigate the overall cost of performing quadrature in the velocity domain, which currently can be several times the cost of the solver time step when extending to 3D problems.

\section*{Acknowledgements}
\label{sec:ack}
This work was partially performed under the auspices of the U.S. Department of Energy by Lawrence Livermore National Laboratory under contract DE--AC52--07NA27344 and the LLNL-LDRD Program under Project tracking No.\ 24--ERD--050. Release number LLNL--JRNL--867251--DRAFT.
\bibliographystyle{unsrtnat}
\bibliography{reference}



\end{document}